\documentclass[11pt,a4paper,twoside]{article}
\usepackage[hmarginratio=1:1,bmargin=1.3in]{geometry}
\usepackage[OT2,T1]{fontenc}
\usepackage{amsmath, amsthm}
\allowdisplaybreaks
\usepackage{amsfonts}
\usepackage{amssymb}
\usepackage{mathrsfs}
\usepackage[all]{xy}
\usepackage{graphicx}
\usepackage[nottoc, notlof, notlot]{tocbibind}
\usepackage{array}
\usepackage{url}
\usepackage{rotating}
\usepackage{fancyhdr}
\usepackage{titlesec}
\usepackage{xfrac} 
\usepackage{multicol}
\usepackage{comment}
\usepackage{cases}
\usepackage{hyperref}
\usepackage[normalem]{ulem}
\usepackage{enumerate}

\title{On the cohomological dimension of Henselian valued fields}
\author{Diego Izquierdo, Giancarlo Lucchini Arteche and Gonzalo Manzano Flores}
\date{}

\titleformat{\section}[hang]{\center\Large\bf}{\thesection.}{0.5cm}{}

\DeclareSymbolFont{cyrletters}{OT2}{wncyr}{m}{n}
\DeclareMathSymbol{\Sha}{\mathalpha}{cyrletters}{"58}
\DeclareMathSymbol{\Brusse}{\mathalpha}{cyrletters}{"42}

\theoremstyle{plain}
\newtheorem{theorem}{Theorem}[section]

\newtheorem{lemma}[theorem]{Lemma}
\newtheorem{proposition}[theorem]{Proposition}
\newtheorem{corollary}[theorem]{Corollary}
\newtheorem{definition}[theorem]{Definition}

\theoremstyle{definition}
\newtheorem{remarque}[theorem]{Remark}

\newtheorem{notation}[theorem]{Notation}

\renewcommand{\cal}[1]{\mathcal{#1}}

\def\Gal{\mathop{\rm Gal}\nolimits}

\newcommand{\K}{\mathrm{K}^\mathrm{M}}
\newcommand{\kk}{\mathrm{k}^\mathrm{M}}

\newcommand \Z {{\mathbb Z}}
\newcommand \Q {{\mathbb Q}}
\newcommand \F {{\mathbb F}}
\newcommand \N {{\mathbb N}}
\newcommand \R {{\mathbb R}}

\renewcommand{\ker}{\mathrm{Ker}}

\newcommand \im {{\mathrm{Im}}}

\newcommand{\unr}{\mathrm{unr}}

\usepackage{marvosym}

\usepackage[dvipsnames,svgnames,table]{xcolor}

\numberwithin{equation}{section}

\begin{document}

\maketitle

\begin{abstract}
Let $(K,v)$ be a Henselian valued field of rank one and of mixed characteristic $(0,p)$, with value group $\Gamma$ dense in $\R$ and residue field $k$. We compute all the graded pieces of the filtration induced by $v$ on the Milnor $\mathrm{K}$-theory modulo $p$ of $K$, in terms of the exterior powers of the $\F_p$-vector space $\Gamma/p\Gamma$ and of the modules of differential forms of $k$. This extends to the non-discrete setting a classical theorem of Kato concerning complete discretely valued fields. We also prove that, when $(K,v)$ is maximally complete, no class escapes this filtration by accumulating at a limit. As an application, we prove the inequality $\dim_p(K)\geq\dim_{\F_p}(\Gamma/p\Gamma)+\dim_p(k)$ between Kato's $p$-dimensions, and we show that it becomes an equality as soon as $(K,v)$ is maximally complete.\\

\textbf{MSC Classes:} 19F15, 19D45, 12J10, 12G10, 12G05.

\end{abstract}

\section{Introduction}

Let $(K,v)$ be a Henselian valued field with residue field $k$ and value group $\Gamma$. It is a guiding principle in the arithmetic of valued fields that the invariants of $K$ ought to be computable in terms of those of $k$ and of $\Gamma$. For the $p$-dimension $\dim_p$ of Kato and Gille, recalled in Definition \ref{def cohdim} below, which is the right substitute for the cohomological dimension when the residue field is imperfect, the expected formula reads
\begin{equation}\label{eqn intro formula}
\dim_p(K)=\dim_{\F_p}(\Gamma/p\Gamma)+\dim_p(k).
\end{equation}
When $p$ is invertible in $k$, this is a classical and essentially formal consequence of the structure of the tamely ramified part of the absolute Galois group of $K$, see for instance \cite[Prop.~1.2]{Soman}. The case where $p$ is the residue characteristic is of a completely different nature, since wild ramification is not governed by $\Gamma$ and $k$ in any naive way.

The first general result in this direction is due to Kato \cite{Kato}, who proved \eqref{eqn intro formula} for complete discretely valued fields of mixed characteristic $(0,p)$; there $\Gamma=\Z$, so that $\dim_{\F_p}(\Gamma/p\Gamma)=1$ and the formula reads $\dim_p(K)=1+\dim_p(k)$. Rather than computing Galois cohomology directly, Kato works with the Milnor $\mathrm{K}$-theory modulo $p$ of $K$, which he endows with the filtration $(u_q^\gamma(K))_{\gamma\geq0}$ induced by the valuation, $u_q^\gamma(K)$ being generated by the symbols $\{1+\pi^\gamma x,y_2,\dots,y_q\}$ with $x\in \mathcal{O}_K$ and $y_2,\dots,y_q\in \mathcal{O}_K^\times$. Assuming that $K$ contains all $p$-th roots of unity, he computes all the graded pieces of this filtration in terms of the modules of differential forms of $k$, and shows that the filtration stops at $e:=\frac{p\,v(p)}{p-1}$; formula \eqref{eqn intro formula} then follows by inspection. A detailed and very readable account of these arguments was given by Colliot-Th\'el\`ene in his Trieste lecture notes \cite{CT}, and our exposition owes a lot to both texts. The passage from Milnor $\mathrm{K}$-theory to Galois cohomology, only available in low degrees at the time, is now provided in full by the Bloch--Kato conjecture (due to Rost and Voevodskii, among others), which we use freely.

The purpose of this article is to prove \eqref{eqn intro formula} for the remaining rank one fields, that is, for Henselian valued fields whose value group is \emph{dense} in $\R$. Note that in this generality the $p$-rank $\dim_{\F_p}(\Gamma/p\Gamma)$ may be any element of $\Z_{\geq 0}\cup\{\infty\}$ --- including $0$, when $\Gamma$ is $p$-divisible --- and that $k$ is typically imperfect.

\paragraph{Motivation}
Rank one valued fields with a non-discrete value group are far from being exotic objects, and they show up in several currently very active areas of research. They are, first of all, the fields underlying non-archimedean analytic geometry in the sense of Tate and Berkovich, and in mixed characteristic they include all perfectoid fields, whose value group is $p$-divisible and whose residue field is perfect \cite{FW,Scholze}; observe that for those fields the right-hand side of \eqref{eqn intro formula} reduces to $\dim_p(k)$. They also appear prominently in resolution of singularities and, more specifically, in local uniformization, where one has to deal with arbitrary valuations on function fields and where the defect --- the very same phenomenon that makes the maximality assumption in Theorem \ref{thm cd} below relevant --- is the main obstruction \cite{KnafKuhlmann,KuhlmannDefect}. They are moreover the natural base fields for the recent extensions, due to Mehmeti \cite{Mehmeti}, of the patching techniques of Harbater, Hartmann and Krashen \cite{HHK} to Berkovich curves over arbitrary complete non-archimedean fields. Finally, they are ubiquitous in the model theory of valued fields, where Ax--Kochen--Ershov principles are available well beyond the discretely valued case, for instance for Kuhlmann's tame fields \cite{KuhlmannTame}. Let us also stress that maximally complete fields, for which our main theorem gives a complete answer, are plentiful: every valued field admits a maximal immediate extension, with the same value group and the same residue field, and explicit constructions in mixed characteristic go back to Kaplansky \cite{Kaplansky} and Poonen \cite{Poonen}. We hope that the computations below may be of use to some of these communities. As an illustration of the kind of applications one may expect, computations of this type over complete discretely valued fields are a basic input in the study of homogeneous spaces and of central simple algebras over such fields, see \cite{GIL}.

\paragraph{Results}
Our first main result, Theorem \ref{thm isom Kato}, is the exact analogue of Kato's computation: we exhibit explicit isomorphisms $\rho_q^\gamma$ describing every graded piece of the filtration $(u_q^\gamma(K))_{\gamma\in\Gamma_{\geq0}}$ of $\kk_q(K)=\K_q(K)/p$, the answer being given by exterior powers of $\Gamma/p\Gamma$ tensored with $\kk_\bullet(k)$, with $\Omega_k^\bullet$, with $\Omega_k^\bullet/\ker(d)$ or with $H_p^\bullet(k)$, according to the position of $\gamma$ with respect to $e$ and to $p\Gamma$. Our second main result, Theorem \ref{th invisible terms}, is of a different nature and has no interesting counterpart in the discrete case: it asserts that, over a maximally complete field, every nontrivial class in $\kk_q(K)$ is detected by some graded piece of the filtration. Together, these two statements yield our third main result, Theorem \ref{thm cd}: the inequality $\geq$ in formula \eqref{eqn intro formula} holds for every Henselian field as above, and equality holds as soon as $(K,v)$ is maximally complete, or equivalently maximal in the sense of Kaplansky. Some hypothesis of this kind cannot be dispensed with: the maximal cyclotomic extension of $\Q_p$ satisfies $\dim_{\F_p}(\Gamma/p\Gamma)=\dim_p(k)=0$ and yet has $p$-dimension $1$, as explained in Remark \ref{rem cyclotomic}.

\paragraph{On the proofs}
We follow the general strategy of \cite{Kato} and \cite{CT}, but the density of $\Gamma$ in $\mathbb{R}$ forces us to depart from it at several places, and in fact allows us to simplify some of the arguments.

The first and main novelty concerns the injectivity of the maps $\rho_q^\gamma$ for $0<\gamma<e$. As in the discrete case, it is obtained by a duality argument: one pairs a given class of $u_q^{\gamma}(K)/u_q^{>\gamma}(K)$ with a well-chosen class of a complementary graded piece and one checks that the resulting cup-product does not vanish. In Kato's setting, the only possible target for such a pairing is the top graded piece, described by $\rho^e$, and the non-triviality statement one needs there is delicate: it is at this point that Kato reduces to the case of a residue field which is finitely generated over $\F_p$, a reduction relying, among other tools, on the existence of Cohen rings. Since $\Gamma$ is dense in $\R$, we may instead fix an auxiliary element $\gamma_\infty\in p\Gamma$ with $\gamma_\infty<e$ but \emph{close enough} to $e$, and use $\rho^{\gamma_\infty}$ as the target of the pairing (see Lemma \ref{lem calcul}); of course, such a choice makes no sense when $v$ is discrete. The benefit is twofold. On the one hand, the non-triviality result that we actually need can be established directly by an argument relying on the norm varieties of Suslin and Joukhovitski (Proposition \ref{lem triv symb} and Corollary \ref{cor non trivial}). In particular, we never have to worry about the possible triviality of the source of $\rho_q^e$. On the other hand, this allows us to completely avoid the reduction to finitely generated residue fields, and hence the use of Cohen rings, which are simply not available in our context because the valuation rings we work with are not Noetherian.

The second novelty is of a different nature and is forced upon us by the density of $\Gamma$ in $\mathbb{R}$. In the discrete case, the filtration $(u_q^n(K))_{n\geq0}$ is indexed by $\Z_{\geq 0}$ and stops at $e$, so that the vanishing of all its graded pieces in a given degree $q$ immediately yields $\kk_q(K)=0$. Here the filtration is indexed by the dense set $\Gamma\cap[0,e]$, and a class $a\in\kk_q(K)$ may a priori escape to a limit: it may very well lie in $u_q^\gamma(K)$ for every $\gamma<\gamma_\infty$ without lying in $u_q^{\gamma_\infty}(K)$, in which case it would remain invisible to all the graded pieces. Ruling out this phenomenon for maximally complete fields is the content of Theorem \ref{th invisible terms}, and the whole of Section \ref{sec invisible} is devoted to its proof. That proof is a transfinite induction, in the course of which one produces pseudo-Cauchy sequences whose pseudo-limits, provided by Kaplansky's characterization of maximally complete fields \cite{Kaplansky}, yield the elements we are looking for. This is the only place where maximality is used, and this is what makes Theorem \ref{thm cd} an equality.

\paragraph{Structure of the article}
Section \ref{subsec setting} below fixes the conventions and the notation, and Section \ref{subsec main results} contains the precise statements of our three main results. Section \ref{sec prelim} gathers a number of elementary computations in Milnor $\mathrm{K}$-theory, which will be used over and over again. Section \ref{sec proof isom} is devoted to the proof of Theorem \ref{thm isom Kato} and is divided into six steps: after checking that the maps $\rho_q^\gamma$ are well defined and surjective (Step 1), we deal successively with the case $\gamma=0$ (Step 2), with the non-triviality statements alluded to above (Step 3), with the injectivity of $\rho_q^\gamma$ for $0<\gamma<e$ (Step 4) and, finally, with the case $\gamma=e$, which is settled by a diagram chase involving the maximal unramified extension of $K$ (Step 5). Section \ref{sec invisible} contains the study of maximally complete fields and the proof of Theorem \ref{th invisible terms}. In Section \ref{sec cd} we put everything together and prove Theorem \ref{thm cd}.

\subsection{Setting, conventions and notation}\label{subsec setting}

\paragraph{Dimension of a field}
The cohomological dimension $\mathrm{cd}(K)$ of a field $K$ is the cohomological dimension of its absolute Galois group. In other words, it is the smallest integer $\delta$ (or $\infty$ if such an integer does not exist) such that $H^n(K,M)=0$ for every $n>\delta$ and each finite Galois module $M$. This notion turns out not to be satisfactory for imperfect fields, in particular because a field with characteristic $p>0$ and whose Galois group is a pro-$p$-group has cohomological dimension at most $1$. Better definitions in that context were provided by Kato in \cite{Kato}.

\begin{definition}\label{def cohdim}
Let $K$ be any field and $p$ a prime number. 
\begin{itemize}
\item[(i)] Assume that the characteristic of $K$ is different from $p$. The \emph{$p$-dimension} $\mathrm{dim}_{p}(K)$ of $K$ is the $p$-cohomological dimension of the absolute Galois group of $K$.
\item[(ii)] Assume now that $K$ has characteristic $p$. Let $\Omega^i_K$ be the $i$-th exterior product over $K$ of the absolute differential module $\Omega^1_{K/\mathbb{Z}}$ and let $H^{i+1}_p(K)$ be the cokernel of the morphism $\wp^i_K\,:\, \Omega^i_K \to \Omega^i_K/d(\Omega^{i-1}_K)$ defined by
\[x \frac{dy_1}{y_1} \wedge\ldots \wedge \frac{dy_i}{y_i} \mapsto (x^p-x) \frac{dy_1}{y_1} \wedge\ldots \wedge \frac{dy_i}{y_i} \pmod {d(\Omega^{i-1}_K)},\]
for $x\in K$ and $y_1,\ldots,y_i \in K^{\times}$. The \emph{$p$-dimension} $\mathrm{dim}_{p}(K)$ of $K$ is the smallest integer $\delta$ (or $\infty$ if such an integer does not exist) such that $[K:K^p]\leq p^\delta$ and $H^{\delta+1}_p(L)=0$ for all finite extensions $L$ of $K$.
\item[(iii)] The \emph{dimension} $\mathrm{dim}(K)$ of $K$ is the supremum of all the $\mathrm{dim}_{\ell}(K)$'s when $\ell$ runs through all prime numbers.
\end{itemize}
\end{definition}

\paragraph{Valued fields}
Let $\Gamma$ be a rank 1 totally ordered abelian group and let $K$ be a Henselian field endowed with a surjective valuation $v: K^\times \rightarrow \Gamma \subseteq\R$. We denote by
$\mathcal{O}_K$ the corresponding valuation ring, by $\mathcal{O}_{K}^{\times}$ its group of units, and by $k$ its residue field. For $x\in\mathcal{O}_K$, we write $\bar{x}$ for its image in $k$. For any algebraic extension $L/K$, the valuation $v$ extends in a unique way, so that we will abusively still denote it by $v$ and use the notations $\cal O_L$, $l$ for the respective valuation ring and residue field.

Assume that $(K,v)$ has mixed characteristic $(0,p)$ and that $K$ contains a primitive $p$-th root of unity $\zeta_p$. Note that this implies that $v(p)$ is divisible by $p-1$, so that we may define $e:=\frac{pv(p)}{p-1}=pv(\zeta_p-1)\in\Gamma$.

We assume throughout that $\Gamma$ is \emph{dense} in $\R$, that is, that $v$ is not discrete. The case of discretely valued fields was already covered by Kato in \cite{Kato}.

Finally, we fix once and for all the two integers
\[g:=\dim_{\F_p}(\Gamma/p\Gamma)\in\Z_{\geq 0}\cup\{\infty\}, \qquad p^r:=[k:k^p], \quad r\in\Z_{\geq 0}\cup\{\infty\},\]
so that $g$ is the $p$-rank of the value group and $r$ is the $p$-degree of the residue field.\\

In order to get an equality for the cohomological dimension of such a field, we will need an extra hypothesis besides it being Henselian. Recall that $(K,v)$ is said to be maximal if every proper extension of $K$ (not necessarily algebraic) must either enlarge the residue field $k$ or the value group $\Gamma$. Such a field is necessarily Henselian, complete and defectless \cite[Thm.~2.12]{BlaszczokKuhlmann}, that is, every finite extension of degree $n$ satisfies the equality $n=ef$, where $e$ is the ramification degree and $f$ is the degree of the residue field extension.

In \cite[Thm.~4]{Kaplansky}, Kaplansky proves that $(K,v)$ is maximal if and only if it is maximally complete, which means that every pseudo-Cauchy sequence (called ``pseudo-convergent set'' in \textit{loc.~cit.}) admits a pseudo-limit. Recall that a sequence $(a_\mu)_{\mu<\lambda}$ in $K$ indexed by ordinals $\mu<\lambda$ is a pseudo-Cauchy sequence if
\[v(a_\tau-a_\sigma)>v(a_\sigma-a_\rho)\quad\text{whenever }\rho<\sigma<\tau<\lambda.\]
A pseudo-limit of a pseudo-Cauchy sequence $(a_\mu)_{\mu<\lambda}$ is an element $a\in K$ such that
\[v(a-a_\mu)>v(a-a_\nu)\quad\text{whenever }\nu<\mu<\lambda.\]
This is the actual property we will be working with in the proofs.

\paragraph{Milnor $\mathrm{K}$-theory of valued fields}
We fix $q$ a positive integer and we set $\kk_q(K):=\K_q(K)/p$ and $\kk_q(k):=\K_q(k)/p$. We let $U_q(K)$ be the subgroup of $\K_q(K)$ generated by symbols of the form $\{ x_{1},\dots,x_{q}\}$ with $x_{1},\dots,x_{q} \in \mathcal{O}_K^\times$, and given some non-negative $\gamma\in \Gamma$, we let $U_q^\gamma(K)$ (resp. $U_q^{>\gamma}(K)$) be the subgroup generated by symbols that lie in $U_q(K)$ and that are  of the form $\{ x_{1},\dots,x_{q}\}$ with $x_{1},\dots,x_{q} \in K^\times$ and $v(x_1-1)\geq \gamma$ (resp. $v(x_1-1)> \gamma$). We denote by $u_q(K)$, $u_q^{\gamma}(K)$ and $u_q^{>\gamma}(K)$ the respective images of $U_q(K)$, $U_q^{\gamma}(K)$ and $U_q^{>\gamma}(K)$ in $\kk_q(K)$. Finally, for $\gamma_\infty\in\R_{>0}$, we set
\[u_q^{<\gamma_\infty}(K) := \bigcap_{\gamma < \gamma_\infty} u_q^\gamma(K) \;\;\;\;\; \text{and} \;\;\;\;\; u_q^{\geq \gamma_\infty}(K) := \bigcup_{\gamma \geq \gamma_\infty} u_q^\gamma(K),\]
where $\gamma$ runs through $\Gamma_{\geq0}$. Since $u_q^\gamma(K)$ decreases with $\gamma$, we always have $u_q^{\geq\gamma_\infty}(K)\subseteq u_q^{<\gamma_\infty}(K)$; the reverse inclusion when $K$ is maximal is the content of Theorem \ref{th invisible terms} below.\\

\subsection{Main results}\label{subsec main results}
The valuation $v$ induces a surjective $\mathbb{F}_p$-linear map:
$$ K^\times/{K^\times}^p\rightarrow \Gamma/p\Gamma.$$
We let $s$ be some $\mathbb{F}_p$-linear section of this map and for each $\bar{\gamma} \in \Gamma/p\Gamma$, we denote by $\pi^{\bar{\gamma}}$ some preimage in $K^\times$ of the class $s(\bar\gamma)\in K^\times/{K^\times}^p$. Also, for $\gamma \in \Gamma$, we denote by $\pi^\gamma$ a fixed element of $K$ with valuation $\gamma$ such that $\pi^\gamma (\pi^{\bar\gamma})^{-1}$ is a $p$-th power in $K^\times$. Note that in general $\pi^\gamma\neq \pi^{\bar\gamma}$.

\begin{theorem}\label{thm isom Kato}
Recall that $\Gamma$ is dense in $\mathbb{R}$ and that $\zeta_p\in K$. For any choice of $s$, $\pi^{\bar\gamma}$ and $\pi^\gamma$ as above, and for every $\gamma \in \Gamma_{\geq 0}$, we have the following well-defined isomorphisms:
\begin{enumerate}[(i)]
     \item if $\gamma=0$:
    \begin{align*}
        \rho_q^0: \bigoplus_{j=0}^{q} \left(\Lambda^j(\Gamma/p\Gamma) \otimes \kk_{q-j}(k)\right) & \rightarrow \kk_q(K)/u_q^{>0}(K)\\
        (\bar{\gamma}_1 \wedge \dots \wedge \bar{\gamma}_j)\otimes\{\bar{x}_{j+1},\dots,\bar{x}_{q}\} & \mapsto \{\pi^{\bar{\gamma}_1},\dots ,\pi^{\bar{\gamma}_j}, x_{j+1},\dots,x_{q}\};
         \end{align*}
\item if $\gamma < e$ and $\gamma\not\in p\Gamma$, for any decomposition $\Gamma/p\Gamma=\langle\bar\gamma\rangle\oplus\Delta$:
          \begin{align*}
        \rho_q^\gamma: \bigoplus_{j=0}^{q-1} \left(\Lambda^j\Delta \otimes \Omega_k^{q-j-1}\right) & \rightarrow u_q^{\gamma}(K)/u_q^{>\gamma}(K)\\
                (\bar{\gamma}_1 \wedge \dots \wedge \bar{\gamma}_j)\otimes \left\{\bar{x} \frac{d \bar{y}_{j+2}}{\bar{y}_{j+2}}\wedge \dots \wedge\frac{d \bar{y}_{q}}{\bar{y}_{q}}\right\} & \mapsto \{\pi^{\bar{\gamma}_1},\dots ,\pi^{\bar{\gamma}_j}, 1+\pi^{\gamma}x, y_{j+2},\dots,y_{q}\};
         \end{align*}
\item if $\gamma < e$ and $\gamma=p\gamma'\in p\Gamma\smallsetminus \{0\}$:
    \begin{align*}
        \rho_q^\gamma: \bigoplus_{j=0}^{q-1} \left(\Lambda^j(\Gamma/p\Gamma) \otimes (\Omega_k^{q-j-1}/\ker (d))\right) & \rightarrow u_q^{\gamma}(K)/u_q^{>\gamma}(K)\\
                (\bar{\gamma}_1 \wedge \dots \wedge \bar{\gamma}_j)\otimes\left\{\bar{x} \frac{d \bar{y}_{j+2}}{\bar{y}_{j+2}}\wedge \dots \wedge\frac{d \bar{y}_{q}}{\bar{y}_{q}}\right\} & \mapsto \{\pi^{\bar{\gamma}_1},\dots ,\pi^{\bar{\gamma}_j}, 1+(\pi^{\gamma'})^px, y_{j+2},\dots,y_{q}\};
    \end{align*}
\item if $\gamma=e$:
 \begin{align*}
        \rho_q^e: \bigoplus_{j=0}^{q-1} \left(\Lambda^j(\Gamma/p\Gamma) \otimes H_p^{q-j}(k)\right) & \rightarrow u_q^{e}(K)/u_q^{>e}(K)\\
                (\bar{\gamma}_1 \wedge \dots \wedge \bar{\gamma}_j)\otimes\left\{\bar{x} \frac{d \bar{y}_{j+2}}{\bar{y}_{j+2}}\wedge \dots \wedge\frac{d \bar{y}_{q}}{\bar{y}_{q}}\right\} & \mapsto \{\pi^{\bar{\gamma}_1},\dots ,\pi^{\bar{\gamma}_j}, 1+(\zeta_p-1)^p x, y_{j+2},\dots,y_{q}\};
    \end{align*}
\item if $\gamma>e$, then $u_q^{\gamma}(K)=0$.
         \end{enumerate}
\end{theorem}

Theorem \ref{thm isom Kato} describes every graded piece of the filtration $(u_q^\gamma(K))_{\gamma\in\Gamma_{\geq0}}$. In contrast with the discrete case, this is not yet enough to control $\kk_q(K)$ itself: as $\Gamma$ is dense in $\R$, the filtration is indexed by a dense set, and a class could a priori lie in $u_q^\gamma(K)$ for every $\gamma<\gamma_\infty$ without lying in $u_q^{\gamma_\infty}(K)$, thus remaining invisible to all the graded pieces. Our second main result rules this out for maximally complete fields.

\begin{theorem}\label{th invisible terms}
Recall that $\Gamma$ is dense in $\mathbb{R}$ and that $\zeta_p\in K$. Assume that $(K,v)$ is maximally complete (or equivalently, maximal). Let $\gamma_\infty \in \mathbb{R}_{\geq 0}$. If $\gamma_\infty\leq e=\frac{p v(p)}{p - 1}$, then $u_q^{<\gamma_\infty}(K) = u_q^{\geq \gamma_\infty}(K)$.
\end{theorem}

This is the only point of the article where maximality is used. Combining Theorems \ref{thm isom Kato} and \ref{th invisible terms}, we obtain our main application: a computation of the $p$-cohomological dimension of $K$ in terms of the $p$-rank $g$ of the value group and of Kato's $p$-dimension of the residue field.

\begin{theorem}\label{thm cd}
Let $(K,v)$ be a Henselian valued field of rank 1 of mixed characteristic $(0,p)$. Then
\[\dim_p(K) \geq \dim_{\F_p}(\Gamma/p\Gamma)+\dim_p(k).\]
Moreover, equality holds if $(K,v)$ is maximally complete (or equivalently, maximal).
\end{theorem}

This last equality is easy to settle in the case of prime-to-$p$ cohomological dimension, without the maximality assumption, as it is done for instance in \cite[Prop.~1.2]{Soman}. Putting this together with Theorem \ref{thm cd}, we get:

\begin{corollary}\label{cor cd}
Let $(K,v)$ be a Henselian valued field of rank 1 of characteristic $0$ and let $p$ be a prime number. Then
\[\dim_p(K) \geq \dim_{\F_p}(\Gamma/p\Gamma)+\dim_p(k),\]
with equality whenever $p\neq\mathrm{char}(k)$ or $(K,v)$ is maximally complete (or equivalently, maximal).
\end{corollary}

\begin{remarque}\label{rem cd}
Theorem \ref{thm cd} and Corollary \ref{cor cd} hold even without the hypotheses $\zeta_p\in K$ and $\Gamma$ dense in $\R$. For the former, this is simply because the extension $K(\zeta_p)/K$ is always of degree prime to $p$: it therefore does not affect $\dim_p(K)$, $\dim_p(k)$ (by a restriction-corestriction argument) or $\dim_{\F_p}(\Gamma/p\Gamma)$ (see Lemma \ref{lem gps}), and the maximality assumption is preserved by \cite[Thm.~2.14]{BlaszczokKuhlmann}. The latter hypothesis can be dropped because the discrete case has already been dealt with by Kato. This shows that it suffices to prove Theorem \ref{thm cd} under those hypothesis. That is what we will do in Section \ref{sec cd}.
\end{remarque}

\begin{remarque}
    Theorem \ref{thm cd} and Corollary \ref{cor cd} can be easily generalized to Henselian valued fields of finite rank by an inductive argument. The only slightly nontrivial observation is that, when the field is maximal of rank $n>1$ and we consider the corresponding rank one valuation, its residue field (which is Henselian of rank $n-1$) is also maximal.
\end{remarque}

\begin{remarque}\label{rem cyclotomic}
An extra hypothesis such as maximality is necessary to get equality in Theorem \ref{thm cd} and in Corollary \ref{cor cd} when $p=\mathrm{char}(k)$, as one can see by considering the maximal cyclotomic extension $K$ of $\Q_p$. This is a Henselian valued field with algebraically closed residue field $k$, so that $\dim_p(k)=0$, and value group
\[\Gamma=\left\{\frac{a}{(p-1)p^r}\mid a,r\in\Z\right\}\subseteq\R,\]
which is clearly $p$-divisible, so that $\dim_{\F_p}(\Gamma/p\Gamma)=0$. However, $\sqrt[p]{p}\not\in K$ (or $\sqrt[4]{2}$ if $p=2$), since the corresponding extension of $\Q_p$ is not abelian. Hence $K$ has at least one extension of degree $p$. This yields
\[\dim_p(K)\geq 1> 0+0=\dim_{\F_p}(\Gamma/p\Gamma)+\dim_p(k),\]
and also proves that the field is not defectless, hence not maximal.
\end{remarque}

\paragraph{Acknowledgements}
The second author was partially supported by ANID via Proyecto FONDECYT Regular N\textsuperscript{o}1240001. The third named author gratefully acknowledges support from ANID Fondecyt Postdoctoral Grant 3240191 and from ANID Concurso de Subvenci\'on a la Instalaci\'on en la Academia 2025 N°85250089.

\section{Preliminary computations}\label{sec prelim}

In this section we carry out straightforward computations in Milnor $\mathrm{K}$-theory that will be used in the proof of the main results. In all the subsequent lemmas, the hypotheses are those of Theorem \ref{thm isom Kato}. In particular, $(K,v)$ is a Henselian valued field of rank 1 with mixed characteristic $(0,p)$ and value group $\Gamma$ dense in $\R$.\\

We start with a description of the generators of $u_q^{\gamma}(K)$.

\begin{lemma}\label{lem generateurs u_q}
    Take $\gamma \in \Gamma_{>0}$. 
    \begin{enumerate}
        \item Every symbol of the form $\{1+\pi^\gamma x,y_2,\dots, y_q\}$ with $x\in \mathcal{O}_K$ and $y_2,\dots, y_q \in \mathcal{O}_K \setminus \{0\}$ can also be written as $\{z_1,z_2,\dots, z_q\}$ with $z_1,\dots,z_q\in\mathcal{O}_K^\times$.
        \item The group $u_q^{\gamma}(K)$ is generated by the symbols of the form $\{1+\pi^\gamma x,y_2,\dots, y_q\}$ with $x\in \mathcal{O}_K$ and $y_2,\dots y_q \in \mathcal{O}_K \setminus \{0\}$.
    \end{enumerate}
\end{lemma}

\begin{proof}\,
\begin{enumerate}
    \item Consider a symbol of the form $\{1+\pi^\gamma x,y_2,\dots, y_q\}$ with $x\in \mathcal{O}_K$ and $y_2,\dots ,y_q \in \mathcal{O}_K \setminus \{0\}$. Since $x\in \cal O_K$, up to replacing $\gamma$ by a larger element in $\Gamma$, we may assume that $x\in\cal O_K^\times$. Now observe that, in $\K_q(K)$:
    $$\{1+\pi^\gamma x,y_2,\dots, y_q\}= \{1+\pi^\gamma x,y_2(-\pi^\gamma x)^{i_2},\dots, y_q(-\pi^\gamma x)^{i_q}\}$$
    for all $i_2,\dots,i_q\in \mathbb{Z}$. By choosing the integers $i_2,\dots,i_q$ conveniently, by replacing $y_2,\dots , y_q$ by $y_2(-\pi^\gamma x)^{i_2},\dots, y_q(-\pi^\gamma x)^{i_q}$ respectively and by reordering the $y_i$'s, we may and do assume that $0\leq v(y_q) \leq \dots \leq v(y_2)<\gamma$. Then, if $v(y_2)>0$, we have the following equalities in $\K_q(K)$:
    \begin{align*}
    \{1+\pi^\gamma x,y_2,\dots, y_q\}&=\{(1+\pi^\gamma x)(1-y_2),y_2,\dots, y_q\}\\&=\{1+\pi^\gamma x-y_2-y_2\pi^\gamma x,y_2(-\pi^\gamma x+y_2+y_2\pi^\gamma x)^{-1},y_3,\dots, y_q\}
    \end{align*}
    where $v(y_2(-\pi^\gamma x+y_2+y_2\pi^\gamma x)^{-1})=0$. Hence, up to replacing $\gamma$ by $v(y_2)$, the element $x$ by $(\pi^{v(y_2)})^{-1}(\pi^\gamma x-y_2-y_2\pi^\gamma x)$ and $y_2$ by $y_2(-\pi^\gamma x+y_2+y_2\pi^\gamma x)^{-1}$, we may and do assume that $y_2 \in \mathcal{O}_K^\times$. By repeating this argument for $y_3,\dots,y_q$ in this order, we may and do assume that $y_3,\dots,y_q$ are all units. This shows that there exist $z_1,\dots,z_q\in \mathcal{O}_K^\times$ such that $\{1+\pi^\gamma x,y_2,\dots, y_q\}=\{z_1,z_2,\dots, z_q\}$, as wished.
\item By (i), the symbols of the form $\{1+\pi^\gamma x,y_2,\dots, y_q\}$ with $x\in \mathcal{O}_K$ and $y_2,\dots y_q \in \mathcal{O}_K \setminus \{0\}$ lie in $u_q^{\gamma}(K)$. Moreover, if we take a symbol of the form $\{1+\pi^\gamma x,y_2,\dots, y_q\}$ with $x\in \mathcal{O}_K$ and $y_2,\dots y_q \in K^\times$, we can find $\gamma_2,\dots,\gamma_q \in \Gamma$ such that the elements $(\pi^{\gamma_2})^p y_2,\dots,(\pi^{\gamma_q})^p y_q$ all lie in $\mathcal{O}_K$. But in $\kk_q(K)$:
$$\{1+\pi^\gamma x,y_2,\dots, y_q\}=\{1+\pi^\gamma x,(\pi^{\gamma_2})^p y_2,\dots, (\pi^{\gamma_q})^p y_q\}.$$
The desired result follows by definition of $u_q^{\gamma}(K)$.
\end{enumerate}
\end{proof}

We now prove the case $\gamma>e$ of Theorem \ref{thm isom Kato}.

\begin{lemma}\label{lem gamma sup e}
Let $\gamma\in\Gamma$ with $\gamma>e$ and let $x\in\cal O_K^\times$. Then $1+\pi^\gamma x\in {K^\times}^p$. Consequently $U_1^{\gamma}(K)\subseteq {K^\times}^p$ and, for every $q\geq 1$, $u_q^{\gamma}(K)=u_q^{>\gamma}(K)=0$.
\end{lemma}

\begin{proof}
The power series for $p$-th roots applied in this context reads
\[(1+\pi^\gamma x)^{1/p}=\sum_{i=0}^\infty \binom{1/p}{i} (\pi^{\gamma})^ix^i.\]
For each term in the series we have
\[v\left(\binom{1/p}{i}(\pi^{\gamma})^ix^i\right)=v\left(\binom{1/p}{i}\right)+i\gamma=v_p\left(\binom{1/p}{i}\right)v(p)+i\gamma,\]
where $v_p$ stands for the usual $p$-adic valuation. Since $\gamma=e+\epsilon$, for some $\epsilon\in\R_{>0}$, we get
\begin{align*}
v\left(\binom{1/p}{i}(\pi^{\gamma})^ix^i\right) &=v(p)\left(v_p\left(\binom{1/p}{i}\right)+\frac{ie}{v(p)}\right)+i\epsilon\\
&=v(p)\left(v_p\left(\binom{1/p}{i}\right)+\frac{ip}{p-1}\right)+i\epsilon.
\end{align*}
Since $v_p(\binom{1/p}{i})\geq -\frac{ip}{p-1}$, this last term tends to infinity when $i$ grows. We see then that the series converges in the completion of $K$ and hence $1+\pi^\gamma x$ admits a $p$-th root since $K$ is Henselian \cite[Cor.~1.2.1]{Moret-Bailly Greenberg}. The last assertion follows from Lemma \ref{lem generateurs u_q}, since every generator of $u_q^\gamma(K)$ has a slot of the form $1+\pi^{\delta}y$ with $\delta\geq\gamma>e$.
\end{proof}

This proves item (v) of Theorem \ref{thm isom Kato}, and we will tacitly use it from now on, assuming by default that all symbols in $u_q^{>e}(K)$ are trivial.\\

The remaining lemmas will be useful in several computations below.

\begin{lemma}\label{lem suma}
    Let $\gamma,\delta\in \Gamma$ with $\gamma \geq \delta$ and $x,y\in \mathcal{O}_K$, and set $\alpha:=\{1+\pi^\gamma x,1+\pi^\delta y\}$. Then:
    $$\alpha\equiv \{1-\pi^\gamma\pi^\delta xy,\pi^\gamma x(-1-\pi^{\delta}y)\}\pmod {u^{2\gamma+\delta}_2(K)}.$$
    In particular, $\alpha\in u^{\gamma+\delta}_2(K)$, and if $\gamma+2\delta>e$, then $\alpha=\{1-\pi^\gamma\pi^\delta xy,-\pi^\gamma x\}$.
\end{lemma}

\begin{proof}
First note that, for $\gamma\in\Gamma_{>0}$ and $z\in\cal O_K$, expanding the inverse gives
\[0=-\{1-\pi^\gamma z,\pi^\gamma z\}=\{(1-\pi^\gamma z)^{-1},\pi^\gamma z\}\equiv\{1+\pi^\gamma z+(\pi^\gamma z)^2,\pi^\gamma z\}\pmod {u^{3\gamma}_2(K)}.\]
Applying this to $z=-x-\pi^\delta xy$, we get
\begin{align*}
    \alpha&=\{1+\pi^\gamma x,1+\pi^\delta y\}\\
    &=\{1+\pi^\gamma x,\pi^\gamma x(-1-\pi^{\delta}y)\}\\
    &\equiv\{(1+\pi^\gamma x)(1-\pi^\gamma x-\pi^{\gamma}\pi^{\delta}xy+(\pi^{\gamma}x)^2),\pi^\gamma x(-1-\pi^{\delta}y)\}\pmod {u^{2\gamma+\delta}_2(K)}\\
    &\equiv \{1-\pi^\gamma\pi^\delta xy,\pi^\gamma x(-1-\pi^{\delta}y)\}\pmod {u^{2\gamma+\delta}_2(K)}.
\end{align*}
In particular, we have $\alpha\in u^{\gamma+\delta}_2(K)$.

Assume now that $\gamma+2\delta>e$. We then have $u^{2\gamma+\delta}_2(K)=0$ by Lemma \ref{lem gamma sup e}, since $2\gamma+\delta\geq\gamma+2\delta>e$, and hence:
$$\alpha=\{1-\pi^\gamma\pi^\delta xy,-\pi^\gamma x(1+\pi^{\delta}y)\}.$$
But $\{1-\pi^\gamma\pi^\delta xy,1+\pi^{\delta}y\}\in u_2^{\gamma+2\delta}(K)=0$. Hence $\alpha=\{1-\pi^\gamma\pi^\delta xy,-\pi^\gamma x\}$.
\end{proof}

\begin{lemma}\label{lem sum symb}
    For $\gamma>0$ and $x,y_1,y_2\in \mathcal{O}_K\setminus \{0\}$ such that $y_1+y_2\in \mathcal{O}_K\setminus\{0\}$:
    $$\{1+\pi^{\gamma}x(y_1+y_2),y_1+y_2\}\equiv \{1+\pi^{\gamma}xy_1,y_1\}+\{1+\pi^{\gamma}xy_2,y_2\} \pmod {u^{2\gamma}_2(K)}.$$
\end{lemma}

\begin{proof}
Modulo $u^{2\gamma}_2(K)$, we have the following equivalences:
    \begin{align*}
        \{1+\pi^{\gamma}x(y_1+y_2),y_1+y_2\}
        &\equiv\{1+\pi^{\gamma}x(y_1+y_2),-\pi^{-\gamma}x^{-1}\}\\
        &\equiv \{1+\pi^{\gamma}xy_1,-\pi^{-\gamma}x^{-1}\}+\{1+\pi^{\gamma}xy_2,-\pi^{-\gamma}x^{-1}\}\\
        &\equiv \{1+\pi^{\gamma}xy_1,y_1\}+\{1+\pi^{\gamma}xy_2,y_2\}\pmod{u^{2\gamma}_2(K)}.
    \end{align*}
\end{proof}

\begin{lemma}\label{lem trivialite symboles}
Let $\gamma\in\Gamma_{>0}$ and $q\geq 2$ and let $\alpha=\{1+\pi^\gamma x,y_1,\ldots,y_{q-1}\}\in u_q^{\gamma}(K)$ with $x,y_i\in\cal O_K^\times$. Assume that either $\bar y_i=1$ for some $i$, or $\bar y_i=\bar y_j$ for some $i\neq j$. Then $\alpha\in u_q^{>\gamma}(K)$.
\end{lemma}

\begin{proof}
The statement under the first assumption follows directly from Lemma \ref{lem suma}, since then $y_i=1+\pi^\delta z$ for some $z\in\cal O_K$.

Let us prove now the statement under the second assumption. It suffices to do this for $q=3$, so that $\alpha=\{1+\pi^\gamma x,y,z\}$ with $\bar y=\bar z$. Write $z=-yu$ with $u\in\cal O_K^\times$ and $\bar u=-1$. Then
\[\alpha=\{1+\pi^\gamma x,y,-y\}+\{1+\pi^\gamma x,y,u\},\]
so that it suffices to prove that these two symbols lie in $u_q^{>\gamma}(K)$. The first one is obvious since $\{y,-y\}=0$. For the second one, if $p=2$, then $\bar u=1$ and we are reduced to the first case. If $p$ is odd, then $-1$ is a $p$-th power and hence we may replace $u$ with $-u$ and we are done as well.
\end{proof}

\section{Proof of Theorem \ref{thm isom Kato}}\label{sec proof isom}

\subsection{Step 1: the morphisms $\rho_q^\gamma$ are well-defined and surjective}

Fix $\gamma \in \Gamma_{\geq 0}$ and $q\geq 0$. We check that the morphism $\rho^\gamma_q$ is well-defined and surjective.\\

We start with $\gamma=0$. Fix some $0 \leq j \leq q$. Then given $x_{j+1},\dots,x_{q} \in K^\times$, the map:
\begin{align*}
    (\Gamma/p\Gamma)^j & \rightarrow \kk_q(K)/u_q^{>0}(K)\\
    (\bar{\gamma}_1 , \dots , \bar{\gamma}_j) & \mapsto \{\pi^{\bar{\gamma}_1},\dots ,\pi^{\bar{\gamma}_j}, x_{j+1},\dots,x_{q}\},
\end{align*}
is well-defined, multilinear since the section $(\Gamma/p\Gamma)\to (K^\times/{K^\times}^p):\bar\gamma\mapsto\pi^{\bar\gamma}$ is linear and alternating since $\{\pi^{\bar\gamma},\pi^{\bar\gamma}\}=\{\pi^{\bar\gamma},-1\}\in u_2^{>0}(K)$. Now, given $\bar{\gamma}_1 , \dots , \bar{\gamma}_j \in \Gamma/p\Gamma$, the map:
\begin{align*}
(k^\times)^{q-j} & \rightarrow \kk_q(K)/u_q^{>0}(K)\\
    (\bar{x}_{j+1},\dots,\bar{x}_{q}) & \mapsto \{\pi^{\bar{\gamma}_1},\dots ,\pi^{\bar{\gamma}_j}, x_{j+1},\dots,x_{q}\},
\end{align*}
is well-defined by Lemma \ref{lem trivialite symboles}. Moreover, it is multilinear and its kernel contains all $(q-j)$-tuples containing two entries adding up to $1$. We therefore get a well-defined map,
\begin{align*}
          \Lambda^j(\Gamma/p\Gamma) \times \kk_{q-j}(k) & \rightarrow \kk_q(K)/u_q^{>0}(K)\\
        (\bar{\gamma}_1 \wedge \dots \wedge \bar{\gamma}_j,\{\bar{x}_{j+1},\dots,\bar{x}_{q}\}) & \mapsto \{\pi^{\bar{\gamma}_1},\dots ,\pi^{\bar{\gamma}_j}, x_{j+1},\dots,x_{q}\},
         \end{align*}
which is bilinear. It therefore induces a morphism:
\begin{align*}
        \rho_q^0: \bigoplus_{j=0}^{q} \left(\Lambda^j(\Gamma/p\Gamma) \otimes \kk_{q-j}(k)\right) & \rightarrow \kk_q(K)/u_q^{>0}(K)\\
        (\bar{\gamma}_1 \wedge \dots \wedge \bar{\gamma}_j)\otimes\{\bar{x}_{j+1},\dots,\bar{x}_{q}\} & \mapsto \{\pi^{\bar{\gamma}_1},\dots ,\pi^{\bar{\gamma}_j}, x_{j+1},\dots,x_{q}\}.
\end{align*}
In order to prove its surjectivity, it suffices to note that every slot in any symbol in $\kk_q$ can be written as $\pi^{\bar\gamma}x$ for some $\bar\gamma\in\Gamma/p\Gamma$ and $x\in\cal O_K^\times$, hence any symbol can be written as a sum of symbols of the form above by multilinearity.\\

We move on to $\gamma>0$. Consider a decomposition $\Gamma/p\Gamma=\langle\bar\gamma\rangle\oplus\Delta$ (note that if $\gamma\in p\Gamma$, then $\Delta=\Gamma/p\Gamma$). As before, fix some $0 \leq j \leq q-1$. Given $x\in\cal O_K$ and $y_{j+2},\dots,y_{q} \in \cal O_K^\times$, the map:
\begin{align*}
    \Delta^j & \rightarrow u_q^\gamma(K)/u_q^{>\gamma}(K)\\
    (\bar{\gamma}_1 , \dots , \bar{\gamma}_j) & \mapsto \{\pi^{\bar{\gamma}_1},\dots ,\pi^{\bar{\gamma}_j}, 1+\pi^\gamma x,y_{j+2},\dots,y_{q}\},
\end{align*}
is again well-defined, multilinear and alternating. And given $\bar{\gamma}_1 , \dots , \bar{\gamma}_j \in \Delta$, the map:
\begin{align*}
(\Omega_k^1)^{q-j-1} & \rightarrow u_q^\gamma(K)/u_q^{>\gamma}(K)\\
    (\bar x'_{j+2}d\bar y_{j+2},\dots,\bar x'_{q}d\bar{y}_{q}) & \mapsto \{\pi^{\bar{\gamma}_1},\dots ,\pi^{\bar{\gamma}_j}, 1+\pi^\gamma (x'_{j+2}y_{j+2}\cdots x'_{q}y_q),y_{j+2},\dots,y_{q}\},
\end{align*}
is well-defined by Lemmas \ref{lem generateurs u_q}, \ref{lem sum symb} and  \ref{lem trivialite symboles}. It is also clearly multilinear, and alternating by Lemma \ref{lem trivialite symboles}. Rewriting $\bar x'_id\bar y_i$ as $\bar x_i\frac{d\bar y_i}{\bar y_i}$ and $\bar x=\bar x_1\cdots\bar x_q$, we get a well-defined map,
\begin{align*}
          \Lambda^j\Delta \times \Omega_k^{q-j-1} & \rightarrow u_q^\gamma(K)/u_q^{>\gamma}(K)\\
        \left(\bar{\gamma}_1 \wedge \dots \wedge \bar{\gamma}_j,\bar{x} \frac{d \bar{y}_{j+2}}{\bar{y}_{j+2}}\wedge \dots \wedge\frac{d \bar{y}_{q}}{\bar{y}_{q}}\right) & \mapsto \{\pi^{\bar{\gamma}_1},\dots ,\pi^{\bar{\gamma}_j}, 1+\pi^\gamma x,y_{j+2},\dots,y_{q}\},
         \end{align*}
which is bilinear. It therefore induces a morphism:
\begin{align}
        \bigoplus_{j=0}^{q-1} \left(\Lambda^j\Delta \otimes \Omega_k^{q-j-1}\right) & \rightarrow u_q^\gamma(K)/u_q^{>\gamma}(K) \label{eqn def rhos} \\ 
        (\bar{\gamma}_1 \wedge \dots \wedge \bar{\gamma}_j)\otimes \left(\bar{x} \frac{d \bar{y}_{j+2}}{\bar{y}_{j+2}}\wedge \dots \wedge\frac{d \bar{y}_{q}}{\bar{y}_{q}}\right) & \mapsto \{\pi^{\bar{\gamma}_1},\dots ,\pi^{\bar{\gamma}_j}, 1+\pi^\gamma x,y_{j+2},\dots,y_{q}\}.\nonumber
\end{align}
This is precisely the morphism $\rho_q^\gamma$ if $\gamma\not\in p\Gamma$. Let us prove that the map is also surjective in this case. Reasoning as in the case of $\rho_q^0$ and using Lemma \ref{lem generateurs u_q}, we see that every class in $u_q^{\gamma}(K)$ can be written as a sum of symbols of the form above, but with $\bar\gamma_i\in \Gamma/p\Gamma$ and not necessarily in $\Delta$. However, for $x\in\cal O_K^\times$ we have
\[\{\pi^{\bar\gamma},1+\pi^\gamma x\}=\{\pi^{\bar\gamma}(-\pi^\gamma x)^{-1},1+\pi^\gamma x\},\]
which proves that we may replace any slot with $\pi^{\bar\gamma}$ with a unit. By multilinearity, we see that we may assume that all the $\bar\gamma_i$'s are in $\Delta$, which proves surjectivity.\\

Assume now that $\gamma\in p\Gamma$, so that $\Delta=\Gamma/p\Gamma$ and we may write $\gamma=p\gamma'$. Note that all the computations above work just as fine when replacing $\pi^\gamma$ by $(\pi^{\gamma'})^p$. In order to define $\rho_q^\gamma$, we need to prove then that the map \eqref{eqn def rhos} (replacing $\pi^\gamma$ by $(\pi^{\gamma'})^p$) is trivial on any element of $\Lambda^j(\Gamma/p\Gamma)\otimes \ker(d)$. This amounts to proving that the map
\begin{align}
        \Omega_k^{r} & \rightarrow u_{r+1}^\gamma(K)/u_{r+1}^{>\gamma}(K) \label{eqn def rhos omega}\\
        \bar{x} \frac{d \bar{y}_{1}}{\bar{y}_{1}}\wedge \dots \wedge\frac{d \bar{y}_{r}}{\bar{y}_{r}} & \mapsto \{1+(\pi^{\gamma'})^p x,y_{1},\dots,y_{r}\},\nonumber
\end{align}
is trivial on $\ker(d)\subset\Omega_k^r$ for any $r\geq 0$. By Cartier's Theorem \cite[Thm.~9.4.3]{GS}, any element in $\ker(d)$ can be written as a sum $\alpha+\beta$, where $\alpha\in\im(d)$ and $\beta$ is a finite sum of elements of the form
\begin{equation}\label{eqn diff im AS}
\bar{x}^p \frac{d \bar{y}_{1}}{\bar{y}_{1}}\wedge \dots \wedge\frac{d \bar{y}_{r}}{\bar{y}_{r}}.
\end{equation}
Let us prove that the map \eqref{eqn def rhos omega} is trivial on \eqref{eqn diff im AS}, and hence on $\beta$. For this, it suffices by Lemma \ref{lem generateurs u_q} to prove that the symbol $\{1+(\pi^{\gamma'})^px^p\}\in u_1^{\gamma}(K)$ actually lies in $u_1^{>\gamma}(K)$. Now, since $\gamma<e$, we have $\gamma=p\gamma'<v(p)+\gamma'$ and hence
\[(1+\pi^{\gamma'}x)^p=(1+(\pi^{\gamma'})^px^p)(1+p\pi^{\gamma'}y)=(1+(\pi^{\gamma'})^px^p)(1+\pi^{\delta}z),\]
for some $y,z\in\cal O_K^\times$ and $\delta>p\gamma'=\gamma$. This proves that $\{1+(\pi^{\gamma'})^px^p\}$ lies indeed in $u_1^\delta(K)\subseteq u_1^{>\gamma}(K)$.

Consider now an element $\alpha\in\im (d)\subseteq\Omega_k^r$. A direct computation tells us that it is a finite sum of elements of the form
\[\bar{x}\frac{d\bar x}{\bar x}\wedge \frac{d \bar{y}_{1}}{\bar{y}_{1}}\wedge \dots \wedge\frac{d \bar{y}_{r-1}}{\bar{y}_{r-1}}.\]
Then, to prove that the map \eqref{eqn def rhos omega} is trivial on $\alpha$, it suffices by Lemma \ref{lem generateurs u_q} to prove that the image of $\bar x\frac{d\bar x}{\bar x}$ in $u_2^{\gamma}(K)$, that is the symbol $\{1+(\pi^{\gamma'})^p x,x\}$, is actually in $u_2^{>\gamma}(K)$. Now, we have
\[\{1+(\pi^{\gamma'})^p x,x\}=\{1+(\pi^{\gamma'})^p x,(-(\pi^{\gamma'})^px)^{-1}x\}=\{1+(\pi^{\gamma'})^p x,-1\}.\]
Thus, if $p=2$, we have that the last slot is $1$ on the residue field. And if $p\neq 2$, $-1$ is a $p$-th power and we may replace it with a $1$, so that this still holds. We are done then by Lemma \ref{lem trivialite symboles}.

All in all, we have proved that we have a well-defined map
\begin{align*}
        \rho_q^\gamma: \bigoplus_{j=0}^{q-1} \left(\Lambda^j(\Gamma/p\Gamma) \otimes (\Omega_k^{q-j-1}/\ker (d))\right) & \rightarrow u_q^{\gamma}(K)/u_q^{>\gamma}(K)\\
        (\bar{\gamma}_1 \wedge \dots \wedge \bar{\gamma}_j)\otimes\left\{\bar{x} \frac{d \bar{y}_{j+2}}{\bar{y}_{j+2}}\wedge \dots \wedge\frac{d \bar{y}_{q}}{\bar{y}_{q}}\right\} & \mapsto \{\pi^{\bar{\gamma}_1},\dots ,\pi^{\bar{\gamma}_j}, 1+(\pi^{\gamma'})^px, y_{j+2},\dots,y_{q}\},
\end{align*}
and we still need to prove surjectivity. But this follows immediately from Lemma \ref{lem generateurs u_q} by reasoning as in the previous cases.\\

Finally, let us consider the case where $\gamma=e$. The map \eqref{eqn def rhos} is still well-defined in this context replacing $\pi^{\gamma}$ by $(\zeta_p-1)^p$. Moreover, since $e\in p\Gamma$, the computations from the last case for elements in $\im(d)$ still apply, replacing now $\pi^{\gamma'}$ by $(\zeta_p-1)$, so that the map \eqref{eqn def rhos omega} is trivial on $\im(d)$. Thus, we only need to check that the well-defined map
\begin{align*}
        \Omega_k^{r}/\im(d) & \rightarrow u_{r+1}^e(K)/u_{r+1}^{>e}(K)\\
        \bar{x} \frac{d \bar{y}_{1}}{\bar{y}_{1}}\wedge \dots \wedge\frac{d \bar{y}_{r}}{\bar{y}_{r}} & \mapsto \{1+(\zeta_p-1)^p x,y_{1},\dots,y_{r}\},
\end{align*}
is trivial on $\wp(\Omega_k^r)\subset\Omega_k^r/\im(d)$ for any $q\geq 0$. In order to do this, it will suffice by Lemma \ref{lem generateurs u_q} to prove that the image in $u_{1}^e(K)$ of $\bar x=\bar z^p-\bar z$, which is $\{1+(\zeta_p-1)^p x\}$, is actually in $u_1^{>e}(K)$.

Now, a direct computation gives
\[(1+(\zeta_p-1)z)^p =(1+p(\zeta_p-1)z+(\zeta_p-1)^p z^p)(1+p(\zeta_p-1)^2y)\]
for some $y\in\cal O_K^\times$ (recall that $v(p)=v((\zeta_p-1)^{p-1})$, so that $v((\zeta_p-1)^p)<v(p(\zeta_p-1)^2)$). And since $\frac{(\zeta_p-1)^{p-1}}{p}\equiv -1\pmod p$, we get
\[(1+(\zeta_p-1)z)^p=(1+(\zeta_p-1)^p(z^p-zu))(1+p(\zeta_p-1)^2y),\]
for some $u\in\cal O_K^\times$ with $\bar u=1\in k$. Finally, since $\bar x=\bar z^p-\bar z=\bar z^p-\bar z\bar u$, we get
\[(1+(\zeta_p-1)z)^p=(1+(\zeta_p-1)^px)(1+\pi^\delta y')(1+p(\zeta_p-1)^2y),\]
for some $y'\in\cal O_K^\times$ and $\delta>e$.
The left term gives a trivial symbol in $u_1^{e}(K)$, while the second and third factors on the right give symbols in $u_1^{>e}(K)$ since 
\[v(p(\zeta_p-1)^2)=v(p)+2v(\zeta_p-1)=\frac{p+1}{p-1}v(p)>e.\]
Thus, the first term on the right also yields a symbol in $u_1^{>e}(K)$, as wished. Since $H_p^{q+1}=(\Omega_k^q/\im(d))/\wp(\Omega_k^q)$, we deduce from all the above that we have a well-defined surjective homomorphism
\begin{align*}
        \rho_q^e: \bigoplus_{j=0}^{q-1} \left(\Lambda^j(\Gamma/p\Gamma) \otimes H_p^{q-j}(k)\right) & \rightarrow u_q^{e}(K)/u_q^{>e}(K)\\
        (\bar{\gamma}_1 \wedge \dots \wedge \bar{\gamma}_j)\otimes\left\{\bar{x} \frac{d \bar{y}_{j+2}}{\bar{y}_{j+2}}\wedge \dots \wedge\frac{d \bar{y}_{q}}{\bar{y}_{q}}\right\} & \mapsto \{\pi^{\bar{\gamma}_1},\dots ,\pi^{\bar{\gamma}_j}, 1+(\zeta_p-1)^p x, y_{j+2},\dots,y_{q}\};
\end{align*}

This concludes the proof of Step 1. Before moving forward, let us state the following consequence of this step for further use.

\begin{lemma}\label{lem trivialite bcp de slots}
Write $[k:k^p]=p^r$. Let $\gamma=p\gamma'\in \Gamma$ and let
\[\alpha=\{\pi^{\bar{\gamma}_1},\dots ,\pi^{\bar{\gamma}_j}, 1+(\pi^{\gamma'})^px, y_{j+2},\dots,y_{q}\}\in u_{q}^\gamma(K).\]
Assume that either $j>\dim_{\F_p}(\Gamma/p\Gamma)$ or $q-j>r$. Then $\alpha=0$ in $u_q^{\gamma}(K)/u_q^{>\gamma}(K)$.
\end{lemma}

\begin{proof}
Since the morphism $\rho_q^{\gamma}$ is well-defined for $\gamma=p\gamma'$, we know that
\[\alpha=\rho_q^{\gamma}\left((\bar{\gamma}_1 \wedge \dots \wedge \bar{\gamma}_j)\otimes\left\{\bar{x} \frac{d \bar{y}_{j+2}}{\bar{y}_{j+2}}\wedge \dots \wedge\frac{d \bar{y}_{q}}{\bar{y}_{q}}\right\}\right).\]
But the hypothesis implies that one of the two factors in the pure tensor is trivial (recall that the second one lives in $\Omega_k^{q-j-1}/\ker(d)$), which concludes the proof.
\end{proof}

\subsection{Step 2: the case $\gamma=0$}

Let us prove that $\rho_q^0$ is an isomorphism. Since we already know that it is well-defined and surjective, it will suffice to present an explicit inverse. To do so, we introduce the map:
\begin{align*}
        \theta_q^0: {K^\times}^q & \rightarrow \bigoplus_{j=0}^{q} \left(\Lambda^j(\Gamma/p\Gamma) \otimes \kk_{q-j}(k)\right)\\
        (\pi^{\bar{\gamma}_1}a_1^pu_1,\dots,\pi^{\bar{\gamma}_q}a_q^pu_q) & \mapsto \sum_{I \sqcup J =  \{1,\dots, q\}} (-1)^{\epsilon(I,J)} ({\bar{\gamma}_{i_1}}\wedge\dots \wedge \bar{\gamma}_{i_r}) \otimes \{\overline{u}_{j_1},\dots,\overline{u}_{j_{q-r}}\},
         \end{align*}
         where $\bar{\gamma}_1,\dots,\bar{\gamma}_q \in \Gamma/p\Gamma$, $a_1,\dots,a_q \in K^\times$, $u_1,\dots,u_q \in \mathcal{O}_K^\times$, $I=\{i_1,\dots,i_r\}$ with $i_1<\dots<i_r$, $J=\{j_1,\dots,j_{q-r}\}$ with $j_1<\dots<j_{q-r}$ and $\epsilon(I,J)$ is the signature of the permutation in $\mathcal{S}_q$ that sends each $i_s$ to $s$ and each $j_s$ to $r+s$. One easily checks that it is well-defined and multilinear, and that $\theta_q^0(y_1,\dots,y_q)=0$ if at least one of the $y_i$'s is a $p$-th power. \\
         
         Let us now check that $\theta_q^0(y_1,\dots,y_q)=0$ if $y_i+y_j=1$ for some $i\neq j$. Since for each $\sigma \in \mathcal{S}_q$ and each $(y_1,\dots,y_q)\in {K^\times}^q$ we have $\theta^0_q(y_{\sigma(1)}, \dots, y_{\sigma(q)})=\epsilon(\sigma)\theta^0_q(y_1,\dots,y_q)$ where $\epsilon(\sigma)$ stands for the signature of $\sigma$, it suffices to check that $\theta_q^0(y_1,\dots,y_q)=0$ if $y_1+y_2=1$ and $v(y_1)\geq v(y_2)$. For each $i$, set $\gamma_i:=v(y_i)$ and write $y_i=\pi^{\bar{\gamma}_i}a_i^pu_i$ for some $a_i\in K^\times$ and $u_i\in \mathcal{O}_K^\times$. Three main cases arise:
         \begin{enumerate}
             \item Case 1: $\gamma_1>0$. In that case, $\gamma_2=0$ and we may assume that $a_2=1$ and $y_2=u_2$. We automatically have $\overline{u}_2=1$. Then:
             \begin{align*}
                 \theta^0_q(y_1,\dots,y_q)= &\sum_{\substack{I \sqcup J =  \{1,\dots, q\}\\ \{1,2\}\subset I}} (-1)^{\epsilon(I,J)} ({\bar{\gamma}_{1}}\wedge 0\wedge\bar{\gamma}_{i_3}\wedge\dots \wedge{\bar{\gamma}_{i_r}})\otimes \{\overline{u}_{j_1},\dots,\overline{u}_{j_{q-r}}\}\\
                 &+ \sum_{\substack{I \sqcup J =  \{1,\dots, q\}\\ \{1,2\}\subset J}} (-1)^{\epsilon(I,J)} (\bar{\gamma}_{i_1}\wedge\dots \wedge{\bar{\gamma}_{i_r}}) \otimes \{\overline{u}_{1},1,\overline{u}_{j_3}\dots,\overline{u}_{j_{q-r}}\}\\
                 &+ \sum_{\substack{I \sqcup J =  \{1,\dots, q\}\\ 1\in I, 2 \in J}} (-1)^{\epsilon(I,J)} (\bar{\gamma}_{i_1}\wedge\dots \wedge{\bar{\gamma}_{i_r}}) \otimes \{1, \overline{u}_{j_2},\dots,\overline{u}_{j_{q-r}}\}\\
                 &+ \sum_{\substack{I \sqcup J =  \{1,\dots, q\}\\ 2\in I, 1 \in J}} (-1)^{\epsilon(I,J)} (0\wedge\bar{\gamma}_{i_2}\wedge\dots \wedge{\bar{\gamma}_{i_r}})\otimes \{\overline{u}_{j_1},\dots,\overline{u}_{j_{q-r}}\}\\
                 = &\;0.
             \end{align*}
             \item Case 2: $\gamma_1<0$. In that case, $\gamma_2=\gamma_1<0$ and there exists $\overline{b}\in k^\times$ such that $\overline{u}_2=-\overline{b}^p\overline{u}_1$. Then:
             \begin{align*}
                 \theta^0_q(y_1,\dots,y_q)= &\sum_{\substack{I \sqcup J =  \{1,\dots, q\}\\ \{1,2\}\subset I}} (-1)^{\epsilon(I,J)} ({\bar{\gamma}_{1}}\wedge\bar{\gamma}_{1}\wedge\bar{\gamma}_{i_3}\wedge\dots \wedge{\bar{\gamma}_{i_r}}) \otimes \{\overline{u}_{j_1},\dots,\overline{u}_{j_{q-r}}\}\\
                 &+ \sum_{\substack{I \sqcup J =  \{1,\dots, q\}\\ \{1,2\}\subset J}} (-1)^{\epsilon(I,J)} (\bar{\gamma}_{i_1}\wedge\dots \wedge{\bar{\gamma}_{i_r}}) \otimes \{\overline{u}_{1},-\overline{b}^p\overline{u}_{1},\overline{u}_{j_3}\dots,\overline{u}_{j_{q-r}}\}\\
                 &+ \sum_{\substack{I \sqcup J =  \{1,\dots, q\}\\ 1\in I, 2 \in J}} (-1)^{\epsilon(I,J)} (\bar{\gamma}_{1}\wedge\bar{\gamma}_{i_2}\wedge\dots \wedge{\bar{\gamma}_{i_r}}) \otimes \{\overline{u}_{2}, \overline{u}_{j_2},\dots,\overline{u}_{j_{q-r}}\}\\
                 &+ \sum_{\substack{I \sqcup J =  \{1,\dots, q\}\\ 2\in I, 1 \in J}} (-1)^{\epsilon(I,J)} (\bar{\gamma}_{1}\wedge\bar{\gamma}_{i_2}\wedge\dots \wedge{\bar{\gamma}_{i_r}}) \otimes \{\overline{u}_{1}, \overline{u}_{j_2},\dots,\overline{u}_{j_{q-r}}\}\\
                 =&\sum_{\substack{I \sqcup J =  \{1,\dots, q\}\\ 1\in I, 2 \in J}} (-1)^{\epsilon(I,J)} (\bar{\gamma}_{1}\wedge\bar{\gamma}_{i_2}\wedge\dots \wedge{\bar{\gamma}_{i_r}}) \otimes \{\overline{u}_{1}^{-1}\overline{u}_2, \overline{u}_{j_2},\dots,\overline{u}_{j_{q-r}}\}\\
                 =&\sum_{\substack{I \sqcup J =  \{1,\dots, q\}\\ 1\in I, 2 \in J}} (-1)^{\epsilon(I,J)} (\bar{\gamma}_{1}\wedge\bar{\gamma}_{i_2}\wedge\dots \wedge{\bar{\gamma}_{i_r}}) \otimes \{-\overline{b}^{p}, \overline{u}_{j_2},\dots,\overline{u}_{j_{q-r}}\}\\
                 = &\;0.
             \end{align*}
             \item Case 3: $\gamma_1=0$. In that case, $\gamma_1=\gamma_2=0$ and we may assume that $a_1=a_2=1$. We automatically have $\overline{u}_1+\overline{u}_2=1$. Then:
             \begin{align*}
                 \theta^0_q(y_1,\dots,y_q)= &\sum_{\substack{I \sqcup J =  \{1,\dots, q\}\\ I \cap  \{1,2\}\neq \emptyset}} (-1)^{\epsilon(I,J)} (0\wedge\bar{\gamma}_{i_2}\wedge\bar{\gamma}_{i_3}\wedge\dots \wedge{\bar{\gamma}_{i_r}}) \otimes \{\overline{u}_{j_1},\dots,\overline{u}_{j_{q-r}}\}\\
                 &+ \sum_{\substack{I \sqcup J =  \{1,\dots, q\}\\ \{1,2\}\subset J}} (-1)^{\epsilon(I,J)} (\bar{\gamma}_{i_1}\wedge\dots \wedge{\bar{\gamma}_{i_r}}) \otimes \{\overline{u}_{1},1-\overline{u}_{1},\overline{u}_{j_3}\dots,\overline{u}_{j_{q-r}}\}\\
                 = &\;0.
             \end{align*}
         \end{enumerate}
         
We deduce that $\theta_q^0$ induces a morphism:
$$\kk_q(K)  \rightarrow \bigoplus_{j=0}^{q} \left(\Lambda^j(\Gamma/p\Gamma) \otimes \kk_{q-j}(k)\right).$$
One easily checks using Lemma \ref{lem generateurs u_q} that the kernel contains $u_q^{>0}(K)$, and hence we get a morphism:
$$\kk_q(K)/u_q^{>0}(K) \rightarrow \bigoplus_{j=0}^{q} \left(\Lambda^j(\Gamma/p\Gamma) \otimes \kk_{q-j}(k)\right),$$
that we still denote $\theta_q^0$. A simple computation shows that $\theta_q^0$ and $\rho_q^0$ are inverse of each other.

We have thus proved item (i) of Theorem \ref{thm isom Kato}, which concludes Step 2.

\subsection{Step 3: nontriviality of $\rho_q^\gamma$ for $\gamma\in p\Gamma$}

The proof of injectivity of $\rho_q^\gamma$ for $\gamma<e$ in the classical case uses a duality argument that needs the nontriviality of $\rho_q^e$ under certain hypotheses. While we follow this approach, we can actually simplify the argument by using instead $\rho_q^{\gamma_\infty}$ for a suitable $\gamma_\infty$ close enough to $e$ (note that this makes no sense in the discrete case). Moreover, this allows us to elude the case where the domain of $\rho_q^{e}$ is trivial. We prove then an (easier) nontriviality result, which is the following.

\begin{proposition}\label{lem triv symb}
Let $(\bar{\gamma}_1,\dots ,\bar{\gamma}_m)$ be a free family in the $\mathbb{F}_p$-vector space $\Gamma/p\Gamma$ and $(\bar x,\bar y_1,\dots,\bar y_{n-1})$ be a $p$-free family in $k$. Then the symbol
\[\{\pi^{\bar{\gamma}_1},\dots ,\pi^{\bar{\gamma}_m}, 1+(\pi^{\gamma'})^px,y_1,\dots,y_{n-1}\}\in \kk_{m+n}(K),\]
is nontrivial if $0<p\gamma'<e$.
\end{proposition}

\begin{proof}
We first deal with the case $m=0$ by induction on $n$. Assume first that $n=1$. If $1+(\pi^{\gamma'})^px=z^p$ for some $z\in K$, then, using the fact that $p\gamma'<e$, one sees that $z$ has the form $1+\pi^{\gamma'}y$ with $\bar y=\sqrt[p]{\bar x}$. But this contradicts that $\bar x\not\in k^p$, and hence $1+(\pi^{\gamma'})^px\not\in K^p$.

For the induction step, we use the norm varieties of Suslin and Joukhovitski \cite{SJ}. Denote by $K_{(p)}$ a maximal prime-to-$p$ extension of $K$ (that is, the fixed field of a pro-$p$-Sylow subgroup of $\Gal(K)$), by $\Gamma_{K_{(p)}}$ its value group and by $k_{(p)}$ its residue field. We are going to prove the nontriviality of the symbol $\{1+(\pi^{\gamma'})^px,y_1,\dots,y_{n-1}\}$ in $\kk_{n}(K_{(p)})$, which implies its nontriviality in $\kk_{n}(K)$. By \cite[Thm.~A.1]{SJ}, we know that this symbol is nontrivial in $\kk_{n}(K_{(p)})$ if and only if $y_{n-1}$ is not a product of norms coming from finite extensions  of $K_{(p)}$ trivializing $\{1+(\pi^{\gamma'})^px,y_1,\dots,y_{n-2}\}$.

Fix then a finite extension $L$ of $K_{(p)}$ trivializing $\{1+(\pi^{\gamma'})^px,y_1,\dots,y_{n-2}\}$  and denote by $l$ its residue field. Define
$$U(\bar x,\bar y_1,\dots ,\bar y_{n-2}):=\{a\in \mathcal{O}_{K_{(p)}}^\times: \bar a \in k^p_{(p)}(\bar x,\bar y_1,\dots ,\bar y_{n-2})\}.$$
We claim that $N_{L/K_{(p)}}(L^\times)$ is contained in $(K_{(p)}^\times)^p\cdot U(\bar x,\bar y_1,\dots ,\bar y_{n-2})$. Since $y_{n-1}$ clearly cannot be written as a product of elements in $(K_{(p)}^\times)^p\cdot U(\bar x,\bar y_1,\dots ,\bar y_{n-2})$, this yields the desired result.\\

Let us prove the claim for $L$. Since the Galois closure of $L/K_{(p)}$ has degree a power of $p$, we can find intermediary field extensions $K_{(p)}=L_0\subset L_1 \subset \dots \subset L_r=L$ such that $[L_{i+1}:L_i]=p$ for every $0\leq i\leq r-1$. We denote by $k_{(p)}=l_0\subset l_1\subset \dots \subset l_r=l$ the corresponding residue fields.

By assumption, we know that $(\bar x, \bar y_1, \dots, \bar y_{n-2},\bar y_{n-1})$ is $p$-free in $k_{(p)}$, while by the inductive assumption, we know that $(\bar x, \bar y_1, \dots, \bar y_{n-2})$ is not $p$-free in $l$. Let $s$ be the smallest integer such that $(\bar x, \bar y_1, \dots, \bar y_{n-2})$ is $p$-free in $l_{s}$ but not in $l_{s+1}$. Given $z\in L_{s+1}^\times$, we aim at proving that $N_{L_{s+1}/K_{(p)}}(z)\in (K_{(p)}^\times)^p\cdot U(\bar x,\bar y_1,\dots ,\bar y_{n-2})$. Since $[l_{s+1}:l_s]=[L_{s+1}:L_s]=p$, we can write $z=\varpi u$ for some $\varpi\in L_s$ and some $u\in \mathcal{O}_{L_{s+1}}^\times$. Then $N_{L_{s+1}/K_{(p)}}(z)= N_{L_s/K_{(p)}}(\varpi)^p\cdot N_{L_{s+1}/K_{(p)}}(u)$. It therefore suffices to prove that $N_{L_{s+1}/K_{(p)}}(u)\in (K_{(p)}^\times)^p\cdot U(\bar x,\bar y_1,\dots ,\bar y_{n-2})$.

For that purpose, we prove for each $i \leq s$ that: 
$$N_{L_{s+1}/L_i}(u)\in U_{L_i}(\bar x,\bar y_1,\dots ,\bar y_{n-2}):=\{a\in \mathcal{O}_{L_i}^\times: \bar a \in l^p_i(\bar x,\bar y_1,\dots ,\bar y_{n-2})\}.$$ We proceed by descending induction on $i$. We start with the case $i=s$. Since $[l_{s+1}:l_s]=[L_{s+1}:L_s]=p$, we have $\overline{N_{L_{s+1}/L_s}(u)}=N_{l_{s+1}/l_s}(\bar u)=\bar u^p$. By definition of $s$, the family $(\bar x, \bar y_1, \dots, \bar y_{n-2})$ is $p$-free in $l_{s}$ but not in $l_{s+1}$. This implies that $l_{s+1}=l_s(\sqrt[p]{\rho})$ for some $\rho\in l_s^p(\bar x, \bar y_1, \dots, \bar y_{n-2}) \setminus l_s^p$. 
Hence $\bar u^p \in l_s^p(\rho)\subset l_s^p(\bar x, \bar y_1, \dots, \bar y_{n-2})$, as wished.

We now consider $i<s$. Three cases arise:
\begin{enumerate}
    \item \emph{$l_{i+1}/l_i$ is purely inseparable of degree $p$.} In this case, we have: $$\overline{N_{L_{s+1}/L_i}(u)}=\overline{N_{L_{i+1}/L_{i}}(N_{L_{s+1}/L_{i+1}}(u))}=N_{l_{i+1}/l_{i}}(\overline{N_{L_{s+1}/L_{i+1}}(u)})=\overline{N_{L_{s+1}/L_{i+1}}(u)}^p.$$ 
    But by the inductive assumption, $N_{L_{s+1}/L_{i+1}}(u)\in  U_{L_{i+1}}(\bar x,\bar y_1,\dots ,\bar y_{n-2})$. Hence $\overline{N_{L_{s+1}/L_{i+1}}(u)}\in l_{i+1}^p(\bar x, \bar y_1,\dots,\bar y_{n-2})\subset l_i$, so that $\overline{N_{L_{s+1}/L_{i+1}}(u)}^p\in l_{i}^p$.
    \item \emph{$l_{i+1}/l_i$ is separable of degree $p$.} In this case, we have: $$\overline{N_{L_{s+1}/L_i}(u)}=\overline{N_{L_{i+1}/L_{i}}(N_{L_{s+1}/L_{i+1}}(u))}=N_{l_{i+1}/l_{i}}(\overline{N_{L_{s+1}/L_{i+1}}(u)}).$$ 
    But by the inductive assumption, $N_{L_{s+1}/L_{i+1}}(u)\in  U_{L_{i+1}}(\bar x,\bar y_1,\dots ,\bar y_{n-2})$. Hence $\overline{N_{L_{s+1}/L_{i+1}}(u)}\in l_{i+1}^p(\bar x, \bar y_1,\dots,\bar y_{n-2})$, so that: 
    \begin{align*}
    N_{l_{i+1}/l_{i}}(\overline{N_{L_{s+1}/L_{i+1}}(u)})&=\prod_{\sigma\in \mathrm{Gal}(l_{i+1}/l_i)} \sigma (\overline{N_{L_{s+1}/L_{i+1}}(u)})\\&\in (l_{i+1}^p \cap l_i)(\bar x, \bar y_1,\dots,\bar y_{n-2})=l_i^p(\bar x, \bar y_1,\dots,\bar y_{n-2}). 
    \end{align*}
    \item \emph{$l_{i+1}=l_i$.} In this case, we have: $$\overline{N_{L_{s+1}/L_i}(u)}=\overline{N_{L_{i+1}/L_{i}}(N_{L_{s+1}/L_{i+1}}(u))}=\overline{N_{L_{s+1}/L_{i+1}}(u)}^p\in l_{i+1}^p=l_i^p.$$ 
\end{enumerate}
We deduce that $N_{L_{s+1}/K_{(p)}}(u)\in U(\bar x,\bar y_1,\dots ,\bar y_{n-2})$, and hence that $N_{L/K_{(p)}}(z)\in (K_{(p)}^{\times})^p \cdot U(\bar x, \bar y_1,\dots,\bar y_{n-2})$, as wished. This finishes the proof in the case $m=0$.\\

We now deal with the case $m\geq 1$, working again by induction, over $K_{(p)}$ and using \cite{SJ}. Let us prove then that $\pi^{\bar\gamma_m}$ is not a product of norms coming from finite extensions of $K_{(p)}$ trivializing $\{\pi^{\bar\gamma_1},\dots,\pi^{\bar \gamma_{m-1}},1+(\pi^{\gamma'})^px,y_1,\dots,y_{n-1}\}$. This will imply that the symbol $\{\pi^{\bar\gamma_1},\dots,\pi^{\bar \gamma_{m}},1+(\pi^{\gamma'})^px,y_1,\dots,y_{n-1}\}$ is nontrivial in $\kk_{m+n}(K_{(p)})$ by \cite[Thm.~A.1]{SJ}.

To prove this, fix a finite extension $L$ of $K_{(p)}$ trivializing $\{\pi^{\bar\gamma_1},\dots,\pi^{\bar \gamma_{m-1}},1+(\pi^{\gamma'})^px,y_1,\dots,y_{n-1}\}$, denote by $\Gamma_L$ its value group and by $l$ its residue field. By the inductive assumption, either $(\bar \gamma_1,\dots,\bar \gamma_{m-1})$ is not $\mathbb{F}_p$-free in $\Gamma_L/p\Gamma_L$ or $(\bar x, \bar y_1,\bar y_2,\dots,\bar y_{n-1})$ is not $p$-free in $l$.
\begin{enumerate}
    \item[$\bullet$] In the first case, there exists some $\mu' \in \Gamma_{K_{(p)}}$ such that $\bar\mu' \in \langle \bar \gamma_1,\dots,\bar \gamma_{m-1} \rangle \setminus \{0\}\subset \Gamma_{K_{(p)}}/p\Gamma_{K_{(p)}}$ and $\mu'\in p\Gamma_L$. Take some element $y \in L^\times$ and set  $\mu:=v(N_{L/K_{(p)}}(y))\in \Gamma_{K_{(p)}}$, so that $v(y)=\frac{\mu}{[L:K_{(p)}]}\in \Gamma_L$. If $\bar \mu\neq 0\in\Gamma_{K_{(p)}}/p\Gamma_{K_{(p)}}$, then $[\langle \Gamma_{K_{(p)}}, \frac{\mu}{[L:K_{(p)}]}\rangle:\Gamma_{K_{(p)}}]=[L:K_{(p)}]$ since the extension is of $p$-primary degree, and hence $\Gamma_L=\langle \Gamma_{K_{(p)}}, \frac{\mu}{[L:K_{(p)}]}\rangle \ni \frac{\mu'}{p}$. This implies that $\mu'\equiv \frac{\mu}{p^\ell} \pmod{p\Gamma_{K_{(p)}}}$ for some $\ell\geq 0$ and hence $\bar \mu \in \langle \bar\gamma_1,\dots,\bar\gamma_{m-1}\rangle$.
\item[$\bullet$] In the second case, we have $[\Gamma_L:\Gamma_{K_{(p)}}] < [L:K_{(p)}]$ and hence $p | \frac{[L:K_{(p)}]}{[\Gamma_L:\Gamma_{K_{(p)}}]} | v(N_{L/K_{(p)}}(y))$ for all $y\in L^\times$. 
\end{enumerate}
In both cases, we get that $v(N_{L/K_{(p)}}(L^\times))\subset \langle \gamma_1,\dots,\gamma_{m-1}\rangle + p\Gamma_{K_{(p)}}$. Then $\pi^{\bar\gamma_m}$ cannot be a product of such norms, proving the claim and concluding the proof of the proposition.
\end{proof}

\begin{corollary}\label{cor non trivial}
    Let $(\bar{\gamma}_1,\dots ,\bar{\gamma}_m)$ be a free family in the $\mathbb{F}_p$-vector space $\Gamma/p\Gamma$ and $(\bar x,\bar y_1,\dots,\bar y_{n-1})$ be a $p$-free family in $k$. Let $\gamma'\in \Gamma$ with $e/2<p\gamma'< e$. Then the symbol
\[\alpha:=\{\pi^{\bar{\gamma}_1},\dots ,\pi^{\bar{\gamma}_m}, 1+(\pi^{\gamma'})^px,y_1,\dots,y_{n-1}\}\in \kk_{m+n}(K),\]
does not belong to $u_{m+n}^{>p\gamma'}(K)$.
\end{corollary}

\begin{proof}
By contradiction, assume that $\alpha\in u_{m+n}^{>p\gamma'}(K)$. Since $p\Gamma$ is dense in $\R$, there exists $\gamma\in p\Gamma$ with $\gamma>p\gamma'$ such that $\alpha\in u_{m+n}^{\gamma}(K)$. Then, by Lemma \ref{lem generateurs u_q}, we can write:
$$\alpha=\sum_{i=1}^\ell \{ 1+(\pi^{\gamma/p})^p x_i,y_{i,1},\dots,y_{i,m+n-1}\}$$
for some $y_{i,j}\in \mathcal{O}_K\setminus \{0\}$ and $x_i \in \mathcal{O}_K$.\\

Consider the field $L=K(t)$ equipped with the Gauss extension of $v$, which has the same value group $\Gamma$ (cf.~\cite[Cor.~2.2.2]{EnglerPrestel}), and fix some $\gamma_{\infty}\in p\Gamma$ such that: $$\max\{(p\gamma'+e)/2,e+p\gamma'-\gamma\}<\gamma_{\infty}<e.$$
We then have:
$$\alpha\cup \{1+(\pi^{\gamma_\infty/p-\gamma'})^pt\}=\sum_{i=1}^\ell \{ 1+(\pi^{\gamma/p})^px_i,y_{i,1},\dots,y_{i,m+n-1},1+(\pi^{\gamma_\infty/p-\gamma'})^pt\}\in \kk_{m+n+1}(K(t)).$$
Since $\gamma+\gamma_{\infty}-p\gamma'>e$, the second assertion of Lemma \ref{lem suma} shows that the right-hand side of the equality is zero. Hence:
$$\{\pi^{\bar{\gamma}_1},\dots ,\pi^{\bar{\gamma}_m}, 1+(\pi^{\gamma'})^px,y_1,\dots,y_{n-1},1+(\pi^{\gamma_\infty/p-\gamma'})^pt\}=0.$$
The last assertion of Lemma \ref{lem suma} allows us to deduce:
$$\{\pi^{\bar{\gamma}_1},\dots ,\pi^{\bar{\gamma}_m}, 1-(\pi^{\gamma'}\pi^{\gamma_\infty/p-\gamma'})^pxt,y_1,\dots,y_{n-1},-x\}=0,$$
and hence
$$\{\pi^{\bar{\gamma}_1},\dots ,\pi^{\bar{\gamma}_m}, 1-(\pi^{\gamma_\infty/p})^pu^pxt,y_1,\dots,y_{n-1},-x\}=0,$$
for some $u\in\cal O_K^\times$. Since $(\bar u ^p\bar{x}\bar{t},\bar y_1,\dots,\bar y_{n-1},\bar x)$ is a $p$-free family in the residue field $k(t)$ of $K(t)$, the last equality contradicts Proposition \ref{lem triv symb}.
\end{proof}

\subsection{Step 4: injectivity of $\rho_q^\gamma$ for $0<\gamma<e$}

Recall from Section \ref{subsec setting} that $p^r=[k:k^p]$ and $g=\dim_{\F_p}(\Gamma/p\Gamma)$. Take some $\gamma \in \Gamma$ with $0< \gamma < e$ and choose a decomposition $\Gamma/p\Gamma=\langle \bar{\gamma} \rangle \oplus \Delta$. We aim at proving that the morphism $\rho_q^\gamma$
is injective. To do so, we fix once and for all some $\gamma_\infty\in p\Gamma$ such that:
$$\max \left\{\frac{e+\gamma}{2},e-\gamma\right\}<\gamma_{\infty}<e.$$
We also introduce the following notation:

\begin{notation}
    For $\gamma_1,\gamma_2\in\Gamma$, we let $u_{\gamma_1,\gamma_2}\in\mathcal{O}_K^{\times}$ be such that $\pi^{\gamma_1}\pi^{\gamma_2}=u_{\gamma_1,\gamma_2}\pi^{\gamma_1+\gamma_2}$ if $\gamma_1+\gamma_2\neq e$ and $\pi^{\gamma_1}\pi^{\gamma_2}=u_{\gamma_1,\gamma_2}(\zeta_p-1)^p$ if $\gamma_1+\gamma_2=e$. If moreover $\gamma=p\gamma'\in p\Gamma$, we let $u'_\gamma\in\mathcal{O}_K^{\times}$ be such that $(\pi^{\gamma'})^p = u'_\gamma\pi^\gamma$.
\end{notation}

\begin{lemma}\label{lem calcul}
\begin{enumerate}
    \item Assume that $p\nmid \gamma$. We have a commutative diagram:
\[\xymatrix{
 \left(\displaystyle\bigoplus_{j}\Lambda^{j}\Delta \otimes \Omega_k^{q-j-1}\right) \times \left( \displaystyle\bigoplus_{j} \Lambda^{g-j-1}\Delta \otimes \Omega_k^{r+j-q}\right) \ar[d]^{\rho_q^{\gamma} \times \rho_{r+g-q}^{\gamma_\infty-\gamma}} \ar[r]^-{\varphi} & \Lambda^g(\Gamma/p\Gamma) \otimes \frac{\Omega_k^{r-1}}{\ker(d)} \ar[d]^{\rho^{\gamma_{\infty}}_{r+g}}\\
u_q^{\gamma}(K)/u_q^{>\gamma}(K)\times  u_{r+g-q}^{\gamma_{\infty}-\gamma}(K)/u_{r+g-q}^{>\gamma_{\infty}-\gamma}(K) \ar[r]^-{\cup} & u_{r+g}^{\gamma_\infty}(K)/u_{r+g}^{>\gamma_\infty}(K)
    }\]
where $\varphi$ sends $((\delta_j\otimes\omega_j)_j,(\delta'_j\otimes\omega'_j)_j)\in  \left(\bigoplus_{j} \Lambda^{j}\Delta \otimes \Omega_k^{q-j-1}\right) \times \left( \bigoplus_{j} \Lambda^{g-j-1}\Delta \otimes \Omega_k^{r+j-q}\right)$ to:
$$
    (\bar u'_{\gamma_\infty})^{-1}\bar u_{\gamma,\gamma_\infty-\gamma}\left(\sum_{j=0}^{q-1} \epsilon_j (\delta_j\wedge\delta'_j\wedge \overline{\gamma})\otimes(\omega_j\wedge\omega'_{j})+ \sum_{j=1}^{q-1} \eta_j (\delta_j\wedge\delta'_{j-1})\otimes(d\omega_j\wedge\omega'_{j-1})\right)
$$
and $\epsilon_j,\eta_j\in\{-1,1\}$ are signs that only depend on $j, q, g$.
\item Assume that $p\mid \gamma$.  We have a commutative diagram:
\[\xymatrix{
 \left(\displaystyle\bigoplus_{j} \Lambda^j(\Gamma/p\Gamma) \otimes \frac{\Omega_k^{q-j-1}}{\ker(d)}\right) \times \left(\displaystyle\bigoplus_{j} \Lambda^{g-j}(\Gamma/p\Gamma) \otimes \frac{\Omega_k^{r+j-q-1}}{\ker(d)}\right) \ar[d]^{\rho_q^{\gamma} \times \rho_{r+g-q}^{\gamma_\infty-\gamma}} \ar[r]^-{\varphi} & \Lambda^g(\Gamma/p\Gamma) \otimes \frac{\Omega^{r-1}_k}{\ker (d)}\ar[d]^{\rho^{\gamma_\infty}_{r+g}}\\
u_q^{\gamma}(K)/u_q^{>\gamma}(K)\times  u_{r+g-q}^{\gamma_\infty-\gamma}(K)/u_{r+g-q}^{>\gamma_\infty-\gamma}(K) \ar[r]^-{\cup} & u_{r+g}^{\gamma_\infty}(K)/u_{r+g}^{>\gamma_\infty}(K)
    }\]
where $\varphi$ sends $$((\delta_j\otimes\omega_j)_j,(\delta'_j\otimes\omega'_j)_j)\in \left(\bigoplus_{j} \Lambda^j(\Gamma/p\Gamma) \otimes \Omega_k^{q-j-1}\right) \times \left( \bigoplus_{j} \Lambda^{g-j}(\Gamma/p\Gamma) \otimes \Omega_k^{r+j-q-1}\right)$$ to:
$$
    \bar u_{\gamma,\gamma_\infty-\gamma} \bar u'_{\gamma}\bar u'_{\gamma_\infty-\gamma} (\bar u'_{\gamma_{\infty}})^{-1} \sum_{j=0}^{q-1} \epsilon_j (\delta_j\wedge\delta'_{j})\otimes(d\omega_j\wedge\omega'_{j})  
$$
and $\epsilon_j\in\{-1,1\}$ are signs that only depend on $j, q, g$.
\end{enumerate}
\end{lemma}

\begin{proof}
\begin{enumerate}
    \item Choose $j,j'\geq 0$. Take: $$\begin{cases}
        (\bar{\mu}_1 \wedge \dots \wedge \bar{\mu}_j)\otimes\left\{\bar{x} \frac{d \bar{y}_{j+2}}{\bar{y}_{j+2}}\wedge \dots \wedge\frac{d \bar{y}_{q}}{\bar{y}_{q}}\right\}\in \Lambda^{j}\Delta \otimes \Omega_k^{q-j-1}\\
        (\bar{\nu}_1 \wedge \dots \wedge \bar{\nu}_{g-j'-1})\otimes\left\{\bar{z} \frac{d \bar{t}_{g-j'+1}}{\bar{t}_{g-j'+1}}\wedge \dots \wedge\frac{d \bar{t}_{r+g-q}}{\bar{t}_{r+g-q}}\right\}\in\Lambda^{g-j'-1}\Delta \otimes \Omega_k^{r+j'-q},
    \end{cases}$$
    and set for simplicity: 
    \begin{gather*}
    \bar{\mu}:=\bar{\mu}_1 \wedge \dots \wedge \bar{\mu}_j,\;\;\;\;\; \bar{\nu}:=\bar{\nu}_1 \wedge \dots \wedge \bar{\nu}_{g-j'-1},\\ 
    \pi^{\bar{\mu}}:=(\pi^{\bar\mu_1},\dots , \pi^{\bar\mu_j}),\;\;\;\;\; \pi^{\bar{\nu}}:=(\pi^{\bar\nu_1},\dots,\pi^{\bar\nu_{g-j'-1}}),\\
        \mathbf{y}:=(y_{j+2},\dots,y_q),\;\;\;\;\; \mathbf{t}:=(t_{g-j'+1},\dots,t_{r+g-q}),\\ \frac{d\bar{\mathbf{y}}}{\bar{\mathbf{y}}}:=\frac{d \bar{y}_{j+2}}{\bar{y}_{j+2}}\wedge \dots \wedge\frac{d \bar{y}_{q}}{\bar{y}_{q}},\;\;\;\;\; \frac{d\bar{\mathbf{t}}}{\bar{\mathbf{t}}}:=\frac{d \bar{t}_{g-j'+1}}{\bar{t}_{g-j'+1}}\wedge \dots \wedge\frac{d \bar{t}_{r+g-q}}{\bar{t}_{r+g-q}}.
    \end{gather*}
    We compute:
    \begin{align*}
         \rho_q^\gamma &\left(\bar{\mu} \otimes\left\{\bar{x}  \frac{d\bar{\mathbf{y}}}{\bar{\mathbf{y}}}\right\}\right) \cup \rho_{r+g-q}^{\gamma_\infty-\gamma}\left(\bar{\nu}\otimes\left\{\bar{z}  \frac{d\bar{\mathbf{t}}}{\bar{\mathbf{t}}}\right\}\right)\\
         & =\pm \left\{ \pi^{\bar\mu},\pi^{\bar\nu},1+\pi^{\gamma}x,1+\pi^{\gamma_\infty-\gamma}z,\mathbf{y},\mathbf{t}\right\}\\
         &=
         \begin{cases}
              \pm \left\{ \pi^{\bar\mu},\pi^{\bar\nu},1- \pi^{\gamma}\pi^{\gamma_\infty-\gamma}xz,-\pi^\gamma x,\mathbf{y},\mathbf{t}\right\} &\text{if $\gamma \geq \gamma_\infty-\gamma$}  \\
               \pm \left\{ \pi^{\bar\mu},\pi^{\bar\nu},1- \pi^{\gamma}\pi^{\gamma_\infty-\gamma}xz,-\pi^{\gamma_\infty-\gamma} z,\mathbf{y},\mathbf{t}\right\} &\text{if $\gamma < \gamma_\infty-\gamma$}
         \end{cases}  \\
          \intertext{\raggedleft (by Lemma \ref{lem suma} and the inequality $\max \{\frac{e+\gamma}{2},e-\gamma\}<\gamma_{\infty}<e$)  }
         &=
         \begin{cases}
             \pm \rho^{\gamma_\infty}_{r+g}\left( (\bar{\mu} \wedge \bar{\nu} \wedge\bar\gamma)\otimes \left(\bar u\bar x\bar z \frac{d\bar{\mathbf{y}}}{\bar{\mathbf{y}}} \wedge\frac{d\bar{\mathbf{t}}}{\bar{\mathbf{t}}}\right)\right)
             &\text{if $j=j'$}\\
              \pm \rho^{\gamma_\infty}_{r+g}\left( (\bar{\mu} \wedge \bar{\nu})\otimes \left(\bar u\bar x\bar z\frac{d \bar{x}}{\bar{x}}\wedge \frac{d\bar{\mathbf{y}}}{\bar{\mathbf{y}}} \wedge\frac{d\bar{\mathbf{t}}}{\bar{\mathbf{t}}}\right)\right) &\text{if $j=j'+1$ and $\gamma \geq \gamma_\infty-\gamma$}\\
             \pm  \rho^{\gamma_\infty}_{r+g}\left( (\bar{\mu} \wedge \bar{\nu})\otimes \left(\bar u \bar x\bar z\frac{d \bar{z}}{\bar{z}}\wedge \frac{d\bar{\mathbf{y}}}{\bar{\mathbf{y}}} \wedge\frac{d\bar{\mathbf{t}}}{\bar{\mathbf{t}}}\right)\right) &\text{if $j=j'+1$ and $\gamma < \gamma_\infty-\gamma$}\\
             0 & \text{otherwise}
         \end{cases}\\
         \intertext{\raggedleft (by Lemma \ref{lem trivialite bcp de slots} and with $\bar u:=(\bar u'_{\gamma_\infty})^{-1}\bar u_{\gamma,\gamma_\infty-\gamma}$)  }
         &=
         \begin{cases}
             \pm \rho^{\gamma_\infty}_{r+g}\left( (\bar{\mu} \wedge \bar{\nu} \wedge\bar\gamma)\otimes \left(\bar u\bar x\bar z \frac{d\bar{\mathbf{y}}}{\bar{\mathbf{y}}} \wedge\frac{d\bar{\mathbf{t}}}{\bar{\mathbf{t}}}\right)\right)
             &\text{if $j=j'$}\\
              \pm \rho^{\gamma_\infty}_{r+g}\left( (\bar{\mu} \wedge \bar{\nu})\otimes \left(\bar u \bar x\bar z\frac{d \bar{x}}{\bar{x}}\wedge \frac{d\bar{\mathbf{y}}}{\bar{\mathbf{y}}} \wedge\frac{d\bar{\mathbf{t}}}{\bar{\mathbf{t}}}\right)\right) &\text{if $j=j'+1$}\\
             0 & \text{otherwise}
         \end{cases}\\
         \intertext{\raggedleft (because $(\bar x d\bar z + \bar z d\bar x) \frac{d\bar{\mathbf{y}}}{\bar{\mathbf{y}}} \wedge\frac{d\bar{\mathbf{t}}}{\bar{\mathbf{t}}} \in \ker (d)$)  }
         &=\pm \rho^{\gamma_\infty}_{r+g}\left(\varphi\left( \bar{\mu}\otimes\bar{x} \frac{d\bar{\mathbf{y}}}{\bar{\mathbf{y}}},
        \bar{\nu}\otimes\bar{z} \frac{d\bar{\mathbf{t}}}{\bar{\mathbf{t}}} \right)\right).
    \end{align*}
    \item Write $\gamma=p\gamma'$. Choose $j,j'\geq 0$. Take: $$\begin{cases}
        (\bar{\mu}_1 \wedge \dots \wedge \bar{\mu}_j)\otimes\left\{\bar{x} \frac{d \bar{y}_{j+2}}{\bar{y}_{j+2}}\wedge \dots \wedge\frac{d \bar{y}_{q}}{\bar{y}_{q}}\right\}\in \Lambda^{j}(\Gamma/p\Gamma) \otimes \Omega_k^{q-j-1}\\
        (\bar{\nu}_1 \wedge \dots \wedge \bar{\nu}_{g-j'})\otimes\left\{\bar{z} \frac{d \bar{t}_{g-j'+2}}{\bar{t}_{g-j'+2}}\wedge \dots \wedge\frac{d \bar{t}_{r+g-q}}{\bar{t}_{r+g-q}}\right\}\in\Lambda^{g-j'}(\Gamma/p\Gamma) \otimes \Omega_k^{r+j'-q-1},
    \end{cases}$$
    and, as before, set for simplicity: 
    \begin{gather*}
    \bar{\mu}:=\bar{\mu}_1 \wedge \dots \wedge \bar{\mu}_j,\;\;\;\;\; \bar{\nu}:=\bar{\nu}_1 \wedge \dots \wedge \bar{\nu}_{g-j'},\\ 
    \pi^{\bar{\mu}}:=(\pi^{\bar\mu_1},\dots , \pi^{\bar\mu_j}),\;\;\;\;\; \pi^{\bar{\nu}}:=(\pi^{\bar\nu_1},\dots,\pi^{\bar\nu_{g-j'}}),\\
        \mathbf{y}:=(y_{j+2},\dots,y_q),\;\;\;\;\; \mathbf{t}:=(t_{g-j'+2},\dots,t_{r+g-q}),\\ \frac{d\bar{\mathbf{y}}}{\bar{\mathbf{y}}}:=\frac{d \bar{y}_{j+2}}{\bar{y}_{j+2}}\wedge \dots \wedge\frac{d \bar{y}_{q}}{\bar{y}_{q}},\;\;\;\;\; \frac{d\bar{\mathbf{t}}}{\bar{\mathbf{t}}}:=\frac{d \bar{t}_{g-j'+2}}{\bar{t}_{g-j'+2}}\wedge \dots \wedge\frac{d \bar{t}_{r+g-q}}{\bar{t}_{r+g-q}}.
    \end{gather*}
    We compute:
    \begin{align*}
         \rho_q^\gamma &\left(\bar{\mu}\otimes\left\{\bar{x}  \frac{d\bar{\mathbf{y}}}{\bar{\mathbf{y}}}\right\}\right) \cup \rho_{r+g-q}^{\gamma_\infty-\gamma}\left(\bar{\nu}\otimes\left\{\bar{z}  \frac{d\bar{\mathbf{t}}}{\bar{\mathbf{t}}}\right\}\right)\\
         & =\pm \left\{ \pi^{\bar\mu},\pi^{\bar\nu},1+(\pi^{\gamma'})^p x,1+(\pi^{\gamma_\infty/p-\gamma'})^pz,\mathbf{y},\mathbf{t}\right\}\\
         &=\begin{cases}
             \pm \left\{ \pi^{\bar\mu},\pi^{\bar\nu},1-(\pi^{\gamma'})^p (\pi^{\gamma_\infty/p-\gamma'})^p xz,-(\pi^{\gamma'})^p x,\mathbf{y},\mathbf{t}\right\}
             &\text{if $\gamma \geq \gamma_\infty-\gamma$}\\
              \pm \left\{ \pi^{\bar\mu},\pi^{\bar\nu},1-(\pi^{\gamma'})^p (\pi^{\gamma_\infty/p-\gamma'})^p xz,-(\pi^{\gamma_\infty/p-\gamma'})^p z,\mathbf{y},\mathbf{t}\right\}
               &\text{if $\gamma < \gamma_\infty-\gamma$}
         \end{cases}\\
         \intertext{\raggedleft (by Lemma \ref{lem suma} and the inequality $\max \{\frac{e+\gamma}{2},e-\gamma\}<\gamma_{\infty}<e$)  }
         &=
         \begin{cases}
              \pm \rho^{\gamma_\infty}_{r+g}\left( (\bar{\mu} \wedge \bar{\nu})\otimes \left(\bar u \bar x\bar z\frac{d \bar{x}}{\bar{x}}\wedge \frac{d\bar{\mathbf{y}}}{\bar{\mathbf{y}}} \wedge\frac{d\bar{\mathbf{t}}}{\bar{\mathbf{t}}}\right)\right) &\text{if $j=j'$ and $\gamma \geq \gamma_\infty-\gamma$}\\
              \pm  \rho^{\gamma_\infty}_{r+g}\left( (\bar{\mu} \wedge \bar{\nu})\otimes \left(\bar u\bar x\bar z\frac{d \bar{z}}{\bar{z}}\wedge \frac{d\bar{\mathbf{y}}}{\bar{\mathbf{y}}} \wedge\frac{d\bar{\mathbf{t}}}{\bar{\mathbf{t}}}\right)\right)&\text{if $j=j'$ and $\gamma < \gamma_\infty-\gamma$}\\
             0& \text{otherwise}
         \end{cases}\\
         \intertext{\raggedleft (by Lemma \ref{lem trivialite bcp de slots} and with $\bar u := \bar u_{\gamma,\gamma_\infty-\gamma} \bar u'_{\gamma}\bar u'_{\gamma_\infty-\gamma} (\bar u'_{\gamma_{\infty}})^{-1}$)  }
         &=\begin{cases}
              \pm \rho^{\gamma_\infty}_{r+g}\left( (\bar{\mu} \wedge \bar{\nu})\otimes \left(\bar u\bar x\bar z\frac{d \bar{x}}{\bar{x}}\wedge \frac{d\bar{\mathbf{y}}}{\bar{\mathbf{y}}} \wedge\frac{d\bar{\mathbf{t}}}{\bar{\mathbf{t}}}\right)\right) & \text{if $j=j'$}\\
             0 & \text{otherwise}
         \end{cases}\\
         \intertext{\raggedleft (because $(\bar x d\bar z + \bar z d\bar x) \frac{d\bar{\mathbf{y}}}{\bar{\mathbf{y}}} \wedge\frac{d\bar{\mathbf{t}}}{\bar{\mathbf{t}}} \in \ker (d)$)  }
         &= \rho^{\gamma_\infty}_{r+g}\left(\varphi\left( \bar{\mu} \otimes\bar{x} \frac{d\bar{\mathbf{y}}}{\bar{\mathbf{y}}},
        \bar{\nu}\otimes\bar{z} \frac{d\bar{\mathbf{t}}}{\bar{\mathbf{t}}} \right)\right).
    \end{align*}
\end{enumerate}
\end{proof}

\begin{corollary}
    The morphism $\rho_q^\gamma$ is injective.
\end{corollary}

\begin{proof}
    We distinguish two cases depending on whether $p \mid \gamma$ or $p\nmid \gamma$. In both cases, we will work with the rational function field $K(t)$ endowed with the Gauss valuation. Its residue field is $k(\bar t)$ (cf.~\cite[Cor.~2.2.2]{EnglerPrestel}).
    \begin{enumerate}
        \item Assume first that $p\nmid \gamma$. Take some non-zero $\alpha=(\alpha_j)_j \in \bigoplus_{j=0}^{q-1} (\Lambda^{j}\Delta \otimes \Omega_k^{q-j-1})$. Let $j_0$ be the largest integer such that $\alpha_{j_0}\neq 0$. We have a perfect pairing of finite dimensional $k$-vector spaces:
    \begin{align*}
        \langle \cdot, \cdot\rangle_k: (\Lambda^{j_0}\Delta \otimes \Omega_k^{q-j_0-1}) \times (\Lambda^{g-j_0-1}\Delta \otimes \Omega_k^{r+j_0+1-q})&\rightarrow \Lambda^{g-1}\Delta \otimes \Omega_k^{r}\\
        (\delta \otimes \omega,\delta'\otimes \omega') &\mapsto (\delta\wedge\delta') \otimes (\omega\wedge \omega'),
    \end{align*}
    We can therefore find $\beta \in \Lambda^{g-j_0-1}\Delta \otimes \Omega_k^{r+j_0+1-q}$ such that $\langle \alpha_{j_0}, \beta\rangle_k$ is of the form $(\bar\gamma_1\wedge \dots \wedge \bar\gamma_{g-1} ) \otimes (d\bar y_1 \wedge \dots d\bar y_r)$ for some $\mathbb{F}_p$-basis $(\bar \gamma_1,\dots,\bar \gamma_{g-1})$ of $\Delta$ and some $p$-basis $(\bar y_1,\dots,\bar y_r)$ of $k$. 

    Now consider the pairing:
    \begin{align*}
        \langle \cdot, \cdot\rangle_{k(\bar t)}: (\Lambda^{j_0}\Delta \otimes \Omega_{k(\bar t)}^{q-j_0-1}) \times (\Lambda^{g-j_0-1}\Delta \otimes \Omega_{k(\bar t)}^{r+j_0+1-q})&\rightarrow \Lambda^{g-1}\Delta \otimes \Omega_{k(\bar t)}^{r}\\
        (\delta \otimes \omega,\delta'\otimes \omega') &\mapsto (\delta\wedge\delta') \otimes (\omega\wedge \omega'),
    \end{align*}
    Setting $\bar u := (\bar u'_{\gamma_\infty})^{-1}\bar u_{\gamma,\gamma_\infty-\gamma}$, we have:
    $$\bar \gamma\wedge \langle \alpha_{j_0}, \bar u^{-1} \bar t \beta\rangle_{k(\bar t)} = (\bar \gamma \wedge \bar\gamma_1\wedge \dots \wedge \bar\gamma_{g-1} ) \otimes (\bar u^{-1} \bar t d\bar y_1 \wedge \dots d\bar y_r). $$
    Hence, by Corollary \ref{cor non trivial}:
    $$\rho_{r+g+1}^{\gamma_{\infty}}(\bar u (\bar \gamma \wedge \langle \alpha_{j_0},\bar u^{-1} \bar t\beta\rangle_{k(\bar t)}))\neq 0.$$ 
    But if we see $\beta$ as an element of $\bigoplus_{j} (\Lambda^{g-j-1}\Delta \otimes \Omega_{k(\bar t)}^{r+j+1-q})$, Lemma \ref{lem calcul} and the choice of $j_0$ show that
    \[\rho_q^{\gamma}(\alpha)\cup \rho_{r+g+1-q}^{\gamma_\infty-\gamma}(\bar u^{-1} \bar t\beta)=\pm\rho_{r+g+1}^{\gamma_\infty}(\bar u (\bar \gamma \wedge \langle \alpha_{j_0},\bar u^{-1} \bar t\beta\rangle_{k(\bar t)}))\neq 0 \in u^{\gamma_{\infty}}_{r+g+1}(K(t))/u^{>\gamma_{\infty}}_{r+g+1}(K(t)).\]
    We deduce that $\rho_q^{\gamma}(\alpha)$ is non-zero in $u^{\gamma}_{q}(K(t))/u^{>\gamma}_{q}(K(t))$, and hence in $u^{\gamma}_{q}(K)/u^{>\gamma}_{q}(K)$.
    \item Assume now that $p\mid \gamma$. Take some non-zero $\alpha=(\alpha_j)_j \in \bigoplus_{j=0}^{q-1} \left(\Lambda^{j}(\Gamma/p\Gamma) \otimes \frac{\Omega_{k(\bar t)}^{q-j-1}}{\ker(d)}\right)$. Let $j_0$ be any integer with $\alpha_{j_0}\neq 0$, so that: $$(\mathrm{id}_{\Lambda^{j_0}(\Gamma/p\Gamma)}\otimes d)(\alpha_{j_0})\neq 0 \in \Lambda^{j_0}(\Gamma/p\Gamma) \otimes \Omega_k^{q-j_0}.$$ We have a perfect pairing of finite dimensional $k$-vector spaces: 
    \begin{align*}
       \langle \cdot , \cdot \rangle_k: \left(\Lambda^{j_0}(\Gamma/p\Gamma) \otimes\Omega^{q-j_0}_k\right) \times \left( \Lambda^{g-j_0}(\Gamma/p\Gamma) \otimes\Omega^{r+j_0-q}_k\right) &\rightarrow \Lambda^{g}(\Gamma/p\Gamma) \otimes\Omega_k^{r}\\
       (\delta\otimes  \omega,\delta'\otimes \omega') &\mapsto (\delta\wedge \delta')\otimes (\omega\wedge \omega').
    \end{align*}
    We can therefore find $\beta \in \Lambda^{g-j_0}(\Gamma/p\Gamma) \otimes\Omega^{r+j_0-q}_k$ such that $\langle (\mathrm{id}_{\Lambda^{j_0}(\Gamma/p\Gamma)}\otimes d)(\alpha_{j_0}), \beta\rangle_{k}$ is of the form $(\bar\gamma_1\wedge \dots \wedge \bar\gamma_{g} ) \otimes (d\bar y_1 \wedge \dots d\bar y_r)$ for some $\mathbb{F}_p$-basis $(\bar \gamma_1,\dots,\bar \gamma_{g})$ of $\Gamma/p\Gamma$ and some $p$-basis $(\bar y_1,\dots,\bar y_r)$ of $k$. 

Now consider the pairing:
    \begin{align*}
       \langle \cdot , \cdot \rangle_{k(\bar t)}: \left(\Lambda^{j_0}(\Gamma/p\Gamma) \otimes\Omega^{q-j_0}_{k(\bar t)}\right) \times \left( \Lambda^{g-j_0}(\Gamma/p\Gamma) \otimes\Omega^{r+j_0-q}_{k(\bar t)}\right) &\rightarrow \Lambda^{g}(\Gamma/p\Gamma) \otimes\Omega_{k(\bar t)}^{r}\\
       (\delta\otimes  \omega,\delta'\otimes \omega') &\mapsto (\delta\wedge \delta')\otimes (\omega\wedge \omega').
    \end{align*}
    Setting $\bar u := \bar u_{\gamma,\gamma_\infty-\gamma} \bar u'_{\gamma}\bar u'_{\gamma_\infty-\gamma} (\bar u'_{\gamma_{\infty}})^{-1}$, we have:
    $$ \langle (\mathrm{id}_{\Lambda^{j_0}(\Gamma/p\Gamma)}\otimes d)(\alpha_{j_0}), \bar u^{-1} \bar t \beta\rangle_{k(\bar t)} = ( \bar\gamma_1\wedge \dots \wedge \bar\gamma_{g} ) \otimes (\bar u^{-1} \bar t d\bar y_1 \wedge \dots d\bar y_r). $$
    Hence, by Corollary \ref{cor non trivial}:
    $$\rho_{r+g+1}^{\gamma_{\infty}}(\bar u \langle (\mathrm{id}_{\Lambda^{j_0}(\Gamma/p\Gamma)}\otimes d)(\alpha_{j_0}),\bar u^{-1} \bar t\beta\rangle_{k(\bar t)})\neq 0.$$ 
    But if we see $\beta$ as an element of $\bigoplus_{j} (\Lambda^{g-j}(\Gamma/p\Gamma) \otimes \Omega_{k(\bar t)}^{r+j-q})$, Lemma \ref{lem calcul} shows that $\rho_q^{\gamma}(\alpha)\cup \rho_{r+g+1-q}^{\gamma_\infty-\gamma}(\bar u^{-1} \bar t\beta)=\pm\rho_{r+g+1}^{\gamma_\infty}(\bar u \langle (\mathrm{id}_{\Lambda^{j_0}(\Gamma/p\Gamma)}\otimes d)(\alpha_{j_0}),\bar u^{-1} \bar t\beta\rangle_{k(\bar t)})\neq 0 \in u^{\gamma_{\infty}}_{r+g+1}(K(t))/u^{>\gamma_{\infty}}_{r+g+1}(K(t))$. We deduce that $\rho_q^{\gamma}(\alpha)$ is non-zero in $u^{\gamma}_{q}(K(t))/u^{>\gamma}_{q}(K(t))$, and hence in $u^{\gamma}_{q}(K)/u^{>\gamma}_{q}(K)$.
    \end{enumerate}
\end{proof}

\subsection{Step 5: the case $\gamma=e$}

Let $K^\unr$ be the unramified extension whose residue field is $k^s$, the separable closure of $k$, and whose value group is $\Gamma$, the same as $K$. The Galois group $\Gal(K^\unr/K)$ is isomorphic to $\Gal(k)=\Gal(k^s/k)$ and has then $p$-cohomological dimension $\leq 1$. In particular, if we consider the Hochschild-Serre spectral sequence
\[H^i(k,H^j(K^\unr,\Z/p\Z))\Rightarrow H^{i+j}(K,\Z/p\Z),\]
then it yields a short exact sequence
\[0\to H^1(k,H^j(K^\unr,\Z/p\Z))\to H^{j+1}(K,\Z/p\Z) \to H^{j+1}(K^\unr,\Z/p\Z)^{\Gal(k)}\to 0,\]
for any $j\geq 0$. Using the Bloch-Kato conjecture, we may translate the first arrow into
\[\iota_q: H^1(k,\kk_q(K^\unr))\hookrightarrow \kk_{q+1}(K),\]
for any $q\geq 0$.

We will prove the injectivity of $\rho_q^e$ by a diagram chase that involves this morphism $\iota_q$ and a second morphism $\varphi_q$ defined as the following composition of morphisms
\begin{eqnarray*}
H^1(k,\kk_q(K^\unr)) &\twoheadrightarrow & H^1(k,\kk_q(K^\unr)/u_q^{>0}(K^\unr)) \\
& \xrightarrow[\sim]{(\rho_q^0)^{-1}_*} & H^1\left(k,\bigoplus_{j=0}^{q} \left(\Lambda^j(\Gamma/p\Gamma) \otimes \kk_{q-j}(k^s)\right)\right)\\
& \xrightarrow{\sim} & \bigoplus_{j=0}^{q} \left(\Lambda^j(\Gamma/p\Gamma)\otimes H^1(k,\kk_{q-j}(k^s))\right)\\
& \xrightarrow[\sim]{((-1)^j\otimes(\psi_{k^s}^{q-j})_*)_j} & \bigoplus_{j=0}^{q} \left(\Lambda^j(\Gamma/p\Gamma)\otimes H^1(k,\nu(q-j)(k^s))\right)\\
& \xrightarrow[\sim]{(\mathrm{id}\otimes \delta_{q-j}^{-1})_j} & \bigoplus_{j=0}^{q} \left(\Lambda^j(\Gamma/p\Gamma)\otimes H_p^{q-j+1}(k)\right)
\end{eqnarray*}
where $\psi_{k^s}^{q-j}:\kk_{q-j}(k^s)\to \nu(q-j)(k^s)$ is the differential symbol (cf.~\cite[\S9.5]{GS}) and \[\delta_{q-j}:H_p^{q-j+1}(k)=\left(\Omega_{k}^{q-j}/\im(d)\right)/\wp(\Omega_k^{q-j})\to H^1(k,\nu(q-j)(k^s)),\]
is the isomorphism induced by the connecting map associated to the short exact sequence of $\Gal(k)$-modules
\[0\to\nu(q-j)(k^s)\to \Omega_{k^s}^{q-j} \xrightarrow{\wp} \Omega_{k^s}^{q-j}/\im(d) \to 0.\]
Let us make some remarks on these maps defining $\varphi_q$:
\begin{itemize}
    \item The first map of the composition is indeed surjective since $\Gal(k)$ has $p$-cohomological dimension $\leq 1$ and hence $H^2(k,u_q^{>0}(K^\unr))=0$.
    \item The second map is well-defined since from the construction of the morphisms $\rho_q^\gamma$ one can immediately see that they are functorial with respect to morphisms of valued fields with same value group. In particular, they are naturally $\Gal(k^s/k)$-equivariant.
    \item The third map is naturally an isomorphism since $\Gamma/p\Gamma$ is an $\F_p$-vector space with trivial action of $\Gal(k)$.
    \item The fourth map is also well-defined since the differential symbol is also functorial and hence $\Gal(k^s/k)$-equivariant.
    \item The second, fourth and fifth maps are isomorphisms since they are induced by isomorphisms.
\end{itemize}
In particular, $\varphi_q$ is surjective.\\

We establish now the injectivity of $\rho_1^e$, which is the stepping stone for the injectivity of $\rho_q^e$ in general.

\begin{lemma}\label{lem comm rho^e}
We have $\iota_0=\rho_1^e\circ \varphi_0$. In other words, the following diagram is commutative
\[\xymatrix{
H^1(k,\Z/p\Z) \ar@{^{(}->}[r]^>>>>>{\iota_0} \ar@{->>}[d]_{\varphi_0} & \kk_{1}(K) \\
H_p^{1}(k) \ar[ur]_{\rho_{1}^e}.
}\]
In particular, since $\iota_0$ is injective, so is $\rho_1^e$.
\end{lemma}

\begin{proof}
Recall that $\iota_0$ is obtained via the Bloch-Kato conjecture from the map
\[H^1(k,H^0(K^\unr,\Z/p\Z))\to H^1(K,\Z/p\Z),\]
obtained from a spectral sequence. But since $H^1(k,\Z/p\Z)=H^1(k,H^0(K^\unr,\Z/p\Z))$, this is just an inflation map, which sends a map $\phi:\Gal(k)\to\Z/p\Z$ representing a class in $H^1(k,\Z/p\Z)$ to the class of the composition $\phi\circ\pi:\Gal(K)\to\Z/p\Z$, where $\pi:\Gal(K)\to\Gal(k)$ is the natural projection.

On the other hand, the map $\varphi_0$ is composed of five maps. But the first four are all the identity on $H^1(k,\Z/p\Z)$. The last one however corresponds to the isomorphism $H^1(k,\Z/p\Z)\to k/\wp(k)$ induced by the Artin-Schreier sequence
\[0\to \Z/p\Z \to k \xrightarrow{\wp} k \to 0.\]
In particular, $\varphi_0$ is invertible and $\varphi_0^{-1}$ sends the class of $\bar x$ in $k/\wp(k)=H^1_p(k)$ to the class of the map $\phi:\Gal(k)\to\Z/p\Z$, where the kernel of $\phi$ corresponds to the extension given by the polynomial $P(\bar T)=\bar T^p-\bar T-\bar x$ and, if $\bar \alpha$ is a root of $P$, then $\sigma(\bar \alpha)=\bar \alpha+\phi(\sigma)$ for $\sigma\in\Gal(k)$.

Finally, $\rho_{1}^e$ sends the class of $\bar x$ in $k/\wp(k)$ to the symbol $\{1+(\zeta_p-1)^p x\}$, where $x\in K$ is any lift of $\bar x\in k$. But this symbol represents via Kummer Theory the cyclic extension $K(\sqrt[p]{1+(\zeta_p-1)^p x})/K$ and, if $\beta$ is a root of $T^p-(1+(\zeta_p-1)^p x)$, then the corresponding morphism $\psi:\Gal(K)\to\Z/p\Z$ is such that $\sigma(\beta)=\zeta_p^{\psi(\sigma)}\beta$.

In order to conclude, note that $\beta-1$ is a root of
\[T^p+\sum_{i=1}^{p-1}\binom{p}{i} T^i-(\zeta_p-1)^p x,\]
and hence $\alpha:=(\beta-1)/(\zeta_p-1)$ is a root of
\[T^p+\sum_{i=1}^{p-1}\binom{p}{i}(\zeta_p-1)^{i-p}T^i- x.\]
Since $v(p)=(p-1)v(\zeta_p-1)$ and, more precisely,
\[\frac{p}{(\zeta_p-1)^{p-1}}\equiv -1\pmod{p},\]
we see that $\bar\alpha$ is a root of $P(\bar T)=\bar T^p-\bar T-\bar x$ and hence $K(\sqrt[p]{1+(\zeta_p-1)^p x})/K$ is unramified and induces the expected extension over $k$. Moreover, we check that, for $\sigma\in\Gal(K)$
\begin{multline*}
(\pi(\sigma))(\bar \alpha)=\overline{\sigma(\alpha)}=\overline{\left(\frac{\zeta_p^{\psi(\sigma)}\beta-1}{\zeta_p-1}\right)}=\overline{\left(\frac{\zeta_p^{\psi(\sigma)}\beta-\zeta_p^{\psi(\sigma)}}{\zeta_p-1}\right)}+\overline{\left(\frac{\zeta_p^{\psi(\sigma)}-1}{\zeta_p-1}\right)}\\
=\bar\zeta_p^{\psi(\sigma)}\overline{\left(\frac{\beta-1}{\zeta_p-1}\right)}+\left(1+\bar \zeta_p+\cdots +\bar \zeta_p^{\psi(\sigma)-1}\right)=\bar\alpha+\psi(\sigma).
\end{multline*}
This tells us that the image $\pi(\sigma)\in\Gal(k)$ of $\sigma\in\Gal(K)$ acts precisely as expected, which proves the commutativity of the diagram.
\end{proof}

In order to get the injectivity of $\rho_q^e$ in general, we need some other commutative diagrams.

\begin{lemma}\label{prop comm rho^e 1}
There is a commutative diagram
\[{\small\xymatrix@C=0pt{
H^1(k,\Z/p\Z)\otimes\kk_q(K) \ar[r]^-{\mathrm{id}\otimes\mathrm{Res}} \ar@{->>}[d]_-{\varphi_0\otimes((\rho_q^{0})^{-1}\circ\pi)} & H^1(k,\Z/p\Z)\otimes H^0(k,\kk_q(K^\unr)) \ar[r]^-{\cup} & H^1(k,\kk_{q}(K^\unr)) \ar@{->>}[d]^{\varphi_{q+1}} \\
H_p^1(k)\otimes \left( \displaystyle\bigoplus_{j=0}^{q} \left(\Lambda^j(\Gamma/p\Gamma)\otimes \kk_{q-j}(k)\right)\right) \ar@{->>}[rr]^{((-1)^j\cup\psi_k^{q-j})_j} & & \displaystyle\bigoplus_{j=0}^{q} \left(\Lambda^j(\Gamma/p\Gamma)\otimes H_p^{q-j+1}(k)\right),
}}\]
where $\pi:\kk_q(K)\to \kk_q(K)/u_q^{>0}(K)$ is the canonical projection and the lower horizontal arrow is induced by the differential symbol and concatenation, that is
\[\xymatrix@R=0pt{
H_p^1(k)\otimes \kk_{q-j}(k) \ar[rr]^-{(-1)^j\otimes \psi_k^{q-j}}_\sim && H_p^1(k)\otimes \nu(q-j)(k) \ar@{->>}[r]^-{\cup} & H_p^{q-j+1}(k),\\
&& \bar x\otimes \left(\frac{d\bar y_1}{\bar y_1}\wedge\cdots\wedge\frac{d\bar y_{q-j}}{\bar y_{q-j}}\right) \ar@{|->}[r] & \bar x\frac{d \bar y_1}{\bar y_1}\wedge\cdots\wedge\frac{d\bar y_{q-j}}{\bar y_{q-j}}.
}
\]
\end{lemma}

\begin{proof}
Since $\varphi_0$ and $\varphi_{q+1}$ are both defined as a composition of five maps, we can check for each one of these that we have a corresponding commutative diagram.

The whole diagram is
\[\xymatrix{
H^1(k,\Z/p\Z)\otimes\kk_q(K) \ar[r]^{\mathrm{id}\cup\mathrm{Res}} \ar@{->>}[d] & H^1(k,\kk_{q}(K^\unr)) \ar@{->>}[d] \\
H^1(k,\Z/p\Z)\otimes\left(\kk_q(K)/u_q^{>0}(K)\right) \ar[r]^{\mathrm{id}\cup\mathrm{Res}} \ar[d]_-{\mathrm{id}\otimes (\rho_q^0)^{-1}} & H^1(k,\kk_{q}(K^\unr)/u_q^{>0}(K^\unr)) \ar[d]^-{(\rho_q^0)^{-1}_*} \\
H^1(k,\Z/p\Z)\otimes \left(\displaystyle\bigoplus_{j=0}^{q} \left(\Lambda^j(\Gamma/p\Gamma) \otimes \kk_{q-j}(k)\right)\right) \ar[r]^-{\mathrm{id}\cup\mathrm{Res}} \ar[d] & H^1\left(k,\displaystyle\bigoplus_{j=0}^{q} \left(\Lambda^j(\Gamma/p\Gamma) \otimes \kk_{q-j}(k^s)\right)\right) \ar[d] \\
\displaystyle\bigoplus_{j=0}^{q} \left(\Lambda^j(\Gamma/p\Gamma) \otimes H^1(k,\Z/p\Z)\otimes \kk_{q-j}(k)\right) \ar[r]^-{\mathrm{id}\cup\mathrm{Res}} \ar[d]_{((-1)^j\otimes(\psi_k^{q-j})_*)j} & \displaystyle\bigoplus_{j=0}^{q} \left(\Lambda^j(\Gamma/p\Gamma) \otimes H^1(k,\kk_{q-j}(k^s))\right) \ar[d]^{((-1)^j\otimes(\psi_{k^s}^{q-j})_*)_j}  \\
\displaystyle\bigoplus_{j=0}^{q} \left(\Lambda^j(\Gamma/p\Gamma) \otimes H^1(k,\Z/p\Z)\otimes \nu(q-j)(k)\right) \ar[r]^-{\mathrm{id}\cup\mathrm{Res}} \ar[d]_{(\mathrm{id}\otimes \delta_0^{-1}\otimes\mathrm{id})_j} & \displaystyle\bigoplus_{j=0}^{q} \left(\Lambda^j(\Gamma/p\Gamma) \otimes H^1(k,\nu(q-j)(k^s))\right) \ar[d]^{(\mathrm{id}\otimes \delta_{q-j}^{-1})_j} \\
\displaystyle\bigoplus_{j=0}^{q} \left(\Lambda^j(\Gamma/p\Gamma) \otimes H^1_p(k)\otimes \nu(q-j)(k)\right) \ar[r]^-{\mathrm{id}\otimes \cup} & \displaystyle\bigoplus_{j=0}^{q} \left(\Lambda^j(\Gamma/p\Gamma)\otimes H_p^{q-j+1}(k)\right).
}\]
The first square clearly commutes since both restriction and cup products are functorial. The functoriality of $\rho_q^0$ with respect to unramified extensions makes it compatible with restriction, so that the second square commutes. The third one is evident. The fourth square is commutative by functoriality of the differential symbol. The fifth diagram is also commutative thanks to the compatibility of connecting maps with cup products.

In order to conclude, note that the upper-right half of the rectangle on the statement of the proposition is exactly the upper-right half of the big rectangle here above. On the other hand, the lower-left parts of both diagrams also coincide by a direct check. This concludes the proof of commutativity. The fact that the vertical arrows are surjective is because so are $\varphi_0$, $\varphi_{q+1}$, $(\rho_q^0)^{-1}$ and $\pi$. Surjectivity of the lower horizontal arrow follows directly from the formula and the bijectivity of the differential symbol.
\end{proof}

\begin{lemma}\label{prop comm rho^e 2}
There is a commutative diagram
\[\xymatrix{
\kk_1(K)\otimes\kk_q(K) \ar@{->>}[rr]^{\cup} && \kk_{q+1}(K)
\\
H^1(k,\Z/p\Z)\otimes\kk_q(K) \ar[rr]^{\mathrm{id}\cup\mathrm{Res}} \ar@{^{(}->}[u]_{\iota_0\otimes\mathrm{id}} && H^1(k,\kk_{q}(K^\unr)). \ar@{^{(}->}[u]^{\iota_{q}}
}\]
\end{lemma}

\begin{proof}
Since the arrows $\iota_0$ and $\iota_q$ are defined by a spectral sequence via the norm residue isomorphism, this is just the compatibility between spectral sequences and cup-products.
\end{proof}

\begin{lemma}\label{prop comm rho^e 3}
There is a commutative diagram
\[\xymatrix{
\kk_1(K)\otimes\kk_q(K) \ar@{->>}[r]^{\cup} &  \kk_{q+1}(K).\\
H^1(k,\Z/p\Z)\otimes\kk_q(K) \ar@{->>}[d]_-{\varphi_0\otimes((\rho_q^0)^{-1}\circ\pi)} \ar@{^{(}->}[u]^{\iota_0\otimes\mathrm{id}} & \\
H_p^1(k)\otimes \left( \displaystyle\bigoplus_{j=0}^{q} \left(\Lambda^j(\Gamma/p\Gamma)\otimes \kk_{q-j}(k)\right)\right) \ar@/_3pc/[uru]_-{\rho_1^e\cup\rho_q^0} &
}\]
\end{lemma}

\begin{proof}
Consider the commutative triangle from Lemma \ref{lem comm rho^e} and take its tensor product with $\kk_q(K)$. After composing with cup product in Milnor $\mathrm{K}$-theory, this yields
\[\xymatrix{
\kk_1(K)\otimes\kk_q(K) \ar@{=}[r] & \kk_1(K)\otimes\kk_q(K) \ar@{->>}[r]^-{\cup} &  \kk_{q+1}(K).\\
H^1(k,\Z/p\Z)\otimes\kk_q(K) \ar@{->>}[d]_-{\varphi_0\otimes\mathrm{id}} \ar@{^{(}->}[u]^{\iota_0\otimes\mathrm{id}} & \\
H_p^1(k)\otimes \kk_q(K) \ar@/_1pc/[uru]_-{\rho_1^e\otimes\mathrm{id}} \ar@/_3pc/[urru]_-{\rho_1^e\cup\mathrm{id}}
}\]
Now, since the image of $\rho_1^e$ lies in $u_1^e(K)$, Lemma \ref{lem suma} tells us that the arrow $\rho_1^e\cup\mathrm{id}$ is trivial on $H^1_p(k)\otimes u_q^{>0}(K)$. Then the diagram becomes
\[\xymatrix{
\kk_1(K)\otimes\kk_q(K) \ar@{->>}[r]^-{\cup} &  \kk_{q+1}(K).\\
H^1(k,\Z/p\Z)\otimes\kk_q(K) \ar@{->>}[d]_-{\varphi_0\otimes\pi} \ar@{^{(}->}[u]^{\iota_0\otimes\mathrm{id}} & \\
H_p^1(k)\otimes \left(\kk_q(K)/u_q^{>0}(K)\right) \ar@/_3pc/[uru]_-{\rho_1^e\cup\mathrm{id}}
}\]
Since $\rho_q^0$ is an isomorphism between $\kk_q(K)/u_q^{>0}(K)$ and $\bigoplus_{j=0}^{q} (\Lambda^j(\Gamma/p\Gamma)\otimes \kk_{q-j}(k))$, the result follows.
\end{proof}

\begin{lemma}\label{prop comm rho^e 4}
There is a commutative diagram
\[\xymatrix{
&& \kk_{q+1}(K)\\
H_p^1(k)\otimes \left( \displaystyle\bigoplus_{j=0}^{q} \left(\Lambda^j(\Gamma/p\Gamma)\otimes \kk_{q-j}(k)\right)\right) \ar[rru]^-{\rho_1^e\cup\rho_q^0} \ar[rr]_-{((-1)^j\cup\psi_k^{q-j})_j}  && \displaystyle\bigoplus_{j=0}^{q} \left(\Lambda^j(\Gamma/p\Gamma)\otimes H_p^{q-j+1}(k)\right). \ar[u]_{\rho_{q+1}^e}
}\]

\end{lemma}

\begin{proof}
   Take: $$\alpha=\bar x \otimes (\bar\gamma_1\wedge\dots\wedge\bar \gamma_j)\otimes \{\bar y_1,\dots,\bar y_{q-j}\} \in H_p^1(k)\otimes \Lambda^j(\Gamma/p\Gamma)\otimes \kk_{q-j}(k) .$$
    We then have:
    \begin{align*}
    \rho^e_{q+1}((-1)^j\cup\psi_k^{q-j}(\alpha))&= (-1)^j\rho^e_{q+1}\left((\bar\gamma_1\wedge\dots\wedge\bar \gamma_j)\otimes \left(\bar x \frac{d\bar y_1}{\bar y_1}\wedge \dots \wedge \frac{d\bar y_{q-j}}{\bar y_{q-j}}\right)\right)\\
    &= (-1)^j\{\pi^{\bar\gamma_1},\dots,\pi^{\bar\gamma_j},1+(\zeta_p-1)^px,y_1,\dots,y_{q-j}\}\\
    &=\{1+(\zeta_p-1)^px,\pi^{\bar\gamma_1},\dots,\pi^{\bar\gamma_j},y_1,\dots,y_{q-j}\}\\
    &=(\rho_1^e\cup \rho_q^0)(\alpha).
    \end{align*}
\end{proof}

\begin{proposition}
Set $I:=\mathrm{Im}\left(\mathrm{id} \cup \mathrm{Res}: H^1(k,\Z/p\Z)\otimes\kk_q(K) \rightarrow H^1(k,\kk_{q}(K^\unr)) \right)$. Then the diagram:
\[\xymatrix{
\kk_1(K)\otimes\kk_q(K) \ar@{->>}[rr]^{\cup} && \kk_{q+1}(K).\\
H^1(k,\Z/p\Z)\otimes\kk_q(K) \ar@{->>}[r]^-{\mathrm{id}\cup\mathrm{Res}} \ar@{^{(}->}[u]^-{\iota_0\otimes\mathrm{id}} \ar@{->>}[d]_-{\varphi_0\otimes((\rho_q^0)^{-1}\circ\pi)} & I \ar[dr]^{\varphi_{q+1}} \ar@{^{(}->}[ur]_{\iota_{q}} \\
H_p^1(k)\otimes \left( \displaystyle\bigoplus_{j=0}^{q} \left(\Lambda^j(\Gamma/p\Gamma)\otimes \kk_{q-j}(k)\right)\right) \ar@{->>}[rr]_-{((-1)^j\cup\psi_k^{q-j})_j} \ar@/_3pc/[rruu]^{\rho_1^e\cup\rho_q^0} && \displaystyle\bigoplus_{j=0}^{q} \left(\Lambda^j(\Gamma/p\Gamma)\otimes H_p^{q-j+1}(k)\right), \ar[uu]_{\rho_{q+1}^e}
}\]
commutes.
\end{proposition}

\begin{proof}
This is a direct consequence of Lemmas \ref{prop comm rho^e 1}, \ref{prop comm rho^e 2}, \ref{prop comm rho^e 3} and \ref{prop comm rho^e 4}.
\end{proof}

\begin{corollary}
    For every $q\geq 1$, the morphism $\rho^e_{q+1}$ is injective.
\end{corollary}

\begin{proof}
    Take $x\in \ker \rho^e_{q+1}$. There exists $\tilde{x}\in H^1(k,\Z/p\Z)\otimes\kk_q(K)$ such that $$x=((-1)^j\cup\psi_k^{q-j})_j(\varphi_0\otimes((\rho_q^0)^{-1}\circ\pi)(\tilde{x})).$$ We then have $0=\rho_{q+1}^e(x)=\iota_q(\mathrm{Id}\cup \mathrm{Res}(\tilde{x}))$. 
    By injectivity of $\iota_q$, we get $\mathrm{Id}\cup \mathrm{Res}(\tilde{x})=0$, and hence:
    $$x=\varphi_{q+1}(\mathrm{Id}\cup \mathrm{Res}(\tilde{x}))=0.$$
\end{proof}

\section{Maximally complete fields}\label{sec invisible}

We keep the notations from above. In particular, $(K,v)$ is a rank-one valued field of mixed characteristic $(0,p)$ with value group $\Gamma \subset \mathbb{R}$ and residue field $k$ of characteristic $p>0$, and the groups $u_q^{<\gamma_\infty}(K)$ and $u_q^{\geq\gamma_\infty}(K)$ are those defined in Section \ref{subsec setting}. This whole section is devoted to the proof of Theorem \ref{th invisible terms}, which we recall: \emph{if $(K,v)$ is maximally complete and $\gamma_\infty\leq e$, then $u_q^{<\gamma_\infty}(K)=u_q^{\geq\gamma_\infty}(K)$.}

\begin{notation}
Throughout this section, for a tuple of integers $\mathbf{i}=(i_1,\ldots,i_r)\in\Z^r$, we will denote by $0\leq \mathbf{i}\leq p-1$ to signal that every coordinate lies in that interval. We also denote by $\mathbf{0}$ the tuple $(0,\ldots,0)$.
\end{notation}

\begin{lemma}\label{lem transf ind}
Let $\bar \gamma=(\bar\gamma_1,\dots,\bar \gamma_g)$ be an $\mathbb{F}_p$-basis of $\Gamma/p\Gamma$ and let $\gamma_1,\ldots,\gamma_g\in\Gamma_{>0}$ be such that $v(\pi^{\bar\gamma_s})=\gamma_s$. Let $\bar\alpha_1,\dots,\bar\alpha_r$ be a $p$-basis of $k$ and let $\alpha_1,\ldots,\alpha_r$ denote lifts of these elements in $\cal O_K$. Assume that $(K,v)$ is maximally complete.
    
For tuples of integers $\mathbf{i}=(i_1,\dots,i_r)\in\Z^{r}$ and $\mathbf{j}=(j_1,\dots,j_g)\in\Z^g$, write $\alpha^{\mathbf{i}}:=\alpha_1^{i_1}\cdots\alpha_r^{i_r}$ and $\pi^{\mathbf{j}\bar\gamma}:=\pi^{j_1\bar\gamma_1}\cdots\pi^{j_g\bar\gamma_g}$. Then every $x\in \mathcal{O}_K$ is of the form:
$$x=\sum_{0\leq\mathbf{i},\mathbf{j}\leq p-1} \lambda_{\mathbf{i},\mathbf{j}}^p \alpha^{\mathbf{i}}\pi^{\mathbf{j}\bar\gamma}+\pi^ex_e,\quad \lambda_{\mathbf{i},\mathbf{j}}\in K, \quad x_e\in \mathcal{O}_K.$$
\end{lemma}

\begin{proof} We claim that it is enough to prove the statement under the following assumption on $x$:
\begin{equation}\label{assumption}\exists \lambda_{\mathbf{i},\mathbf{0}}^{(0)}\in \mathcal{O}_K, \;\;\; v\left( x-\sum_{0\leq\mathbf{i}\leq p-1} (\lambda_{\mathbf{i},\mathbf{0}}^{(0)})^p \alpha^{\mathbf{i}}\right)>\sum_{s=1}^g (p-1)\gamma_s.\end{equation}
Indeed, assume we know the lemma holds under \eqref{assumption}. Since $x\in\cal O_K$ and $\alpha$ is a $p$-basis of $k$, one can find $\lambda_{\mathbf{i},\mathbf{0}}^{(0)}\in \mathcal{O}_K$ (for $0\leq \mathbf{i} \leq p-1$) such that:
$$\gamma_0:=v\left( x-\sum_{0\leq\mathbf{i}\leq p-1} (\lambda_{\mathbf{i},\mathbf{0}}^{(0)})^p \alpha^{\mathbf{i}}\right)>0.$$
By density of $\Gamma$ in $\mathbb{R}$, one can also find $(c_s)_s\in (K^\times)^g$ such that $\gamma_s+pv(c_s)>0$ for each $1\leq s\leq g$ and:
\begin{equation}\label{assumption satisfied} \gamma_0> \sum_{s=1}^g (p-1)(\gamma_s+pv(c_s)).\end{equation}
Set $\varpi^{\bar \gamma_s}:=c_s^p\pi^{\bar \gamma_s}$. The valuations of the $\varpi^{\bar \gamma_s}$ are all positive (even though the valuations of the $c_s$ are typically negative) and lift a basis of $\Gamma/p\Gamma$. Hence the $\varpi^{\bar \gamma_s}$ satisfy the same hypothesis as the $\pi^{\bar \gamma_s}$. Since the lemma is assumed to hold under \eqref{assumption}, we know that $x$ is of the form:
$$x=\sum_{0\leq\mathbf{i},\mathbf{j}\leq p-1} \lambda_{\mathbf{i},\mathbf{j}}^p \alpha^{\mathbf{i}}\varpi^{\mathbf{j}\bar\gamma}+\pi^ex_e,\quad \lambda_{\mathbf{i},\mathbf{j}}\in K, \quad x_e\in \mathcal{O}_K.$$
Hence:
$$x=\sum_{0\leq\mathbf{i},\mathbf{j}\leq p-1} (\lambda_{\mathbf{i},\mathbf{j}}c_1^{j_1}\dots c_g^{j_g})^p \alpha^{\mathbf{i}}\pi^{\mathbf{j}\bar\gamma}+\pi^ex_e,$$
as wished.

We henceforth work under the assumption \eqref{assumption}. For all $\mathbf{i}\in\Z^r$ and $\mathbf{j}\in\Z^g$ with $0\leq\mathbf{i},\mathbf{j}\leq p-1$, we inductively construct the following sequences indexed by ordinals $\mu$:
    \begin{itemize}
        \item[\textbullet]  $(\gamma_\mu)_\mu$, with values in $\Gamma$,
        \item[\textbullet]  $(\lambda_{\mathbf{i},\mathbf{j}}^{(\mu)})_\mu$, with values in $ \mathcal{O}_K$,
        \item[\textbullet]  $(x_{\mu})_\mu$, with values in $ \mathcal{O}_K$,
    \end{itemize}
    and satisfying the following properties:
\begin{gather}
   x = \sum_{0\leq\mathbf{i},\mathbf{j}\leq p-1} (\lambda_{\mathbf{i},\mathbf{j}}^{(\mu)})^p \alpha^{\mathbf{i}}\pi^{\mathbf{j}\bar\gamma}+\pi^{\gamma_\mu}x_{\mu},\label{equ1}\\
   \text{the sequence $(\gamma_\mu)_\mu$ is strictly increasing,} \label{equ2}\\
\forall \mu'>\mu, \; v(\lambda_{\mathbf{i},\mathbf{j}}^{(\mu')}-\lambda_{\mathbf{i},\mathbf{j}}^{(\mu)})\in \left\{ \frac{\gamma_\mu-\mathbf{j}\bar\gamma}{p},+\infty\right\}, \label{equ3}
\end{gather}
where $\mathbf{j}\bar\gamma:=\sum_{s=1}^g j_s\gamma_s$. As we will see, this can be done until $\gamma_\mu\geq e$.\\

It is important to remark that, while the lemma states that $\lambda_{\mathbf{i},\mathbf{j}}\in K$, here, under assumption \eqref{assumption}, we are actually constructing a sequence of $\lambda_{\mathbf{i},\mathbf{j}}^{(\mu)}$ in $\cal O_K$. In fact, this hypothesis will be crucial in the induction step. It is when restricting from the general case to this particular case that we lose the control of the $\lambda_{\mathbf{i},\mathbf{j}}$ due to the (typically negative) valuation of the $c_s$.

\paragraph{Base case.} For $\mu=0$, we can simply take $\lambda^{(0)}_{\mathbf{i},\mathbf{j}}:=0$ for $\mathbf{j}\neq\mathbf{0}$ and $$\begin{cases}
\gamma_0=v\left( x-\sum_{0\leq\mathbf{i}\leq p-1}(\lambda_{\mathbf{i},\mathbf{0}}^{(0)})^p \alpha^{\mathbf{i}}\right) \\
x_0=\frac{x-\sum_{0\leq\mathbf{i}\leq p-1} (\lambda_{\mathbf{i},\mathbf{0}}^{(0)})^p \alpha^{\mathbf{i}}}{\pi^{\gamma_0}}.
\end{cases}$$
Note that, by assumption, we have $\gamma_0>\sum_{s=1}^g(p-1)\gamma_s\geq \mathbf{j}\bar\gamma$ for every $0\leq \mathbf{j}\leq p-1$.

\paragraph{Successor step.} Assume we have already constructed all the desired sequences up to the ordinal $\mu$, with $\gamma_\mu < e$. Now, since $\pi^\gamma(\pi^{\bar\gamma})^{-1}\in K^p$ for every $\gamma\in\Gamma$, we may write $\pi^{\gamma_\mu}$ as $\pi^{\mathbf{j}_\mu\bar\gamma}\rho_\mu^p$ for some $0\leq \mathbf{j}_\mu\leq p-1$ and some $\rho_\mu \in \mathcal{O}_K$ (recall that $\gamma_\mu>\gamma_0>\mathbf{j}_\mu\bar\gamma$). Write also $\bar x_\mu$ as $\sum_{0\leq \mathbf{i} \leq p-1} (\bar\xi_{\mathbf{i}}^{(\mu)})^p\bar\alpha^\mathbf{i}$. For each $\mathbf{i}$, denote by $\xi_{\mathbf{i}}^{(\mu)}$ a lifting of $\bar \xi_{\mathbf{i}}^{(\mu)}$ in $\mathcal{O}_K$ with $\xi_{\mathbf{i}}^{(\mu)}=0$ if $\bar \xi_{\mathbf{i}}^{(\mu)}=0$, and set: 
$$\lambda_{\mathbf{i},\mathbf{j}}^{(\mu+1)}:=\begin{cases}
    \lambda_{\mathbf{i},\mathbf{j}}^{(\mu)} & \text{if $\mathbf{j}\neq\mathbf{j}_\mu$}\\
    \lambda_{\mathbf{i},\mathbf{j}}^{(\mu)}+\rho_\mu \xi_{\mathbf{i}}^{(\mu)} & \text{if $\mathbf{j}=\mathbf{j}_\mu$.}
\end{cases}$$
Setting $K^{>\gamma_\mu}:=\{x\in K: v(x) >\gamma_\mu\}$, we get:
\begin{align*}
    x-&\sum_{0\leq\mathbf{i},\mathbf{j}\leq p-1} (\lambda_{\mathbf{i},\mathbf{j}}^{(\mu+1)})^p \alpha^{\mathbf{i}}\pi^{\mathbf{j}\bar\gamma}\\
    & = x-\sum_{\substack{0\leq\mathbf{i},\mathbf{j}\leq p-1\\\mathbf{j}\neq\mathbf{j}_\mu }} (\lambda_{\mathbf{i},\mathbf{j}}^{(\mu)})^p \alpha^{\mathbf{i}}\pi^{\mathbf{j}\bar\gamma} - \sum_{0\leq\mathbf{i}\leq p-1}(\lambda_{\mathbf{i},\mathbf{j}_\mu}^{(\mu)}+\rho_\mu \xi_{\mathbf{i}}^{(\mu)})^p \alpha^{\mathbf{i}}\pi^{\mathbf{j}_\mu\bar\gamma}\\
    &\equiv x-\sum_{\substack{0\leq\mathbf{i},\mathbf{j}\leq p-1\\\mathbf{j}\neq\mathbf{j}_\mu }} (\lambda_{\mathbf{i},\mathbf{j}}^{(\mu)})^p \alpha^{\mathbf{i}}\pi^{\mathbf{j}\bar\gamma} - \sum_{0\leq\mathbf{i}\leq p-1}((\lambda_{\mathbf{i},\mathbf{j}_\mu}^{(\mu)})^p+(\rho_\mu \xi_{\mathbf{i}}^{(\mu)})^p) \alpha^{\mathbf{i}}\pi^{\mathbf{j}_\mu\bar\gamma} &\pmod {K^{>\gamma_\mu}}\\
    \intertext{\raggedleft (since $\lambda_{\mathbf{i},\mathbf{j}}^{(\mu)}\in\cal O_K$ and $0<v(\rho_\mu)=\frac{1}{p}(\gamma_\mu-\mathbf{j}_\mu\bar\gamma)<\frac{\gamma_\mu}{p}<\frac{e}{p}$)}
    &\equiv \pi^{\gamma_\mu}x_\mu-\sum_{0\leq\mathbf{i}\leq p-1}(\rho_\mu \xi_{\mathbf{i}}^{(\mu)})^p \alpha^{\mathbf{i}}\pi^{\mathbf{j}_\mu\bar\gamma}&\pmod {K^{>\gamma_\mu}}\\
    &\equiv \pi^{\gamma_\mu}\left(x_\mu-\sum_{0\leq\mathbf{i}\leq p-1}(\xi_{\mathbf{i}}^{(\mu)})^p \alpha^{\mathbf{i}}\right)&\pmod{K^{>\gamma_\mu}}\\
    &\equiv 0 &\pmod{K^{>\gamma_\mu}}.
\end{align*}
Hence there exist $\gamma_{\mu+1}>\gamma_\mu$ and $x_{\mu+1}\in\mathcal{O}_K$ such that:
$$x = \sum_{0\leq\mathbf{i},\mathbf{j}\leq p-1} (\lambda_{\mathbf{i},\mathbf{j}}^{(\mu+1)})^p \alpha^{\mathbf{i}}\pi^{\mathbf{j}\bar\gamma}+\pi^{\gamma_{\mu+1}}x_{\mu+1}.$$
Moreover, since $v(\rho_\mu)=\frac{\gamma_\mu-\mathbf{j}_\mu\bar \gamma}{p}$, we have, for every $\mathbf{i}$ and every $\mathbf{j}$:
$$v(\lambda_{\mathbf{i},\mathbf{j}}^{(\mu+1)}-\lambda_{\mathbf{i},\mathbf{j}}^{(\mu)})\in \left\{ \frac{\gamma_\mu-\mathbf{j}\bar\gamma}{p},+\infty\right\}. $$

\paragraph{Limit step.} Let $\mu_\infty$ be a limit ordinal and assume that the sequences $(\gamma_\mu)_\mu$, $(\lambda^{(\mu)}_{\mathbf{i},\mathbf{j}})_\mu$, and  $(x_\mu)_\mu$ have all been already defined for $\mu<\mu_{\infty}$ and that $\gamma_\mu<e$ for every $\mu<\mu_\infty$.

Now let us inductively construct sequences $(\tilde \lambda_{\mathbf{i},\mathbf{j}}^{(\nu)})_\nu$ indexed by an ordinal $\nu$ in the following way. We first set $\tilde \lambda_{\mathbf{i},\mathbf{j}}^{(0)}:= \lambda_{\mathbf{i},\mathbf{j}}^{(0)}$. Assuming that $\nu$ is some ordinal such that $\tilde \lambda_{\mathbf{i},\mathbf{j}}^{(\nu')}$ has already been constructed for every $\nu'<\nu$, we set $\tilde \lambda_{\mathbf{i},\mathbf{j}}^{(\nu)}:=\tilde \lambda_{\mathbf{i},\mathbf{j}}^{(\nu'')}$ where $\nu''$ is the smallest ordinal such that $\tilde \lambda_{\mathbf{i},\mathbf{j}}^{(\nu'')}\neq \tilde \lambda_{\mathbf{i},\mathbf{j}}^{(\nu')}$ for all $\nu'<\nu$.
This definition makes sense because ordinals are well-ordered sets.

Now, since the sequences $(\gamma_\mu)_\mu$, $(\lambda^{(\mu)}_{\mathbf{i},\mathbf{j}})_\mu$, and  $(x_\mu)_\mu$ have only been already defined for $\mu<\mu_{\infty}$, there exist some ordinals $ \nu_{\mathbf{i},\mathbf{j},\infty}\leq \mu_\infty$ such that each of the sequences $(\tilde \lambda_{\mathbf{i},\mathbf{j}}^{(\nu)})_\nu$ is only defined for $\nu<\nu_{\mathbf{i},\mathbf{j},\infty}$. If $\nu_{\mathbf{i},\mathbf{j},\infty}=\nu_0+1$ is a successor ordinal, we set $\lambda_{\mathbf{i},\mathbf{j}}^{(\mu_\infty)}:=\tilde \lambda_{\mathbf{i},\mathbf{j}}^{(\nu_0)}$. If $\nu_{\mathbf{i},\mathbf{j},\infty}$ is a limit ordinal, by \eqref{equ3} the sequence $(\tilde \lambda_{\mathbf{i},\mathbf{j}}^{(\nu)})_\nu$ satisfies
\[\forall \nu'>\nu, \; v(\tilde\lambda_{\mathbf{i},\mathbf{j}}^{(\nu')}-\tilde\lambda_{\mathbf{i},\mathbf{j}}^{(\nu)})=\frac{\tilde\gamma_\nu-\mathbf{j}\bar\gamma}{p},\]
for some subsequence $(\tilde\gamma_\nu)_{\nu}$ of $(\gamma_\mu)_{\mu}$, which must be strictly increasing by \eqref{equ2}. We see then that $(\tilde \lambda_{\mathbf{i},\mathbf{j}}^{(\nu)})_\nu$ is a pseudo-Cauchy sequence, so we
let $\lambda_{\mathbf{i},\mathbf{j}}^{(\mu_\infty)}$ be a pseudo-limit of this sequence.

By construction and definition of pseudo-limits, we have:
$$\forall \mu< \mu_{\infty},\; v(\lambda_{\mathbf{i},\mathbf{j}}^{(\mu_\infty)}- \lambda_{\mathbf{i},\mathbf{j}}^{(\mu)})\in \left\{ \frac{\gamma_\mu-\mathbf{j}\bar\gamma}{p},+\infty\right\},$$
which settles \eqref{equ3} up to rank $\mu_{\infty}$. This also implies that, for $\mu< \mu_{\infty}$:
\[x - \sum_{0\leq\mathbf{i},\mathbf{j}\leq p-1} (\lambda_{\mathbf{i},\mathbf{j}}^{(\mu_\infty)})^p \alpha^{\mathbf{i}}\pi^{\mathbf{j}\bar\gamma} \equiv x - \sum_{0\leq\mathbf{i},\mathbf{j}\leq p-1} (\lambda_{\mathbf{i},\mathbf{j}}^{(\mu)})^p \alpha^{\mathbf{i}}\pi^{\mathbf{j}\bar\gamma}
 \equiv 0 \pmod{\pi^{\gamma_\mu}}.\]
Hence, if we set:
$$\gamma_{\mu_\infty}:= v\left( x - \sum_{0\leq\mathbf{i},\mathbf{j}\leq p-1} (\lambda_{\mathbf{i},\mathbf{j}}^{(\mu_\infty)})^p \alpha^{\mathbf{i}}\pi^{\mathbf{j}\bar\gamma} \right), $$
we deduce that $\gamma_{\mu_\infty}>\gamma_\mu$ for every $\mu< \mu_{\infty}$ and that there exists $x_{\mu_\infty}\in \mathcal{O}_K$ such that: $$x = \sum_{0\leq\mathbf{i},\mathbf{j}\leq p-1} (\lambda_{\mathbf{i},\mathbf{j}}^{(\mu_\infty)})^p \alpha^{\mathbf{i}}\pi^{\mathbf{j}\bar\gamma}+\pi^{\gamma_{\mu_\infty}}x_{\mu_\infty},$$
thus establishing assertions \eqref{equ1} and \eqref{equ2} at rank $\mu_\infty$.

\paragraph{Termination of the inductive process.} Since $\Gamma$ is a set, our inductive construction of the sequences $(\gamma_\mu)_\mu$, $(\lambda^{(\mu)}_{\mathbf{i},\mathbf{j}})_\mu$, and  $(x_\mu)_\mu$ has to stop. This happens when $\gamma_{\mu_\infty}\geq e$ for some ordinal $\mu_\infty$. We then have:
$$x \equiv \sum_{0\leq\mathbf{i},\mathbf{j}\leq p-1} (\lambda_{\mathbf{i},\mathbf{j}}^{(\mu_\infty)})^p \alpha^{\mathbf{i}}\pi^{\mathbf{j}\bar\gamma}\pmod{\pi^{e}}.$$
\end{proof}

\begin{lemma}\label{lem common slots}
    Keep the previous assumptions and notations. Given a subset $S$ of $\{1,2,\dots,g\}$ and a subset $T$ of $\{1,2,\dots,r\}$, we set $\pi^{\bar\gamma_S}:=(\pi^{\bar\gamma_s})_{s\in S}$ and $ \alpha_T:=(\alpha_t)_{t\in T}$.  For $a\in u^\gamma_q(K)$, there exist $x_{S,T}\in \mathcal{O}_K$ such that:
    $$a\equiv\sum_{|S|+|T|=q-1} \{\pi^{\bar\gamma_S},1+\pi^{\gamma}x_{S,T},\alpha_T\} \pmod {u^{\gamma'}_q(K)},$$
    with $\gamma'=\min \{2\gamma,e\}$.
\end{lemma}

\begin{proof}
By Lemma \ref{lem generateurs u_q}, we may and do assume that $a$ is of the form $\{\pi^{\bar\gamma_S},1+\pi^\gamma x,y_1,\dots,y_{q-1-|S|}\}$ for some $S\subset\{1,\dots,g\}$, $x\in \mathcal{O}_K$ and $y_1,\dots, y_{q-1-|S|}\in \mathcal{O}_K^\times$. By Lemma \ref{lem transf ind}, we can find elements $\lambda_{\mathbf{i},\mathbf{j}} \in K$ such that: $$y_1\equiv\sum_{0\leq\mathbf{i},\mathbf{j}\leq p-1 } \lambda_{\mathbf{i},\mathbf{j}}^p \alpha^\mathbf{i}\pi^{\mathbf{j}\bar \gamma} \pmod {\pi^e},$$
where $\alpha^\mathbf{i}=\alpha_1^{i_1}\dots\alpha_r^{i_r}$ and $\pi^{\mathbf{j}\bar \gamma}=\pi^{j_1\bar \gamma_1}\dots \pi^{j_g\bar \gamma_g}$. Since $y_1\in \mathcal{O}_K$ and $\bar \gamma$ is a basis of $\Gamma/p\Gamma$, for every $1\leq\mathbf{j}\leq p-1$, $\left(\sum_{0\leq\mathbf{i}\leq p-1 } \lambda_{\mathbf{i},\mathbf{j}}^p \alpha^\mathbf{i}\right)\pi^{\mathbf{j}\bar \gamma}$ is in $\mathcal{O}_K$. But since $\bar \alpha$ is a $p$-basis of $k$, $v\left(\sum_{0\leq\mathbf{i}\leq p-1 } \lambda_{\mathbf{i},\mathbf{j}}^p \alpha^\mathbf{i}\right)=\min \{v(\lambda_{\mathbf{i},\mathbf{j}}^p \alpha^\mathbf{i}): 0\leq \mathbf{i} \leq p-1\}$ and hence, for every $0\leq \mathbf{i},\mathbf{j}\leq p-1$, $\lambda_{\mathbf{i},\mathbf{j}}^p \alpha^\mathbf{i}\pi^{\mathbf{j}\bar \gamma}$ lies in $\mathcal{O}_K$ (even though $\lambda_{\mathbf{i},\mathbf{j}}$ might not).

Applying Lemma \ref{lem sum symb}, we get the following equalities modulo $u_q^{\gamma'}(K)$:
\begin{align*}
    a &= \{\pi^{\bar\gamma_S},1+\pi^\gamma (xy_1^{-1})y_1,y_1,\dots,y_{q-1-|S|}\}\\
    &\equiv\sum_{\substack{0\leq\mathbf{i},\mathbf{j}\leq p-1\\ \lambda_{\mathbf{i},\mathbf{j}}\neq 0}}  \{\pi^{\bar\gamma_S},1+\pi^\gamma (xy_1^{-1})\lambda_{\mathbf{i},\mathbf{j}}^p \alpha^\mathbf{i}\pi^{\mathbf{j}\bar \gamma},\lambda_{\mathbf{i},\mathbf{j}}^p \alpha^\mathbf{i}\pi^{\mathbf{j}\bar \gamma},y_2,\dots,y_{q-1-|S|}\}\\
    &\equiv\sum_{0\leq\mathbf{i},\mathbf{j}\leq p-1 }  \{\pi^{\bar\gamma_S},1+\pi^\gamma x_{\mathbf{i},\mathbf{j}},\alpha^\mathbf{i}\pi^{\mathbf{j}\bar \gamma},y_2\dots,y_{q-1-|S|}\}\pmod{u_q^{\gamma'}(K)},
\end{align*}
for some $x_{\mathbf{i},\mathbf{j}}\in \mathcal{O}_K$. 
    
Repeating the previous argument for $y_2,\dots,y_{q-1-|S|}$ instead of $y_1$, we see that $a \pmod{u^{\gamma'}_q(K)}$ can be written as a sum of symbols of the form: $$\{\pi^{\bar\gamma_S},1+\pi^\gamma x',\alpha^{\mathbf{i}^{(1)}}\pi^{\mathbf{j}^{(1)}\bar \gamma},\dots,\alpha^{\mathbf{i}^{(q-1-|S|)}}\pi^{\mathbf{j}^{(q-1-|S|)}\bar \gamma}\}$$
for some $\mathbf{i}^{(1)},\ldots,\mathbf{i}^{(q-1-|S|)}\in\Z^r$, $\mathbf{j}^{(1)},\ldots,\mathbf{j}^{(q-1-|S|)}\in\Z^g$ and $x'\in\mathcal{O}_K$. But each of those symbols is itself a sum of symbols of the form $\{\pi^{\bar\gamma_S},1+\pi^{\gamma}x'',\alpha_T\}$ for some $S\subset \{1,2,\dots,g\}$ and $T\subset \{1,2,\dots,r\}$ with $|S|+|T|=q-1$ and some $x''\in \mathcal{O}_K$. The result immediately follows.
\end{proof}

\begin{lemma}\label{lem annulation}
    Fix $\gamma\in p\Gamma$ and set $\gamma_{\infty}=\min \{2\gamma,v(p)+\gamma\}$. Fix $t\in\{0,\dots r\}$. For each subset $T$ of $\{1,\dots,r\}$ with $t$ elements, each subset $U$ of $\{1,\dots,r\}$ with $t-1$ elements and each $\mathbf{i}\in\Z^r$ with $0\leq\mathbf{i}\leq p-1$, take elements $\lambda_{T,\mathbf{i}}$ and $\rho_{U,\mathbf{i}}$ in $\mathcal{O}_K$ such that:
$$\lambda_{T,\mathbf{i}} = \sum_{j\in T} (-1)^{t(j)+1}i_j\rho_{T\setminus\{j\},\mathbf{i}},$$
    where $t(j)$ stands for the position of $j$ in $T$. Then:
    $$\sum_{|T|=t}\{1+\pi^{\gamma}\sum_{0\leq\mathbf{i}\leq p-1}\lambda_{T,\mathbf{i}}^p  \alpha^\mathbf{i}, \alpha_T\}\in u^{\gamma_{\infty}}_{t+1}(K).$$
\end{lemma}

\begin{proof} Since $v(p)+\gamma\geq\gamma_\infty$, we have
\[\pi^\gamma\lambda_{T,\mathbf{i}}^p\equiv \pi^\gamma \sum_{j\in T} (-1)^{p(t(j)+1)}i_j^p\rho_{T\setminus\{j\},\mathbf{i}}^p\pmod{\pi^{\gamma_\infty}}.\]
We compute then modulo $u^{\gamma_{\infty}}_{t+1}(K)$.
    \begin{align*}
        \sum_{|T|=t}\{&1+\pi^{\gamma}\sum_{0\leq\mathbf{i}\leq p-1}\lambda_{T,\mathbf{i}}^p  \alpha^\mathbf{i}, \alpha_T\} \\
        &\equiv \sum_{|T|=t}\{1+\pi^{\gamma}\sum_{0\leq\mathbf{i}\leq p-1}\left(
        \sum_{j\in T} (-1)^{p(t(j)+1)}i_j^p\rho_{T\setminus\{j\},\mathbf{i}}^p  \right)\alpha^\mathbf{i},\alpha_T\}\pmod {u^{\gamma_{\infty}}_{t+1}(K)},  \\
        \intertext{and using that $2\gamma\geq \gamma_{\infty}$, we get}
        &\equiv \sum_{|T|=t}\sum_{0\leq\mathbf{i}\leq p-1}
        \sum_{j\in T} (-1)^{p(t(j)+1)}i_j^p\{1+\pi^{\gamma}
        \rho_{T\setminus\{j\},\mathbf{i}}^p  \alpha^\mathbf{i},\alpha_T\} \\
         &\equiv \sum_{|T|=t}\sum_{0\leq\mathbf{i}\leq p-1}
        \sum_{j\in T} (-1)^{t(j)+1}i_j\{1+\pi^{\gamma}
        \rho_{T\setminus\{j\},\mathbf{i}}^p  \alpha^\mathbf{i},\alpha_T\} \\
        &\equiv \sum_{|T|=t}\sum_{0\leq\mathbf{i}\leq p-1}
        \sum_{j\in T} i_j\{1+\pi^{\gamma}
        \rho_{T\setminus\{j\},\mathbf{i}}^p  \alpha^\mathbf{i},\alpha_j,\alpha_{T\setminus\{j\}}\} \\
        &\equiv \sum_{|T|=t}\sum_{0\leq\mathbf{i}\leq p-1}
        \sum_{j\in T} \{1+\pi^{\gamma}
        \rho_{T\setminus\{j\},\mathbf{i}}^p  \alpha^\mathbf{i},\alpha_j^{i_j},\alpha_{T\setminus\{j\}}\} \\
        &\equiv \sum_{|U|=t-1}\sum_{0\leq\mathbf{i}\leq p-1}
        \sum_{j\not\in U} \{1+\pi^{\gamma}
        \rho_{U,\mathbf{i}}^p  \alpha^\mathbf{i},\alpha_j^{i_j},\alpha_{U}\}\pmod {u^{\gamma_{\infty}}_{t+1}(K)} \\
        \intertext{but since $\{\alpha,\alpha\}=\{\alpha,-1\}$ for any $\alpha\in K^\times$, $-1$ is a $p$-th power for $p\neq 2$ and $-1=1-p$ for $p=2$, we see by Lemma \ref{lem suma} that $\{\alpha_j,\alpha_U\}\in u_t^{v(p)}$ if $j\in U$, and hence}
        &\equiv \sum_{|U|=t-1}\sum_{0\leq\mathbf{i}\leq p-1}
        \sum_{j\in \{1,\dots,r\}} \{1+\pi^{\gamma}
        \rho_{U,\mathbf{i}}^p  \alpha^\mathbf{i},\alpha_j^{i_j},\alpha_{U}\} \\
        &\equiv \sum_{|U|=t-1}\sum_{0\leq\mathbf{i}\leq p-1}
         \{1+\pi^{\gamma}
        \rho_{U,\mathbf{i}}^p  \alpha^\mathbf{i},\alpha^\mathbf{i},\alpha_{U}\} \\
        &\equiv \sum_{|U|=t-1}\sum_{0\leq\mathbf{i}\leq p-1}
         \{1+\pi^{\gamma}
        \rho_{U,\mathbf{i}}^p  \alpha^\mathbf{i},(\pi^{\gamma}
        \rho_{U,\mathbf{i}}^p)^{-1},\alpha_{U}\}\pmod {u^{\gamma_{\infty}}_{t+1}(K)},
    \end{align*}
and since $\gamma\in p\Gamma$, $\pi^\gamma$ is a $p$-th power and hence this last sum is equal to $0$.
\end{proof}

\begin{proof}[Proof of Theorem \ref{th invisible terms}]
Take some $a\in u_q^{<\gamma_\infty}(K)$ and fix some $\gamma\in p\Gamma$ with:
$$\max\left\{\frac{\gamma_{\infty}}{2},\gamma_{\infty}-v(p)\right\} < \gamma< \gamma_\infty.$$
Since $a\in u_q^{\gamma}(K)$, Lemma \ref{lem common slots} shows that:
\begin{equation}\label{eq0}
a\equiv\sum_{|S|+|T|=q-1} \{\pi^{\bar\gamma_S},1+\pi^{\gamma} x_{S,T},\alpha_T\} \pmod {u^{\gamma_\infty}_q(K)}.
\end{equation}
for some $x_{S,T} \in \mathcal{O}_K$. 
    
    For each $\mathbf{i}\in\Z^r\setminus\{\mathbf{0}\}$ with $0\leq\mathbf{i}\leq p-1$, each $S\subset \{1,\dots,g\}$ and each $T,U\subset \{1,\dots,r\}$ with $|S|+|T|=q-1$ and $|T|=|U|+1$, we inductively construct the following sequences indexed by ordinals $\mu$:
    \begin{itemize}
        \item[\textbullet]  $(\theta_\mu)_\mu$, with values in $\mathcal{O}_K^\times$,
        \item[\textbullet]  $(\gamma_\mu)_\mu$, with values in $\Gamma$,
        \item[\textbullet]  $(x_{S,T}^{(\mu)})_\mu$, with values in $ \mathcal{O}_K$,
        \item[\textbullet]  $(\lambda^{(\mu)}_{S,T,\mathbf{i}})_\mu$, with values in $\mathcal{O}_K$,
        \item[\textbullet] $(\rho_{S,U,\mathbf{i}}^{(\mu)})_\mu$, with values in $\mathcal{O}_K$,
    \end{itemize}
    and satisfying the following properties:
\begin{gather}
   1+\pi^{\gamma} x_{S,T} = \theta_\mu^p \left(1+\pi^{\gamma}\sum_{\substack{0\leq\mathbf{i}\leq p-1\\ \mathbf{i}\neq \mathbf{0}}}(\lambda^{(\mu)}_{S,T,\mathbf{i}})^p  \alpha^\mathbf{i}\right) (1 + \pi^{\gamma_\mu} x_{S,T}^{(\mu)}),\label{eq1}\\
   \lambda^{(\mu)}_{S,T,\mathbf{i}} = \sum_{j\in T} (-1)^{t(j)+1}i_j\rho_{S,T\setminus \{j\},\mathbf{i}}^{(\mu)},\label{eq2}\\
   \text{the sequence $(\gamma_\mu)_\mu$ is strictly increasing,} \label{eq3}\\
\forall \mu'>\mu, \; v(\theta_{\mu'}-\theta_{\mu})\in \left\{ \frac{\gamma_\mu}{p},+\infty\right\}, \label{eq4}\\
\forall \mu'>\mu, \; v(
\rho_{S,U,\mathbf{i}}^{(\mu')}-\rho_{S,U,\mathbf{i}}^{(\mu)}) \in \left\{ \frac{\gamma_\mu-\gamma}{p},+\infty\right\}. \label{eq5}
\end{gather}
As we will see, this can be done until $\gamma_\mu\geq \gamma_{\infty}$.

\paragraph{Base case.} For $\mu=0$, we can simply take $\theta_0=1$, $\gamma_0=\gamma$, $x_{S,T}^{(0)}=x_{S,T}$, $\lambda_{S,T,\mathbf{i}}^{(0)}=0$ and $\rho_{S,U,\mathbf{i}}^{(0)}=0$.

\paragraph{Successor step.} Assume we have already constructed all the desired sequences up to the ordinal $\mu$, with $\gamma_\mu < \gamma_\infty$. By \eqref{eq1} and Lemma \ref{lem annulation} (which uses \eqref{eq2}), the congruence \eqref{eq0} becomes:
$$a\equiv\sum_{|S|+|T|=q-1} \{\pi^{\bar\gamma_S},1+\pi^{\gamma_\mu} x_{S,T}^{(\mu)},\alpha_T\} \pmod {u^{\gamma_\infty}_q(K)}.$$
Now, for each $S$ and $T$, write:
\begin{equation}\label{eq6}
\bar x^{(\mu)}_{S,T} = \sum_{0\leq\mathbf{i}\leq p-1} (\bar\xi^{(\mu+1)}_{S,T,\mathbf{i}})^p \bar \alpha^\mathbf{i}
\end{equation}
    for some $\bar\xi^{(\mu+1)}_{S,T,\mathbf{i}}\in k$. 
    We distinguish two cases.
    \begin{enumerate}
    \item \emph{Case 1: $p\nmid \gamma_\mu$.} Since $a\in u_q^{<\gamma_\infty}(K)$ and $\gamma_\mu<\gamma_\infty$, by Theorem \ref{thm isom Kato}:
    $$\forall S\subset \{1,\dots, g\}, \sum_{\substack{|T|=q-1-|S|\\ 0\leq\mathbf{i}\leq p-1}} (\bar\xi^{(\mu+1)}_{S,T,\mathbf{i}})^p \bar \alpha^\mathbf{i}\frac{d\bar\alpha_T}{\bar\alpha_T}=0,$$
    and hence $\bar \xi^{(\mu+1)}_{S,T,\mathbf{i}}=0$ for each $S$, $T$ and $\mathbf{i}$. We may therefore write $\pi^{\gamma_\mu}x_{S,T}^{(\mu)}=\pi^{\gamma_{\mu+1}}x_{S,T}^{(\mu+1)}$ for some $\gamma_{\mu+1}>\gamma_\mu$ and some $x_{S,T}^{(\mu+1)}\in \mathcal{O}_K$, and we set $\theta_{\mu+1}=\theta_\mu$,  $\lambda_{S,T,\mathbf{i}}^{(\mu+1)}=\lambda_{S,T,\mathbf{i}}^{(\mu)}$ and $\rho_{S,U,\mathbf{i}}^{(\mu+1)}=\rho_{S,U,\mathbf{i}}^{(\mu)}$ for each $S,T,U,\mathbf{i}$. 
        \item \emph{Case 2: $p\mid \gamma_\mu$.} As before, by Theorem \ref{thm isom Kato} and by Cartier's Theorem \cite[Thm.~9.4.3]{GS}:
    \begin{align*}
&\forall S\subset \{1,\dots, g\},
    \\
    &\sum_{\substack{|T|=q-1-|S|\\ 0\leq\mathbf{i}\leq p-1}} (\bar\xi^{(\mu+1)}_{S,T,\mathbf{i}})^p \bar \alpha^\mathbf{i}\frac{d\bar\alpha_T}{\bar\alpha_T}\in \ker d = \sum_{\substack{|U|=q-2-|S|\\ 0\leq\mathbf{i}\leq p-1}} k^p d\left(\bar\alpha^{\mathbf{i}}\frac{d\bar\alpha_U}{\bar\alpha_U}\right) + \sum_{|T|=q-1-|S|} k^p\frac{d\bar\alpha_T}{\bar\alpha_T}.\end{align*}
    We can therefore find $\bar\chi_{S,U,\mathbf{i}}^{(\mu+1)}$ and $\bar\chi_{S,T}^{(\mu+1)}$ in $k$ such that:
    \begin{align*}
&\forall S\subset \{1,\dots, g\},
    \\
    &\sum_{\substack{|T|=q-1-|S|\\ 0\leq\mathbf{i}\leq p-1}} (\bar\xi^{(\mu+1)}_{S,T,\mathbf{i}})^p \bar \alpha^\mathbf{i}\frac{d\bar\alpha_T}{\bar\alpha_T} = \sum_{\substack{|U|=q-2-|S|\\ 0\leq\mathbf{i}\leq p-1}} (\bar\chi_{S,U,\mathbf{i}}^{(\mu+1)})^p d\left(\bar\alpha^{\mathbf{i}}\frac{d\bar\alpha_U}{\bar\alpha_U}\right) + \sum_{|T|=q-1-|S|} (\bar\chi_{S,T}^{(\mu+1)})^p\frac{d\bar\alpha_T}{\bar\alpha_T}.\end{align*}
    The right-hand side of this equality can be rewritten as:
    \begin{align*}\sum_{\substack{|U|=q-2-|S|\\ 0\leq\mathbf{i}\leq p-1}} &(\bar\chi_{S,U,\mathbf{i}}^{(\mu+1)})^p\bar\alpha^\mathbf{i}\left( \sum_{j=1}^r i_j \frac{d\alpha_j}{\alpha_j} \right) \wedge \frac{d\bar\alpha_U}{\bar\alpha_U} + \sum_{|T|=q-1-|S|} (\bar\chi_{S,T}^{(\mu+1)})^p\frac{d\bar\alpha_T}{\bar\alpha_T}\\
    &=\sum_{|T|=q-1-|S|} \left((\bar\chi_{S,T}^{(\mu+1)})^p + \sum_{j\in T,\, \mathbf{i}\neq\mathbf{0}} (-1)^{t(j)+1}i_j(\bar\chi_{S,T\setminus \{j\},\mathbf{i}}^{(\mu+1)})^p\bar\alpha^\mathbf{i}\right)\frac{d\bar\alpha_T}{\bar \alpha_T}.
    \end{align*}
    Hence, for $S\subset \{1,\dots, g\}$ and $T\subset \{1,\dots,r\}$ with $|S|+|T|=q-1$:
    $$
    \begin{cases}
    (\bar\xi^{(\mu+1)}_{S,T,\mathbf{0}})^p = (\bar\chi_{S,T}^{(\mu+1)})^p,\\    
    \forall \,\mathbf{i}\neq \mathbf{0}, \;(\bar\xi^{(\mu+1)}_{S,T,\mathbf{i}})^p = \sum_{j\in T} (-1)^{t(j)+1}i_j(\bar\chi_{S,T\setminus \{j\},\mathbf{i}}^{(\mu+1)})^p.
    \end{cases}
    $$
    Equivalently:
    $$
    \begin{cases}
    \bar\xi^{(\mu+1)}_{S,T,\mathbf{0}} = \bar\chi_{S,T}^{(\mu+1)},\\    
    \forall \mathbf{i}\neq \mathbf{0}, \;\bar\xi^{(\mu+1)}_{S,T,\mathbf{i}} = \sum_{j\in T} (-1)^{t(j)+1}i_j\bar\chi_{S,T\setminus \{j\},\mathbf{i}}^{(\mu+1)}.
    \end{cases}
    $$
    For $\mathbf{i}\neq \mathbf{0}$ and for $S\subset \{1,\dots, g\}$ and $T,U\subset \{1,\dots,r\}$ with $|S|+|T|=q-1$ and $|T|=|U|+1$, we may then consider liftings $\xi^{(\mu+1)}_{S,T,\mathbf{i}}$ and $\chi_{S,U,\mathbf{i}}^{(\mu+1)}$ of $\bar \xi^{(\mu+1)}_{S,T,\mathbf{i}}$ and $\bar \chi_{S,U,\mathbf{i}}^{(\mu+1)}$ respectively such that:
    $$   
   \begin{cases}
   \chi_{S,U,\mathbf{i}}^{(\mu+1)}=0 \Leftrightarrow \bar \chi_{S,U,\mathbf{i}}^{(\mu+1)}=0\\
    \xi^{(\mu+1)}_{S,T,\mathbf{i}} = \sum_{j\in T} (-1)^{t(j)+1}i_j\chi_{S,T\setminus \{j\},\mathbf{i}}^{(\mu+1)}.
    \end{cases}
    $$
    We set:
    $$\begin{cases}
    \rho_{S,U,\mathbf{i}}^{(\mu+1)}:=\rho_{S,U,\mathbf{i}}^{(\mu)}+\left(\frac{\pi^{\gamma_\mu}}{\pi^\gamma}\right)^{1/p}\chi_{S,U,\mathbf{i}}^{(\mu+1)}\\
    \lambda_{S,T,\mathbf{i}}^{(\mu+1)}:=\lambda_{S,T,\mathbf{i}}^{(\mu)}+\left(\frac{\pi^{\gamma_\mu}}{\pi^\gamma}\right)^{1/p}\xi_{S,T,\mathbf{i}}^{(\mu+1)}.
    \end{cases}$$

    Moreover, if $\xi^{(\mu+1)}_{S,T,\mathbf{0}}=\chi_{S,T}^{(\mu+1)}$ stands for a lifting of $\bar \xi^{(\mu+1)}_{S,T,\mathbf{0}}=\bar \chi_{S,T}^{(\mu+1)}$ in $\mathcal{O}_K$ with $\xi^{(\mu+1)}_{S,T,\mathbf{0}}=\chi_{S,T}^{(\mu+1)}=0$ whenever $\bar \chi_{S,T}^{(\mu+1)}=0$, we set: $$\theta_{\mu+1}=\theta_\mu \left(1+(\pi^{\gamma_\mu})^{1/p}\chi_{S,T}^{(\mu+1)}\right).$$

    Recalling that $\gamma_\mu>\gamma>\gamma_\infty-v(p)>\gamma_\mu-v(p)$, we then have the following equalities in the quotient $K^\times/U_1^{>\gamma_\mu}(K)$:
    \begin{align*}
        \theta_{\mu+1}^p &\left(1+\pi^{\gamma}\sum_{\mathbf{i}\neq \mathbf{0}}(\lambda^{(\mu+1)}_{S,T,\mathbf{i}})^p  \alpha^\mathbf{i}\right) \\
        & = \theta_\mu^p \left(1+(\pi^{\gamma_\mu})^{1/p}\chi_{S,T}^{(\mu+1)}\right)^p \left(1+\pi^{\gamma}\sum_{\mathbf{i}\neq \mathbf{0}}\left(\lambda_{S,T,\mathbf{i}}^{(\mu)}+\left(\frac{\pi^{\gamma_\mu}}{\pi^\gamma}\right)^{1/p}\xi_{S,T,\mathbf{i}}^{(\mu+1)}\right)^p  \alpha^\mathbf{i}\right) \\
        &= \theta_\mu^p \left( 1+\pi^{\gamma_\mu}(\chi_{S,T}^{(\mu+1)})^p\right) \left(1+\pi^{\gamma}\sum_{\mathbf{i}\neq \mathbf{0}}\left((\lambda_{S,T,\mathbf{i}}^{(\mu)})^p+\frac{\pi^{\gamma_\mu}}{\pi^\gamma}(\xi_{S,T,\mathbf{i}}^{(\mu+1)})^p\right) \alpha^\mathbf{i}\right) \\
        &= \theta_\mu^p \left(1+\pi^{\gamma}\sum_{\mathbf{i}\neq \mathbf{0}}(\lambda_{S,T,\mathbf{i}}^{(\mu)})^p \alpha^\mathbf{i}\right)\left( 1+\pi^{\gamma_\mu}(\xi_{S,T,\mathbf{0}}^{(\mu+1)})^p\right) \left(1+\pi^{\gamma_\mu}\sum_{\mathbf{i}\neq \mathbf{0}}(\xi_{S,T,\mathbf{i}}^{(\mu+1)})^p \alpha^\mathbf{i}\right) \\
        &= (1+\pi^{\gamma} x_{S,T}) (1+\pi^{\gamma_{\mu}} x_{S,T}^{(\mu)})^{-1}  \left(1+\pi^{\gamma_\mu}\sum_{\mathbf{i}}(\xi_{S,T,\mathbf{i}}^{(\mu+1)})^p \alpha^\mathbf{i}\right) \\
        &= 1+\pi^{\gamma} x_{S,T} \in K^\times/U_1^{>\gamma_\mu}(K),
    \end{align*}
where the last step follows from \eqref{eq6}. Hence we can find $\gamma_{\mu+1}>\gamma_\mu$ and $x_{S,T}^{(\mu+1)} \in \mathcal{O}_K$ such that:
$$1+\pi^{\gamma} x_{S,T} = \theta_{\mu+1}^p \left(1+\pi^{\gamma}\sum_{\mathbf{i}\neq \mathbf{0}}(\lambda^{(\mu+1)}_{S,T,\mathbf{i}})^p  \alpha^\mathbf{i}\right) (1 + \pi^{\gamma_{\mu+1}} x_{S,T}^{(\mu+1)}),$$
settling equation \eqref{eq1} at rank $\mu+1$.
    
    Moreover, for every $S,T$ and $\mathbf{i}\neq\mathbf{0}$, the equalities :
    $$\lambda^{(\mu+1)}_{S,T,\mathbf{i}} = \sum_{j\in T} (-1)^{t(j)+1}i_j\rho_{S,T\setminus \{j\},\mathbf{i}}^{(\mu+1)}$$
    and the inequalities:
    \[\forall\, \nu\leq\mu,\quad v(\theta_{\mu+1}-\theta_{\nu})\in \left\{ \frac{\gamma_{\nu}}{p}, +\infty\right\}, \quad v(
\rho_{S,U,\mathbf{i}}^{(\mu+1)}-\rho_{S,U,\mathbf{i}}^{(\nu)})\in \left\{ \frac{\gamma_\nu-\gamma}{p},+\infty\right\}.\]
    are obvious by construction, so equations \eqref{eq2}, \eqref{eq4} and \eqref{eq5} hold up to rank $\mu+1$.
    \end{enumerate}
    
\paragraph{Limit step.} Let $\mu_\infty$ be a limit ordinal and assume that the sequences $(\theta_\mu)_\mu$, $(\gamma_\mu)_\mu$, $(x_{S,T}^{(\mu)})_\mu$, $(\lambda^{(\mu)}_{S,T,\mathbf{i}})_\mu$, and  $(\rho_{S,U,\mathbf{i}}^{(\mu)})_\mu$ have all been already defined for $\mu<\mu_{\infty}$ and that $\gamma_\mu<\gamma_\infty$ for every $\mu<\mu_\infty$.

Now let us inductively construct sequences $(\tilde \theta_\nu)_\nu$ and $(\tilde \rho_{S,U,\mathbf{i}}^{(\nu)})_\nu$ indexed by an ordinal $\nu$ in the following way. We first set $\tilde \theta_0:= \theta_0$ and $\tilde \rho_{S,U,\mathbf{i}}^{(0)}:=\rho_{S,U,\mathbf{i}}^{(0)}$. Assuming that $\nu$ is some ordinal such that $\tilde \theta_{\nu'}$ and $\tilde \rho_{S,U,\mathbf{i}}^{(\nu')}$ have already been constructed for every $\nu'<\nu$, we set:
\begin{itemize}
    \item[\textbullet] $\tilde \theta_{\nu}:=\theta_{\nu''}$ where $\nu''$ is the smallest ordinal such that $\theta_{\nu''}\neq \theta_{\nu'}$ for all $\nu'<\nu$;
    \item[\textbullet] $\tilde \rho_{S,U,\mathbf{i}}^{(\nu)}:=\rho_{S,U,\mathbf{i}}^{(\nu'')}$ where $\nu''$ is the smallest ordinal such that $\rho_{S,U,\mathbf{i}}^{(\nu'')}\neq \rho_{S,U,\mathbf{i}}^{(\nu')}$ for all $\nu'<\nu$.
\end{itemize}
These definitions make sense because ordinals are well-ordered sets.

Now, since the sequences $(\theta_\mu)_\mu$, $(\gamma_\mu)_\mu$, $(x_{S,T}^{(\mu)})_\mu$, $(\lambda^{(\mu)}_{S,T,\mathbf{i}})_\mu$ and  $(\rho_{S,U,\mathbf{i}}^{(\mu)})_\mu$ have only been already defined for $\mu<\mu_{\infty}$, there exist some ordinals $\nu_\infty, \nu_{S,U,\mathbf{i},\infty}\leq \mu_\infty$ such that the sequence $(\tilde \theta_\nu)_\nu$ is only defined for $\nu<\nu_\infty$ and each of the sequences $(\tilde \rho_{S,U,\mathbf{i}}^{(\nu)})_\nu$ is only defined for $\nu<\nu_{S,U,\mathbf{i},\infty}$. If $\nu_{\infty}=\nu_0+1$ is a successor ordinal, we set $\theta_{\mu_{\infty}}:=\tilde\theta_{\nu_0}$. If $\nu_{\infty}$ is a limit ordinal, by \eqref{eq4} the sequence $(\tilde \theta_\nu)_\nu$ satisfies
\[\forall \nu'>\nu, \; v(\tilde\theta_{\nu'}-\tilde\theta_{\nu})=\frac{\tilde\gamma_\nu}{p},\]
for some subsequence $(\tilde\gamma_\nu)_{\nu}$ of $(\gamma_\mu)_{\mu}$, which must be strictly increasing by \eqref{eq3}. We see then that $(\tilde \theta_\nu)_\nu$ is a pseudo-Cauchy sequence, so we let $\theta_{\mu_{\infty}}$ be a pseudo-limit of this sequence. Similarly, for each $S$, $U$ and $\mathbf{i}$, if $\nu_{S,U,\mathbf{i},\infty}=\nu_0+1$ is a successor ordinal, we set $\rho_{S,U,\mathbf{i}}^{(\mu_{\infty})}:=\tilde\rho_{S,U,\mathbf{i}}^{(\nu_0)}$. If $\nu_{S,U,\mathbf{i},\infty}$ is a limit ordinal, the sequence $(\tilde \rho_{S,U,\mathbf{i}}^{(\nu)})_\nu$ is a pseudo-Cauchy sequence by \eqref{eq5} and \eqref{eq3}, and we let $\rho_{S,U,\mathbf{i}}^{(\mu_{\infty})}$ be a pseudo-limit of this sequence. In any case, for each $S$, $T$ and $\mathbf{i}$, we set:
$$\lambda^{(\mu_\infty)}_{S,T,\mathbf{i}} := \sum_{j\in T} (-1)^{t(j)+1}i_j\rho_{S,T\setminus \{j\},\mathbf{i}}^{(\mu_\infty)},$$
so that equation \eqref{eq2} at rank $\mu_{\infty}$ is automatically satisfied by definition.

By construction and definition of pseudo-limits, we have:
$$\forall \mu< \mu_{\infty},\; v(\theta_{\mu_\infty}- \theta_\mu)\in \left\{ \frac{\gamma_\mu}{p},+\infty\right\} \;\;\; \text{and} \;\;\;  v(
\rho_{S,U,\mathbf{i}}^{(\mu_\infty)}-\rho_{S,U,\mathbf{i}}^{(\mu)}) \in \left\{ \frac{\gamma_\mu-\gamma}{p},+\infty\right\},$$
which settles equations \eqref{eq4} and \eqref{eq5} up to rank $\mu_{\infty}$. This also implies that, for $\mu< \mu_{\infty}$:
\begin{align*}
    (1+\pi^{\gamma} x_{S,T})^{-1}&\theta_{\mu_\infty}^p \left(1+\pi^{\gamma}\sum_{\substack{0\leq \mathbf{i}\leq p-1\\ \mathbf{i}\neq\mathbf{0}}}(\lambda^{(\mu_\infty)}_{S,T,\mathbf{i}})^p  \alpha^\mathbf{i}\right)\\
    &\equiv (1+\pi^{\gamma} x_{S,T})^{-1}\theta_{\mu}^p \left(1+\pi^{\gamma}\sum_{\substack{0\leq \mathbf{i}\leq p-1\\ \mathbf{i}\neq\mathbf{0}}}(\lambda^{(\mu)}_{S,T,\mathbf{i}})^p  \alpha^\mathbf{i}\right) \pmod {\pi^{\gamma_\mu}}\\
    &\equiv 1 \pmod {\pi^{\gamma_\mu}}.
\end{align*}
Hence, if we set: $$\gamma_{\mu_\infty}:= v\left( (1+\pi^{\gamma} x_{S,T})^{-1}\theta_{\mu_\infty}^p \left(1+\pi^{\gamma}\sum_{\substack{0\leq \mathbf{i}\leq p-1\\ \mathbf{i}\neq\mathbf{0}}}(\lambda^{(\mu_\infty)}_{S,T,\mathbf{i}})^p  \alpha^\mathbf{i}\right)\right),$$
we deduce that $\gamma_{\mu_\infty}>\gamma_\mu$ for every $\mu< \mu_{\infty}$ and that there exists $x_{S,T}^{(\mu_\infty)}\in \mathcal{O}_K$ such that: $$1+\pi^{\gamma} x_{S,T} = \theta_{\mu_\infty}^p \left(1+\pi^{\gamma}\sum_{\substack{0\leq \mathbf{i}\leq p-1\\ \mathbf{i}\neq\mathbf{0}}}(\lambda^{(\mu_\infty)}_{S,T,\mathbf{i}})^p  \alpha^\mathbf{i}\right) (1 + \pi^{\gamma_{\mu_\infty}} x_{S,T}^{(\mu_\infty)}),$$
thus establishing assertions \eqref{eq1} and \eqref{eq3} at rank $\mu_\infty$.

\paragraph{Termination of the inductive process.} Since $\Gamma$ is a set, our inductive construction of the sequences $(\theta_\mu)_\mu$, $(\gamma_\mu)_\mu$, $(x_{S,T}^{(\mu)})_\mu$, $(\lambda^{(\mu)}_{S,T,\mathbf{i}})_\mu$, and  $(\rho_{S,U,\mathbf{i}}^{(\mu)})_\mu$ has to stop. This happens when $\gamma_{\mu_\infty}\geq \gamma_\infty$ for some ordinal $\mu_\infty$. We then have, for each $S$ and $T$:
$$v\left(1+\pi^{\gamma} x_{S,T} - \theta_{\mu_\infty}^p \left(1+\pi^{\gamma}\sum_{\substack{0\leq \mathbf{i}\leq p-1\\ \mathbf{i}\neq\mathbf{0}}}(\lambda^{(\mu_\infty)}_{S,T,\mathbf{i}})^p  \alpha^\mathbf{i}\right)\right)\geq \gamma_\infty.$$
Lemma \ref{lem annulation} then shows that $a\in u^{\geq \gamma_\infty}_q(K)$, as wished.
\end{proof}

\section{Proof of Theorem \ref{thm cd}}\label{sec cd}

In this last section we deduce Theorem \ref{thm cd} from Theorem \ref{thm isom Kato} and Theorem \ref{th invisible terms}. We keep the notation and assumptions of Section \ref{subsec setting}. In particular, as we mention in Remark \ref{rem cd}, we assume that $\zeta_p\in K$ and that $\Gamma$ is dense in $\R$. Moreover, recall that $g=\dim_{\F_p}(\Gamma/p\Gamma)$ and $p^r=[k:k^p]$, so that $\dim_p(k)\in\{r,r+1\}$ in this case.

\begin{proposition}\label{prop cd lower}
We have $\dim_p(K)\geq g+\dim_p(k)$.
\end{proposition}

\begin{proof}
By the Bloch-Kato conjecture it suffices to exhibit a nonzero element in $\kk_{g+\dim_p(k)}(K)$.

If $g=\infty$ or $r=\infty$, then for every $N\geq0$ one finds an $\F_p$-free family $\bar\gamma_1,\dots,\bar\gamma_a$ in $\Gamma/p\Gamma$ and a $p$-free family $\bar x,\bar y_1,\dots,\bar y_{b-1}$ in $k$ with $a+b=N$, so that Proposition \ref{lem triv symb} gives $\kk_N(K)\neq0$ and hence $\dim_p(K)=\infty$. We may therefore assume that $g$ and $r$ are both finite.

Assume first that $\dim_p(k)=r$. Fix an $\F_p$-basis $\bar\gamma_1,\dots,\bar\gamma_g$ of $\Gamma/p\Gamma$ and a $p$-basis $\bar\alpha_1,\dots,\bar\alpha_r$ of $k$. By the Bloch-Kato-Gabber theorem \cite[Thm.~9.5.2]{GS}, the group $\kk_r(k)$ is isomorphic to $\nu(r)(k)\subseteq\Omega_k^r$ and the symbol $\{\bar\alpha_1,\dots,\bar\alpha_r\}$ corresponds to $\frac{d\bar\alpha_1}{\bar\alpha_1}\wedge\cdots\wedge\frac{d\bar\alpha_r}{\bar\alpha_r}$, which is nonzero since $\bar\alpha_1,\dots,\bar\alpha_r$ is a $p$-basis; hence $\kk_r(k)\neq0$. Setting $q:=g+r$, the summand $\Lambda^{g}(\Gamma/p\Gamma)\otimes\kk_{r}(k)$ of the source of $\rho_q^0$ is then nonzero, and item (i) of Theorem \ref{thm isom Kato} shows that $\kk_q(K)/u_q^{>0}(K)\neq0$. In particular $\kk_{g+r}(K)\neq0$ and $\dim_p(K)\geq g+r=g+\dim_p(k)$.

Assume now that $\dim_p(k)=r+1$, that is, $H_p^{r+1}(k)\neq0$. Set $q:=g+r+1$ and note that $g\leq q-1$. The summand $\Lambda^{g}(\Gamma/p\Gamma)\otimes H_p^{r+1}(k)$ of the source of $\rho_q^e$ is then nonzero, and item (iv) of Theorem \ref{thm isom Kato} shows that $u_q^{e}(K)/u_q^{>e}(K)\neq0$. In particular $\kk_{g+r+1}(K)\neq0$ and $\dim_p(K)\geq g+r+1=g+\dim_p(k)$.
\end{proof}

The point of the following lemma is that all the groups appearing in Theorem \ref{thm isom Kato} vanish above the expected range.

\begin{lemma}\label{lem vanishing graded}
Assume $g,r<\infty$ and let $q>g+\dim_p(k)$. Then $\kk_q(K)=u_q^{>0}(K)$ and $u_q^{\gamma}(K)=u_q^{>\gamma}(K)$ for every $\gamma\in\Gamma_{\geq0}$.
\end{lemma}

\begin{proof}
We have $\dim_p(k)\in\{r,r+1\}$ and $H_p^m(k)=0$ for every $m>\dim_p(k)$. Recall also that $\Lambda^j(\Gamma/p\Gamma)=0$ for $j>g$, that $\Omega_k^m=0$ for $m>r$, that $\Omega_k^m/\ker(d)=0$ for $m\geq r$ and that $\kk_m(k)=0$ for $m>r$. We examine the four items of Theorem \ref{thm isom Kato} separately.
\begin{enumerate}[(i)]
\item For $0\leq j\leq q$, either $j>g$ and $\Lambda^j(\Gamma/p\Gamma)=0$, or $j\leq g$ and then $q-j>r$, so that $\kk_{q-j}(k)=0$. Hence the source of $\rho_q^0$ vanishes and $\kk_q(K)=u_q^{>0}(K)$.
\item For $0\leq j\leq q-1$, either $j>g-1$ and $\Lambda^j\Delta=0$ (as $\dim_{\F_p}\Delta=g-1$), or $j\leq g-1$ and then $q-j-1\geq q-g>r$, so that $\Omega_k^{q-j-1}=0$. Hence the source of $\rho_q^\gamma$ vanishes for $0<\gamma<e$ with $\gamma\notin p\Gamma$.
\item For $0\leq j\leq q-1$, either $j>g$ and $\Lambda^j(\Gamma/p\Gamma)=0$, or $j\leq g$ and then $q-j-1\geq q-g-1\geq r$, so that $\Omega_k^{q-j-1}/\ker(d)=0$. Hence the source of $\rho_q^\gamma$ vanishes for $0<\gamma<e$ with $\gamma\in p\Gamma$.
\item For $0\leq j\leq q-1$, either $j>g$ and $\Lambda^j(\Gamma/p\Gamma)=0$, or $j\leq g$ and then $q-j>\dim_p(k)$, so that $H_p^{q-j}(k)=0$. Hence the source of $\rho_q^e$ vanishes.
\end{enumerate}
Since the $\rho_q^\gamma$ are isomorphisms, we get $u_q^{\gamma}(K)=u_q^{>\gamma}(K)$ for every $\gamma\in\Gamma$ with $0\leq\gamma\leq e$, while for $\gamma>e$ both groups are zero by Lemma \ref{lem gamma sup e}.
\end{proof}

\begin{proposition}\label{prop cd upper}
Assume that $(K,v)$ is maximally complete and that $g,r<\infty$. Then $\kk_q(K)=0$ for every $q>g+\dim_p(k)$.
\end{proposition}

\begin{proof}
Fix such a $q$ and take $a\in\kk_q(K)$. By Lemma \ref{lem vanishing graded} we have $a\in u_q^{>0}(K)$, and since $u_q^{>\gamma}(K)=\bigcup_{\gamma'>\gamma}u_q^{\gamma'}(K)$ for every $\gamma\in\Gamma_{\geq0}$, the set $S:=\{\gamma\in\Gamma_{>0}: a\in u_q^{\gamma}(K)\}$ is a nonempty interval in $\Gamma_{>0}$. Let $\gamma_\infty:=\sup S\in\R_{>0}\cup\{+\infty\}$.

Suppose first that $\gamma_\infty>e$. Then $a\in u_q^{\gamma}(K)$ for some $\gamma>e$, hence $a=0$ by Lemma \ref{lem gamma sup e}, as wished.

Suppose now that $\gamma_\infty\leq e$. By definition of $\gamma_\infty$ we have $a\in u_q^{\gamma}(K)$ for every $\gamma<\gamma_\infty$, that is, $a\in u_q^{<\gamma_\infty}(K)$. Since $(K,v)$ is maximally complete, Theorem \ref{th invisible terms} gives $a\in u_q^{\geq\gamma_\infty}(K)$, so that $a\in u_q^{\gamma_0}(K)$ for some $\gamma_0\in\Gamma$ with $\gamma_0\geq\gamma_\infty$. By definition of $\gamma_\infty$ as a supremum we get $\gamma_0=\gamma_\infty\in\Gamma$ and $a\in u_q^{\gamma_\infty}(K)$. But Lemma \ref{lem vanishing graded} gives $u_q^{\gamma_\infty}(K)=u_q^{>\gamma_\infty}(K)$, so $a\in u_q^{\gamma}(K)$ for some $\gamma>\gamma_\infty$, contradicting the definition of $\gamma_\infty$.
\end{proof}

\begin{proof}[Proof of Theorem \ref{thm cd}]
The inequality $\dim_p(K)\geq g+\dim_p(k)$ is Proposition \ref{prop cd lower}. Assume now that $(K,v)$ is maximally complete and set $N:=g+\dim_p(k)$; we may assume $N<\infty$, so that $g$ and $r$ are finite. By \cite[I, \S4, Prop.~21']{SerreCohGal} it suffices to prove that $H^{N+1}(L,\Z/p\Z)=0$ for every finite separable extension $L/K$.

Let $L/K$ be such an extension, let $\Gamma_L$ be its value group and $l$ its residue field. Then $(L,v)$ is again Henselian of rank $1$ with value group dense in $\R$, contains $\zeta_p$, and is again maximal (hence maximally complete) by \cite[Thm.~2.14]{BlaszczokKuhlmann}. Moreover $[l:l^p]=[k:k^p]=p^r$ and $\dim_p(l)=\dim_p(k)$, while $\dim_{\F_p}(\Gamma_L/p\Gamma_L)=\dim_{\F_p}(\Gamma/p\Gamma)=g$ since $\Gamma$ has finite index in $\Gamma_L$, see Lemma \ref{lem gps} here below.

Applying Proposition \ref{prop cd upper} to $(L,v)$ we get $\kk_{N+1}(L)=0$, hence $H^{N+1}(L,\Z/p\Z)\cong\kk_{N+1}(L)=0$ by the Bloch-Kato conjecture. Therefore $\dim_p(K)\leq N$, which together with Proposition \ref{prop cd lower} gives the desired equality.
\end{proof}

Since we were unable to find a reference for the equality $\dim_{\F_p}(\Gamma_L/p\Gamma_L)=\dim_{\F_p}(\Gamma/p\Gamma)$ used here above, we state and prove the following elementary lemma.

\begin{lemma}\label{lem gps}
Let $\Gamma$ be a subgroup of $\R$ and let $\Delta\subseteq\Gamma$ be a subgroup of finite index. Then $\dim_{\F_p}(\Gamma/p\Gamma)=\dim_{\F_p}(\Delta/p\Delta)$.
\end{lemma}

\begin{proof}
We may clearly assume that $[\Gamma:\Delta]$ is a power of $p$. Moreover, since $\R$ is torsion-free, the multiplication morphism $\Gamma\to n\Gamma:\gamma\mapsto n\gamma$ is an isomorphism for any $n\in\Z$. In particular, $\Gamma/p\Gamma\simeq p^a\Gamma/p^{a+1}\Gamma$ for any $a\in\Z$.

If $\dim_{\F_p}(\Gamma/p\Gamma)$ is infinite, then clearly so is $\dim_{\F_p}(\Delta/p\Delta)$. Assume then that $g=\dim_{\F_p}(\Gamma/p\Gamma)$ is finite. Since the index is a power of $p$, there exists $a\in\N$ such that $p^a\Gamma\subseteq \Delta$ and hence $p^{a+1}\Gamma\subseteq p\Delta$. In particular, $\Delta$ and $p^{a}\Gamma$ satisfy the hypotheses of the lemma as well. Then if $\dim_{\F_p}(\Delta/p\Delta)$ is infinite, so is $g=\dim_{\F_p}(\Gamma/p\Gamma)=\dim_{\F_p}(p^a\Gamma/p^{a+1}\Gamma)$, which is a contradiction. We may thus assume that both $p$-ranks are finite.

Write now $\Gamma/p^{a+1}\Gamma\simeq\bigoplus_{s\in S}\Z/p^{a_s}\Z$ for some set $S$ and $1\leq a_s\leq a+1$. Then clearly $\Gamma/p\Gamma\simeq \bigoplus_{s\in S}\Z/p\Z$,
so that we may assume that $|S|=g$. By the theory of finitely generated abelian groups, any subgroup (and hence any subquotient) of $\Gamma/p^{a+1}\Gamma$ has at most $g$ generators. Then this is the case for $\Delta/p\Delta$, which proves that $\dim_{\F_p}(\Gamma/p\Gamma)\geq \dim_{\F_p}(\Delta/p\Delta)$. Then, by the same reasoning, we get
\[\dim_{\F_p}(\Gamma/p\Gamma)\geq \dim_{\F_p}(\Delta/p\Delta)\geq \dim_{\F_p}(p^a\Gamma/p^{a+1}\Gamma)=\dim_{\F_p}(\Gamma/p\Gamma),\]
which yields the equality.
\end{proof}

\end{document}